\documentclass[11pt]{article}

\usepackage{amscd}      
\usepackage{amssymb}
\usepackage[intlimits]{amsmath}
\usepackage{xypic}      
\LaTeXdiagrams          
\usepackage[all,v2]{xy}
\xyoption{2cell}
\UseAllTwocells
\xyoption{frame}
\CompileMatrices

\usepackage{theorem}

\newtheorem{prop}{Proposition}[section]
\newtheorem{lem}[prop]{Lemma}

\newtheorem{cor}[prop]{Corollary}
\newtheorem{them}[prop]{Theorem}

\theorembodyfont{\upshape}

\newtheorem{defn}[prop]{Definition}

\newtheorem{rmk}[prop]{Remark}

\newtheorem{numex}[prop]{Example}

\newenvironment{pf}{\begin{trivlist}\item[]{\sc Proof.}}%
            {\nolinebreak $\Box$ \end{trivlist}}
\newenvironment{pff}{\begin{trivlist}\item[]{\sc Proof of Theorem \ref{thm:vb smooth}.}}%
            {\nolinebreak $\Box$ \end{trivlist}}

\newcommand{\noprint}[1]{}

\newcommand{\equal}{=}

\renewcommand{\tilde}{\widetilde}

\newcommand{\toto}{\rightrightarrows}

\newcommand{\upst}{^{\ast}}

\newcommand{\com}{{\scriptscriptstyle\bullet}}

\newcommand{\upcom}{^{\scriptscriptstyle\bullet}}

\newcommand{\XX}{{\mathfrak X}}

\newcommand{\RR}{{\mathfrak R}}
\newcommand{\HH}{{\Bbb H}}

\newcommand{\zz}{{\mathbb Z}}
\newcommand{\hh}{{\mathbb H}}

\newcommand{\nn}{{\mathbb N}}

\newcommand{\cc}{{\mathbb C}}
\newcommand{\rr}{{\mathbb R}}
\newcommand{\rrr}{{\cal R}}

\newcommand{\aA}{{\cal A}}
\newcommand{\bB}{{\cal B}}
\newcommand{\dD}{{\cal D}}
\newcommand{\cC}{{\cal C}}
\newcommand{\eE}{{\cal E}}

\newcommand{\fF}{{\cal F}}
\newcommand{\kK}{{\cal K}}
\newcommand{\lL}{{\cal L}}

\newcommand{\hH}{{\cal H}}
\newcommand{\sS}{{\cal S}}

\newcommand{\del}{\partial}

\newcommand{\pr}{\mathop{\rm pr}\nolimits}

\newcommand{\id}{\mathop{\rm id}\nolimits}

\newcommand{\smalcirc}{\mbox{\tiny{$\circ $}}}

\newcommand{\Gam}{{\Bbb E}}
\newcommand{\ldiag}[1]%
       {\makebox[0cm]{${\scriptstyle#1}\downarrow\phantom{\scriptstyle#1}$}}
\newcommand{\ldiagup}[1]%
       {\makebox[0cm]{${\scriptstyle#1}\uparrow\phantom{\scriptstyle#1}$}}
\newcommand{\rdiag}[1]%
       {\makebox[0cm]{$\phantom{\scriptstyle#1}\downarrow{\scriptstyle#1}$}}
\newcommand{\sediagr}[1]%
       {\makebox[0cm]{$\phantom{\scriptstyle#1}\searrow{\scriptstyle#1}$}}
\newcommand{\nediagr}[1]%
       {\makebox[0cm]{$\phantom{\scriptstyle#1}\nearrow{\scriptstyle#1}$}}
\newcommand{\rdiagup}[1]%
       {\makebox[0cm]{$\phantom{\scriptstyle#1}\uparrow{\scriptstyle#1}$}}
\newcommand{\swdiag}[1]%
       {\makebox[0cm]{$\phantom{\scriptstyle#1}\swarrow{\scriptstyle#1}$}}
\newcommand{\sediag}[1]%
       {\makebox[0cm]{${\scriptstyle#1}\searrow\phantom{\scriptstyle#1}$}}
\newcommand{\nediag}[1]%
       {\makebox[0cm]{${\scriptstyle#1}\nearrow\phantom{\scriptstyle#1}$}}

\newcommand{\longiso}{\stackrel{\textstyle\sim}{\longrightarrow}}

\newcommand{\doublearrowstack}[2]%
 {{{{\scriptstyle#1}\atop{\textstyle\longrightarrow}}\atop{{\textstyle\longright
arrow}\atop{\scriptstyle#2}}}}
\newcommand{\rightleftarrowstack}[2]%
 {{{{\scriptstyle#1}\atop{\textstyle\longrightarrow}}\atop{{\textstyle\longlefta
rrow}\atop{\scriptstyle#2}}}}
\newcommand{\leftrightarrowstack}[2]%
 {{{{\scriptstyle#1}\atop{\textstyle\longleftarrow}}\atop{{\textstyle\longrighta
rrow}\atop{\scriptstyle#2}}}}

\newcommand{\overtoparrow}%
{\makebox[0cm]{\beginpicture
\setcoordinatesystem units <.8cm,.4cm> point at 0 0
\setplotarea x from -3 to 3, y from 0 to 1
\setquadratic
\plot -3 0 0 1 3 0 /
\put{\vector(3,-1){0}}[Bl] at 3 0
\endpicture}}

\newcommand{\underbottomarrow}%
{\makebox[0cm]{\beginpicture
\setcoordinatesystem units <.8cm,.4cm> point at 0 0
\setplotarea x from -3 to 3, y from 0 to 1
\setquadratic
\plot -3 1 0 0 3 1 /
\put{\vector(3,1){0}}[Bl] at 3 1
\endpicture}}

\newcommand{\ses}[5]%
{0\longrightarrow#1\stackrel{#2}{ \longrightarrow}#3\stackrel{#4}{
\longrightarrow}#5\longrightarrow0}

\newcommand{\dt}[6]%
{#1\stackrel{#2}{\longrightarrow}#3 \stackrel{#4}{\longrightarrow}#5
\stackrel{#6}{\longrightarrow} #1[1]}

\newcommand{\cat}[1]%
{(\mbox{\rm #1})}

\newcommand{\gm}{\Gamma}
\newcommand{\lon }{\longrightarrow }

\newcommand{\tp}{\tilde{p}}

\newcommand{\be }{\begin{eqnarray*}}
\newcommand{\ee }{\end{eqnarray*}}

\def\gpd{\,\lower1pt\hbox{$\longrightarrow$}\hskip-.24in\raise2pt
             \hbox{$\longrightarrow$}\,}

\title{Twisted $K$-theory of Differentiable Stacks} 
\author{ Jean-Louis Tu\\ 
University Pierre et Marie Curie\\
 175 rue du Chevaleret\\
 75013 Paris, France\\
 {\sf email: tu@math.jussieu.fr} 
 \\\\
Ping Xu \thanks{ Research partially supported by NSF
       grants DMS00-72171 and DMS03-06665. }\\
        Department of Mathematics\\
         Pennsylvania State University \\
         University Park, PA 16802, USA\\
{\sf email: ping@math.psu.edu }
\\\\
AND\\
Camille Laurent-Gengoux\\
D\'epartement de math\'ematiques\\
Universit\'e de Poitiers\\
86962 Futuroscope  Chasseneuil, France\\
{\sf email:  laurent@math.univ-poitiers.fr}}

\begin{document}
\sloppy
\maketitle
\begin{abstract}
In this paper, we develop twisted $K$-theory for stacks, where the twisted class
is given by an $S^1$-gerbe over the stack. General properties, including the
Mayer-Vietoris property, Bott periodicity, and the product structure $K^i_\alpha
\otimes K^j_\beta \to K^{i+j}_{\alpha +\beta}$ are derived. Our approach
provides a uniform framework for studying various twisted $K$-theories including
the usual twisted $K$-theory of topological spaces, twisted equivariant
$K$-theory, and the twisted $K$-theory of orbifolds. We also present a Fredholm
picture, and discuss the conditions under which twisted $K$-groups can be
expressed by so-called ``twisted vector bundles". 

Our approach is to work on presentations of stacks, namely \emph{groupoids}, and
relies heavily on the machinery of $K$-theory ($KK$-theory) of $C^*$-algebras.

\par\bigskip
\begin{center}
{\bf $K$-th\'eorie tordue des champs diff\'erentiables}

{\bf R\'esum\'e}
\end{center}
Dans cet article, nous d\'eveloppons la $K$-th\'eorie tordue pour les
champs diff\'erentiables, o\`u la torsion s'effectue
par une $S^1$-gerbe sur le champ en question.
Nous en \'etablissons les propri\'et\'es g\'en\'erales telles que
les suites exactes de Mayer-Vietoris, la p\'eriodicit\'e de Bott,
et le produit $K^i_\alpha \otimes K^j_\beta \to K^{i+j}_{\alpha +\beta}$.
Notre approche fournit un cadre g\'en\'eral permettant d'\'etudier diverses
$K$-th\'eories tordues, en particulier la $K$-th\'eorie tordue usuelle
des espaces topologiques, la $K$-th\'eorie tordue \'equivariante,
et la $K$-th\'eorie tordue des orbifolds. Nous donnons \'egalement
une d\'efinition \'equivalente utilisant des op\'erateurs de Fredholm,
et nous discutons les conditions sous lesquelles les groupes de
$K$-th\'eorie tordue peuvent \^etre r\'ealis\'es \`a partir de ``fibr\'es
vectoriels tordus''.

Notre approche consiste \`a travailler sur les r\'ealisations
concr\`etes des champs, \`a savoir les \emph{groupo\"{\i}des},
et s'appuie de fa\c{c}on importante sur les techniques
de $K$-th\'eorie ($KK$-th\'eorie) des $C^*$-alg\`ebres.
\end{abstract}

{\small \tableofcontents}

\section{Introduction}
Recently, motivated by $D$-branes in 
string theory, there has been a great deal of
interest in the study of twisted $K$-theory \cite{BM, MM, M, Witten}. The
$K$-theory of a topological space $M$ twisted by a torsion class in $H^3 (M,
\zz)$ was first studied by Donovan-Karoubi \cite{don-kar70} in the early 1970s,
and, in the 1980s, Rosenberg \cite{rosenberg1} introduced $K$-theory
twisted by a general element of $H^3 (M, \zz)$. More recently, twisted
$K$-theory has enjoyed renewed vigor due to the discovery of its close
connection with string theory \cite{Witten, Witten1}. See also \cite{Atiyah,
Segal, BCMMS} and references therein.

A very natural problem which arises is the development of other types
of twisted $K$-theory. 
In particular, twisted equivariant $K$-theory and twisted 
$K$-theory for orbifolds should be developed. 
Indeed, various definitions of  such theories have been offered.
For instance, Adem--Ruan introduced a version of twisted
$K$-theory of an orbifold by  a discrete torsion element  \cite{AR}.
Others, for example \cite{LU, MS}, offer various related 
(but unsupported) definitions.
We also remark that
Freed--Hopkins--Teleman announced \cite{F} the amazing result that the
twisted equivariant $K$-theory of a compact Lie group is related to the
Verlinde algebra.


It is important that twisted $K$-theory is a cohomology theory and, in
particular, satisfies the Mayer-Vietoris property. One also expects that, like
ordinary $K$-theory, it should satisfy Bott periodicity. 
The purpose of this paper is to develop a twisted $K$-theory for stacks,
the idea being that this is general enough
to include all the above cases, including twisted equivariant $K$-theory and
twisted $K$-theory of orbifolds.
As far as we know, except for the special case of manifolds, there has been no twisted
$K$-theory for general stacks for which all such properties have been 
established (as far as we know, this is the case even for twisted
equivariant $K$-theory).

Rather than working directly with stacks, we will work on presentations
of  stacks, namely groupoids. 
Indeed there is a dictionary in which 
a stack corresponds to a {\em Morita equivalence class} of groupoids \cite{BX, BX1}. 
In this paper, we will deal with
differentiable stacks which are more relevant to string theory.
They correspond to Lie groupoids.

An advantage of working with Lie groupoids, for a differential
geometer, is that one can still do differential geometry even though
the spaces they represent do not usually allow such a possibility. 

The notion of a \emph{groupoid} is a standard generalization of the
concepts of \emph{spaces} and \emph{groups}. 
In the theory of groupoids, spaces and groups are treated on equal footing.
Simplifying somewhat, one could say that a groupoid is a mixture of a space and
a group; it has space-like and group-like properties that interact in
a delicate way. In a certain sense, groupoids provide a uniform
framework for many different geometric objects. 
For instance, when a Lie group acts on a manifold
properly, the corresponding equivariant cohomology theories, including
$K$-theory, can be treated using the transformation groupoid $G\times
M\toto M$. 
On the other hand, an orbifold can be represented by an \'etale groupoid
\cite{MP, Moer}.

The problem of computing the $K$-theory of groupoids has been studied by many
authors. For instance, given a locally compact groupoid $\Gamma$, the
Baum-Connes map $\mu_r\colon K_*^{\mathrm{top}}(\Gamma)\to K_*(C^*_r(\Gamma))$
can be used to study the $K$-theory groups of $C^*_r(\Gamma)$. The above map
generalizes both the assembly map for groups \cite{bch94} and the coarse
assembly map \cite{sty02}. Of course, many techniques used to study the
Baum-Connes conjecture for groups \cite{hig98} can be extended to groupoids such
as foliation groupoids \cite{tu99,tu99a}. However, recent counterexamples
\cite{hls02} show that other ways of attacking the problem need to be
discovered. Applications of the $K$-theory of groupoids include: tilings and gap
labeling (see for instance \cite{kel95}), index theorems, and pseudodifferential
calculi \cite{mon01, pat01}.

By twisted $K$-theory, in this paper we mean $K$-theory twisted by an $S^1$-gerbe. All
$S^1$-gerbes over a groupoid $\gm$ (or more precisely a stack $\XX_\gm$) form
an Abelian group which can be identified with $H^2 (\XX_\gm, {\mathcal{S}}^1 )$
\cite{BX, BX1}. Unlike the manifold case, this is not always isomorphic to the
third integer cohomology group $H^3 (\XX_\gm, \zz )$. Indeed, this fails to be
the case even when $\gm$ is a non-compact group. Not enough attention seems to
have been paid to this in the literature. However, for a proper Lie groupoid
$\gm$,
these two groups are always isomorphic, and it therefore makes sense to talk
about its $K$-theory twisted by an integer class in $H^3 (\XX_\gm, \zz )$. In
particular, when a Lie group $G$ acts on a smooth manifold $M$ properly,
 one can define
the equivariant $K$-theory twisted by an element in $H^3 _G (M, \zz)$. The same
situation applies to orbifolds since their corresponding groupoids are always
proper.

Our approach in developing twisted $K$-theory is to utilize operator algebras,
where many sophisticated $K$-theoretic techniques have been developed. An
$S^1$-central extension $S^1\to R\to \Gamma \toto M$ of groupoids gives rise to
an $S^1$-gerbe $\RR$ over the differentiable stack $\XX_\Gamma$ associated to
the groupoid $\gm$ \cite{BX}, and Morita equivalent $S^1$-central extensions
correspond to isomorphic gerbes. Therefore, given a Lie groupoid $\gm \toto M$,
one may identify an $S^1$-gerbe over the stack $\XX_\gm$ as a Morita equivalence
class of $S^1$-central extensions $S^1\to R'\to \Gamma' \toto M'$, where
$\Gamma' \toto M'$ is a Lie groupoid Morita equivalent to $\gm \toto M$. Given
an $S^1$-central extension of Lie groupoids $S^1\to R\to \Gamma\toto M$, its
associated complex line bundle $L\equal R\times_{S^1}\cc$ can be considered as a
Fell bundle of $C^*$-algebras over the groupoid $\Gamma\toto M$. Therefore, from
this one can construct the reduced $C^*$-algebra $C^*_r(\Gamma,R)$. The
$K$-groups are simply defined to be the $K$-groups of this $C^*$-algebra. 

This
definition yields several advantages. First, since it is standard that Morita
equivalent central extensions yield Morita equivalent $C^*$-algebras, the
$K$-groups indeed only depend on the stack and the $S^1$-gerbe, instead of on
any particular groupoid $S^1$-central extension. Such a viewpoint is quite interesting
already, even when dealing with untwisted $K$-theory. For instance, some
classical results of Segal on equivariant $K$-theory \cite{Segal} may be
reinterpreted as a consequence of the fact that equivariant $K$-theory only
depends on the stack $M/G$, i.e. the Morita equivalence class of the
transformation groupoid $G \times M \toto M$. Secondly, important properties of
$K$-theory, such as the Mayer-Vietoris property and Bott periodicity, are
immediate consequences of this definition.

A drawback of this definition, however, is that it is too abstract and
algebraic. Our second goal in this paper
is to connect it with the usual topological approach of
$K$-theory in terms of Fredholm bundles \cite{Atiyah1, Segal1}. 
As in the
manifold case, an $S^1$-central extension of a groupoid naturally gives rise
to a canonical  principal $PU(\HH )$-bundle over the groupoid, which in turn induces
associated Fredholm bundles over the groupoid. We show that the $K$-groups can
be interpreted as homotopy classes of invariant sections of these Fredholm
bundles (assuming a certain appropriate continuity). This picture fits with
the usual definition of twisted $K$-theory \cite{Atiyah} when the groupoid
reduces to a space.

Geometrically, it is always desirable to describe $K$-groups in terms of vector
bundles. For twisted $K$-groups, a natural candidate is to use twisted vector
bundles over the groupoid. This is a natural generalization, in the context of
groupoids, of projective representations of a group. More precisely, given an
$S^1$-central extension of Lie groupoids
 $S^1\to R\stackrel{\pi}{\to}\Gamma\toto M$, a
twisted vector bundle is a vector bundle $E$ over the groupoid $R$ where $\ker
\pi \cong M\times S^{1}$ acts on $E$ by scalar multiplication. When $\gm$ is a
groupoid Morita equivalent to a manifold, they correspond to the so-called
bundle gerbe modules in \cite{BCMMS}. However, note that twisted vector bundles
do not always exist. In fact, a necessary condition for their existence is that
the twisted class $\alpha \in H^2 (\gm\upcom, {\mathcal{S}}^1 )$ must be a   
torsion. Another main theme of this paper is to explore the conditions under
which the twisted $K_0$-group is isomorphic to the Grothendieck group of twisted
vector bundles. 

As is already the case for manifolds, twisted $K$-groups no longer
admit a ring structure \cite{Atiyah}. 
It is expected, however, that there exists a bilinear
product $K_\alpha^i\otimes K_\beta^j \to K_{\alpha+\beta}^{i+j}$. For twisted
vector bundles, such a product is obvious and corresponds to the tensor product
of vector bundles. However, in general,
 twisted $K$-groups cannot be expressed by
twisted vector bundles as discussed above. The main difficulty in constructing
such a product using the Fredholm picture of twisted $K$-theory is obtaining a
Fredholm operator $T$ out of two Fredholm operators $T_1$ and $T_2$. This is
very similar to the situation of the Kasparov product in $KK$-theory where a
non-constructive method must be used. Motivated by $KK$-theory, our approach is
to develop a generalized version of Le Gall's groupoid equivariant $KK$-theory
and interpret our twisted $K$-groups as such $KK$-groups, which allows us to
obtain such a product.


The paper is organized as follows. Section 2 is devoted to the basic theory of
$S^1$-gerbes over stacks in terms of the groupoid picture; related cohomology
theory and characteristic classes are reviewed briefly. In Section 3, we
introduce the definition of twisted $K$-groups and outline some basic
properties. Section 4 is devoted to the study of the $K$-groups of the
$C^*$-algebra of a Fell bundle over a proper groupoid, which includes our
$C^*$-algebras of groupoid  $S^1$-central extensions
 as a special case. In particular,
we give the Fredholm picture of the $K$-groups. In Section 5, we investigate the
conditions under which the twisted $K_0$-group can be expressed in terms of
twisted vector bundles. In Section 6, we discuss the construction of the
$K$-group product as outlined above. In the Appendix, we review some basic
material concerning Fell bundles over groupoids which we use frequently in the
paper.

We would like to point out that there are many interesting and
important questions that we are not able to address in this paper. One of them
is the study of the Chern character in twisted $K$-theory, in which Connes'
noncommutative differential geometry \cite{con} will play a prominent role due
to the nature of our algebraic definition. This subject will be discussed in a
separate paper.

Finally, we note that after our paper was submitted a paper by Atiyah
and Segal \cite{Ati-Seg} appeared, in which twisted equivariant
$K$-theory (for a compact group acting on a space) was introduced
independently using a different method.  It is not hard to check that 
at least in our case of interest, that of a compact Lie group acting on a
manifold, our twisted $K$-theory coincides with theirs (using the 
remark in Appendix A1 of \cite{Ati-Seg} that, in the metrizable case, the
compact-open topology is the same as the strong topology).

\paragraph{Notations}

Finally, we list the notation used throughout the paper. $\Gamma$ will denote a
groupoid (all groupoids considered are Hausdorff, locally compact, and second
countable). We denote by $s$ and $t$ the source and target maps
of $\Gamma$, respectively. $\Gamma^{(0)}$ will denote the unit space of
$\Gamma$, and $\Gamma^{(n)}$ will denote the set of strings of length $n$
\[ g_1\leftarrow g_2\leftarrow\cdots\leftarrow g_n, \]
i.e., the set of $n$-tuples
$(g_1,\ldots,g_n)\in \Gamma\times\cdots \times \Gamma$ such that $s(g_i)\equal
t(g_{i+1})$ for all $i\equal 1,\ldots, n-1$.

We will commonly use the expression ``Let $\Gamma\toto M$ be a Lie groupoid
\ldots'' to indicate that $\Gamma$ is a Lie groupoid and $\Gamma^{(0)}\equal M$.

For all $K$, $L\subset \Gamma^{(0)}$, we let $\Gamma_K\equal s^{-1}(K)$,
$\Gamma^L\equal t^{-1}(L)$, and $\Gamma_K^L\equal \Gamma_K\cap\Gamma^L$. If
$K\equal \{x\}$ and $L\equal \{y\}$, we will use the notation $\Gamma_x$,
$\Gamma^y$, and $\Gamma_x^y$, respectively.

If $Y$ is a space, then $Y\times Y$ will be endowed with the pair groupoid
structure: $(Y\times Y)^{(0)}\equal Y$, $s(y_1,y_2)\equal y_2$,
$t(y_1,y_2)\equal y_1$, and $(y_1,y_2)(y_2,y_3) \equal (y_1,y_3)$.

If $Y$ is a space and $f\colon Y\to \Gamma^{(0)}$ is a map, we denote by
$\Gamma[Y]$ the subgroupoid of $(Y\times Y)\times \Gamma$ consisting of
$\{(y_1,y_2,\gamma)|\; \gamma\in\Gamma_{f(y_2)}^{f(y_1)}\}$.  Then $\Gamma[Y]$
is called the pullback of the groupoid $\Gamma$ by $f$.

In particular, if ${\mathcal{U}}\equal (U_i)$ is an open cover of
$\Gamma^{(0)}$, then the pullback of $\Gamma$ by the canonical map $\coprod U_i
\to \Gamma^{(0)}$ is denoted either by $\Gamma[{\mathcal{U}}]$ or by
$\Gamma[U_i]$.

\par\medskip

If $\Gamma\toto M$ is a locally compact groupoid (resp., a Lie groupoid), a Haar
system for $\Gamma$ will usually be denoted by $\lambda\equal (\lambda^x)_{x\in
M}$, where $\lambda^x$ is a measure with support $\Gamma^x$ such that for all
$f\in C_c(\Gamma)$ (resp., $C_c^\infty(\Gamma)$), $x\mapsto \int_{g\in\Gamma^x}
f(g)\,\lambda^x(dg)$ is continuous (resp., smooth).

\par\bigskip

Let $\HH$ be the separable Hilbert space. We denote by ${\mathcal{K}}(\HH)$, or
even ${\mathcal{K}}$, the algebra of compact operators on $\HH$; we denote by
${\mathcal{L}}(\HH)$ the algebra of linear bounded operators.

For a $C^*$-algebra $A$, $M(A)$ denotes its multiplier algebra \cite[Section
3.12]{ped79}. Recall that $M(A)$ is a unital $C^*$-algebra containing $A$ as an
essential ideal, and, moreover, if a $C^*$-algebra $B$ also contains $A$ as an
essential ideal, then $A\subset B\subset M(A)$. For instance, if $X$ is a
locally compact space and $A\equal C_0(X)$, then $M(A)\equal C_b(X)$ is the
space of continuous bounded functions on $X$. On the other hand, if $A\equal
{\mathcal{K}}$, then $M(A)\equal {\mathcal{L}}(\HH)$.

For a Hilbert $C^*$-module ${\mathcal{E}}$ over $A$ (see \cite{wegge93}), we
denote by ${\mathcal{L}}({\mathcal{E}})$ the algebra of ($A$-linear bounded)
adjointable operators on ${\mathcal{E}}$. For all $\xi$, $\eta\in
{\mathcal{E}}$, let $T_{\xi,\eta}$ be the operator $T_{\xi,\eta}(\zeta)\equal
\xi\langle\eta,\zeta\rangle$. Then $T_{\xi,\eta}$ is called a rank-one operator.
The closed linear span of rank-one operators is called the algebra of compact
operators on ${\mathcal{E}}$ and will be denoted by
${\mathcal{K}}({\mathcal{E}})$; this is an ideal of
${\mathcal{L}}({\mathcal{E}})$.

We gather below the most frequently used notations and terminology:
\par\medskip
\begin{small}
\noindent
$\cC(\eE)$\dotfill
Equation (\ref{eqn:C(E)})

\par\noindent
$C^*_r(\Gamma)$ [reduced $C^*$-algebra of a groupoid]\dotfill
Reference~\cite{ren80}

\par\noindent
$C^*_r(\Gamma;E)$ [reduced $C^*$-algebra of a Fell bundle]\dotfill
Subsection~\ref{subsec:reduced}

\par\noindent
$C^*_r(\Gamma;R)$ [$C^*$-algebra of an $S^1$-central extension]\dotfill
Definition~\ref{def:CGR}

\par\noindent
$C^*_r(R)^{S^1}$\dotfill
Equation (\ref{eqn:CS1})

\par\noindent
$S^1$-Central extension\dotfill
Definition~\ref{def:central extension}

\par\noindent
$G$-bundle over a groupoid\dotfill
Definition~\ref{def:G bundle}

\par\noindent
$\Gamma$-space\dotfill
Definition~\ref{def:principal}

\par\noindent
Generalized homomorphism\dotfill
Definition~\ref{def:generalized morphism}

\par\noindent
$\eE^{sm}(\Gamma)$, $\mbox{Ext}^{sm}(\Gamma,S^1)$\dotfill
Subsection~\ref{subsec:central extensions}

\par\noindent
$\fF_\alpha^i$\dotfill
Theorem~\ref{thm:fredholm proper}

\par\noindent
Fell bundle over a groupoid $\Gamma$\dotfill
Definition~\ref{def:fell over groupoid}

\par\noindent
$K_\alpha^i(\Gamma\upcom)$ [twisted $K$-theory]\dotfill
Definition~\ref{def:twisted K theory}

\par\noindent
$\kK_\Gamma(\eE)$\dotfill
Equation (\ref{eqn:K(E)})

\par\noindent
$L^1(\Gamma;E)$, $L^2(\Gamma;E)$\dotfill
Subsection~\ref{subsec:reduced}

\par\noindent
Morita equivalent extension\dotfill
Definition~\ref{def:morita equivalent extension}

\par\noindent
Strictly trivial extension\dotfill
Proposition~\ref{prop:strictly trivial}

\par\noindent
$T^\Gamma$ [``average'' of $T$]\dotfill
Equation (\ref{eqn:TGamma})
\par\noindent
Trivial central extension\dotfill
Proposition~\ref{prop:trivial extensions}

\par\noindent
Twisted vector bundle\dotfill
Definition~\ref{def:twisted vb}
\end{small}

{\bf Acknowledgments.}
We  would like to thank several institutions
for their hospitality while work on this project was being done: 
Penn State University
(Laurent-Gengoux and Tu),  and RIMS, Ecole Polytechnique (Xu).
 We  also  wish to thank many people for  useful discussions and comments,
including  Paul Baum,  Kai Behrend, Jean-Luc Brylinski,
 Andre Henriques, Nigel Higson, Pierre-Yves Le Gall,  Jeff Raven,
Jean Renault, Jonathan Rosenberg, 
Jim Stasheff, especially Bin Zhang, who was participating in this
project at its earlier stage.

\section{$S^1$-gerbes and  central extensions of groupoids}

\subsection{Generalized homomorphisms}

In this subsection, we will review some basic facts
concerning generalized homomorphisms. Here groupoids
are assumed to be Lie groupoids although most of
the discussions can be easily adapted to general locally compact groupoids.
 Let us recall the definition
below \cite{hae84,hilsum-skandalis87,mrcun}.

\begin{defn}\label{def:generalized morphism}
A generalized groupoid homomorphism
from $\Gamma$
to $G$ is given by a manifold $Z$, two smooth maps
$$\Gamma^{(0)}\stackrel{\tau }{\leftarrow} Z
\stackrel{\sigma}{\rightarrow} G^{(0)},$$
a left action of $\Gamma$ with respect to $\tau$,
a right action of $G$ with respect to $\sigma$, such that the two
actions commute, and  $Z$ is a locally trivial $G$-principal
bundle over $\Gamma^{(0)}\stackrel{\tau}{\leftarrow} Z$.
\end{defn}

To explain the terminology,
if $f\colon \Gamma\to G$ is a strict homomorphism
(i.e. a smooth map satisfying $f(gh)\equal f(g)f(h)$) then
$Z_f\equal \Gamma^{(0)}\times_{f,G^{(0)},t} G$, with $\tau(x,g)\equal x$,
$\sigma(x,g)\equal s(g)$, and the actions 
$\gamma\cdot(x,g)\equal (t(\gamma),f(\gamma)g)$ and
$(x,g)\cdot g'\equal (x,gg')$, is a generalized homomorphism
from $\Gamma$ to $G$.

Generalized homomorphisms can be composed  just like the usual
groupoid homomorphisms.

\begin{prop}\label{prop:composition homomorphisms}
Let $Z$ and $Z'$ be generalized homomorphisms from
$\gm\toto \gm^{(0)} $ to $G\toto G^{(0)}$, and from
$G\toto G^{(0)}$ to $H\toto H^{(0)}$ respectively. Then
$$Z''\equal Z\times_G Z':\equal {(Z\times_{\sigma,G^{(0)},\tau'}Z')}
/_{(z,z')\sim(zg,g^{-1}z')}$$
is a generalized groupoid homomorphism
from $\gm\toto \gm^{(0)}   $ to $H\toto H^{(0)}$. Moreover, the composition
of generalized homomorphisms is associative, and thus there is a
category ${\mathcal{G}}$ whose objects are Lie groupoids and morphisms
are isomorphism classes of generalized
homomorphisms\footnote{Two generalized homomorphisms
$Z_1$ and $Z_2$ are isomorphic whenever they are $\Gamma$, $G$-equivariantly
diffeomorphic.}. There is a functor
$${\mathcal{G}}_s\to {\mathcal{G}}$$
where ${\mathcal{G}}_s$ is the category of Lie groupoids  with
strict homomorphisms given by $f\mapsto Z_f$ as described above.
\end{prop}
\begin{pf}
All the assertions are easy to check. For instance, to show that
$\Gamma^{(0)}\leftarrow Z''$ is a locally trivial $H$-principal bundle,
note that $Z$ and $Z'$ are locally of the form $Y\times_{G^{(0)}} G$
and $Y'\times_{H^{(0)}} H$ respectively. Therefore
$Z''$ is locally of the form $Y''\times_{G^{(0)}} H$
where $Y''\equal Y\times_{H^{(0)}} Y'$.
\end{pf}

Note that isomorphism in the category ${\mathcal{G}}$ is just Morita
equivalence \cite{muhly-renault-williams87, Xu:90}.

\begin{prop}\label{prop:homomorphism cover}
(see \cite[Definition 1.1]{hilsum-skandalis87})
Any generalized homomorphism
$$\Gamma^{(0)}\stackrel{\tau }{\leftarrow} Z
\stackrel{\sigma}{\rightarrow} G^{(0)}$$
is obtained by composition
of the canonical Morita equivalence between $\Gamma$ and $\Gamma[U_i]$,
where $(U_i)$ is an open cover of $\Gamma^{(0)}$, with a strict
homomorphism $\Gamma[U_i]\to G$.
\end{prop}

Consequently, giving a generalized homomorphism $\Gamma\to G$
is equivalent to giving a   Morita equivalence $\Gamma\sim_{morita}\Gamma'$
 together with
 a strict homomorphism $\Gamma'\to G$.

\begin{pf}
Denoting by $\Gamma[Z]$ the pull-back of $\Gamma$ via the
surjective submersion $Z\stackrel{\tau}{\to}\Gamma^{(0)}$, i.e. the
groupoid $Z\times_{\gm^{(0)}, t}\Gamma\times_{\gm^{(0)},s} Z$ with multiplication law
$(z_1,g,z_2)(z_2,h,z_3)\equal (z_1,gh,z_3)$. Then the canonical
strict homomorphism $\Gamma[Z]\to\Gamma$ is a Morita equivalence.

Moreover, it is not hard to check that
$$\Gamma[Z]\cong\{(z,z',\gamma,g)\in (Z\times Z)\times \Gamma\times G|
\;\gamma z'\equal zg\}. $$
Thus there is a strict homomorphism $f\colon \Gamma[Z]\to G$
given by the fourth projection. One can then verify that
the following diagram is commutative (in the category ${\mathcal{G}}$):
$$\xymatrix{
\Gamma\drto_Z&\Gamma[Z]\ar[l]_\cong\ar[d]^f\\
&G}$$
Now, since $Z\to \Gamma^{(0)}$
is a submersion, it admits local sections. Hence there
exists an open cover $(U_i)$ of $\Gamma^{(0)}$
and maps $s_i\colon U_i\to Z$ such that $\tau\circ s_i\equal \mathrm{Id}$,
and therefore a map $\tilde s\colon
\Gamma[U_i]\to\Gamma[Z]$ such that the composition
$\Gamma[U_i]\to\Gamma[Z]\to \Gamma$ is the canonical map.
Then, $f\circ\tilde s\colon \Gamma[U_i]\to G$ is the desired strict
homomorphism.
\end{pf}

\begin{lem}\label{lem:strict homomom cohom}
Let $f_1$, $f_2\colon \Gamma\to G$ be two strict homomorphisms.
Then $f_1$ and $f_2$  define the isomorphic generalized homomorphisms
if and only if there exists a smooth map $\varphi\colon \Gamma^{(0)}
\to G$ such that
$f_2(\gamma)\equal \varphi(t(\gamma))f_1(\gamma)\varphi(s(\gamma))^{-1}$.
\end{lem}

\begin{pf}
Suppose that there exists a smooth $\Gamma$, $G$-equivariant
map $Z_{f_1}\to Z_{f_2}$. Then it is necessary of the form
$(x,g)\mapsto (x,\varphi(x)g)$. Using $\Gamma$-equivariance,
we get 
$(t(\gamma),\varphi(t(\gamma))f_1(\gamma))\equal (t(\gamma),f_2(\gamma)
\varphi(s(\gamma)))$. The converse is proved by working backwards.
\end{pf}

The following  result is useful when dealing with generalized
homomorphisms (see also \cite{GZ}).

\begin{prop}\label{prop:functor Morita}
Let $\mathcal{C}$ be a category, and  $\Phi\colon {\mathcal{G}}_s\to
{\mathcal{C}}$ be a functor. The following are
equivalent:
\begin{itemize}
\item[(i)] for every smooth groupoid $\Gamma$ and every open cover
$(U_i)$, $\Phi(\pi)$ is an isomorphism, where $\pi$ is the canonical
map $\Gamma[U_i]\to \Gamma$.
\item[(ii)] The functor $\Phi$ factors through the category
${\mathcal{G}}$ (and thus $\Phi(G)\cong \Phi(H)$ if
$G$ and $H$ are Morita equivalent).
\end{itemize}
\end{prop}

\begin{pf}
The only non-trivial implication is (i)$\implies$(ii).
Let $\Gamma^{(0)}\leftarrow Z\rightarrow G^{(0)}$ be a generalized
homomorphism. From Proposition~\ref{prop:homomorphism cover},
there exists a strict homomorphism $f\colon\Gamma[U_i]\to G$
which is the same morphism in the category ${\mathcal{G}}$.
We define $\Phi(Z): \Phi (\gm )\to \Phi (G) $ to be the composition
$$\Phi(\Gamma)\stackrel{\cong}{\leftarrow}\Phi(\Gamma[U_i])
\stackrel{\Phi(f)}{\rightarrow} \Phi(G).$$
To check that this is well-defined, suppose that
$f_1\colon \Gamma[U_i]\to G$ and $f_2\colon\Gamma[V_j]\to G$
define the same generalized homomorphism. We need to  show that
$\Phi(f_1)\equal \Phi(f_2)$ via the identification
$\Phi(\Gamma[U_i])\equal \Phi(\Gamma[V_j])$. Using the cover $(U_i\cap V_j)$,
we may assume that $(U_i)\equal (V_j)$, and 
$f_1$, $f_2$ are strict homomorphisms from $\Gamma$ to $G$.

From Lemma~\ref{lem:strict homomom cohom}, there exists
$\varphi\colon \Gamma^{(0)}\to\Gamma$ such that
$f_2(\gamma)\equal \varphi(t(\gamma)) f_1(\gamma)\varphi(s(\gamma))^{-1}$.
Let $\tilde\Gamma\equal \{1,2\}^2\times\Gamma$ and let
$\tilde f\colon \tilde\Gamma\to G$ be the morphism
\begin{eqnarray*}
(1,1,\gamma)&\mapsto&f_1(\gamma)\\
(1,2,\gamma)&\mapsto&f_1(\gamma)\varphi(s(\gamma))^{-1}\\
(2,1,\gamma)&\mapsto&\varphi(t(\gamma))f_1(\gamma)\\
(2,2,\gamma)&\mapsto&\varphi(t(\gamma))f_1(\gamma)\varphi(s(\gamma))^{-1}.
\end{eqnarray*}
Let $i_j\colon \Gamma\to \tilde\Gamma$ be defined by $i_j(\gamma)
\equal (j,j,\gamma)$
and $\pi\colon\tilde\Gamma\to \Gamma$, be the map $\pi(i,j,\gamma)\equal \gamma$.
Since $\tilde\Gamma\equal \Gamma[W_k]$ with $W_1\equal W_2\equal \Gamma^{(0)}$,
$\Phi(\pi)$ is an isomorphism. Now, from $\pi\circ i_1\equal \pi\circ i_2$,
we get $\Phi(i_1)\equal \Phi(\pi)^{-1}\equal \Phi(i_2)$, and therefore 
$\Phi(f_1)\equal \Phi(\tilde f\circ i_1)\equal \Phi(\tilde f\circ i_2)\equal \Phi(f_2)$.
\end{pf}

\begin{rmk}
Given two Lie groupoids $ \Gamma_i \toto \gm_i^{(0)}, \ i\equal 1, 2$,  
 a generalized homomorphism from $ \Gamma_1 \toto \gm_1^{(0)}$ to
$ \Gamma_2 \toto \gm_2^{(0)}$  induces a morphism
 between their associated differential stacks
$\XX_1\to \XX_2$, and vice versa. And  a generalized isomorphism, i.e.
a  Morita equivalence,  corresponds to an isomorphism of  stacks.
Therefore the category ${\mathcal{G}}$ and the category of differentiable
stacks are equivalent categories (see \cite{BX1} for details). 
\end{rmk}

\subsection{$S^1$-central extensions of groupoids}
\label{subsec:central extensions}
\begin {defn}\label{def:central extension}
Let ${\gm}\toto {M}$ be a Lie groupoid. An \emph{$S^1$-central
extension (or ``twist'')} of $\gm\toto M$ consists of
 
1. a  Lie groupoid ${R}\toto {M}$, together with a morphism of
Lie
groupoids $(\pi,\id):[R\toto M]\to[\Gamma\toto M]$ which restricts to
the identity on $M$,
 
2. a left $S^1$-action on $R$, making $\pi:R\to \Gamma$ a (left)
principal $S^1$-bundle.
\noindent These two structures are compatible in the sense that
$(\lambda_1 \cdot x) (\lambda_2 \cdot y)\equal \lambda_1\lambda_2
 \cdot (xy )$,
for all  $ \lambda_1 , \lambda_2  \in S^{1}$ and 
$(x, y) \in R^{(2)}\equal  R\times_{s, M, t} R $.
\end {defn}                                     

We denote by $Tw^{sm}(\gm)$ the set of $S^1$-central extensions
of $\gm$ (the superscript ``sm'' stands for ``smooth'').

Note that $R$ being restricted to $\epsilon_0 (M)$ is  a
trivial $S^1$-bundle, where $\epsilon_0 : M\to \gm$ is
the unit  map.  In fact, it   admits a canonical
trivialization since $R|_{\epsilon_0 (M)}$ admits
a smooth section, namely,  the base space of the groupoid
$R$.  By $\ker \pi $, we denote this trivial bundle
$R|_{\epsilon_0 (M)}$, i.e., $\ker \pi\cong M\times S^1$.
It is clear that $\ker \pi$, as a bundle of groups,
is a normal subgroupoid of $R\toto M$, and lies in
the center. Indeed its quotient groupoid is isomorphic
to $\gm \toto M$. This  coincides with the usual
definition of Lie group $S^1$-central extensions.

When  $\pi:R\to \Gamma$  is topologically trivial (for instance, this
is true if as a  space
  the 2nd cohomology of $\Gamma$  vanishes), then 
$R\cong \gm \times S^1$ and the central extension is determined 
by a groupoid 2-cocycle of $\gm$ valued in $S^1$,  i.e.,
a smooth map $c: \gm^{(2)}\equal
\{(x, y)|s(x)\equal t(y), x, y \in \gm\}\to S^1$
satisfying the cocycle condition:
\begin{equation}
\label{eq:cocycle}
c(x, y)c(xy, z)c(x, yz)^{-1}c(y, z)^{-1}\equal 1,
\ \ \ \ \forall (x, y, z)\in \gm^{(3)}.
\end{equation} 
The groupoid structure on $R$ is given by
\begin{equation} 
\label{eq:extension}
(x, \lambda_1)\cdot (y, \lambda_2 )\equal (xy, \lambda_1 \lambda_2 c(x, y)), \ \ \forall
(x, y)\in \gm^{(2)}.
\end{equation}

For  every locally compact groupoid $\gm$ with a Haar system (thus for
every Lie groupoid), Kumjian, Muhly, Renault and Williams \cite{kum98}
constructed a group, called the Brauer group of $\Gamma$.
Some of the constructions below is an
easy adaptation of their results  to our context,
so we will omit most of proofs.
                                                           
Note that  $Tw^{sm}(\Gamma)$ admits  an abelian group structure in a
canonical way:
if $S^1\to R\stackrel{\pi}{\to} \Gamma\toto M$ and
$S^1\to R'\stackrel{\pi'}{\to} \Gamma\toto M$
 are $S^1$-central extensions, then the
addition of $R$ and $R'$, called the {\em tensor of  $R$ and $R'$} and
denoted by $R\otimes R'$,  is
\begin{equation}
\label{eq:Rsum}
(R\times_{\gm } R')/{S^1}:\equal \{(r,r')\in R\times_{\gm } R'\}
/_{(r,r')\sim(\lambda r,\lambda^{-1}r')}
\end{equation}
($\lambda\in S^1$), and the inverse of $R$ is $\bar R$
(where the action of $S^1$ on $\bar R$ is $\lambda \bar r
\equal \overline{\overline{\lambda} r}$
and $\bar r\in \bar R$ denotes the same element
$r\in R$).
 
The zero element is the strictly trivial extension, i.e., the
extension satisfying the following equivalent conditions.

\begin{prop}\label{prop:strictly trivial}
Let $S^1\to R\stackrel{\pi}{\to} {\Gamma}\toto M$
be an $S^1$-central extension.
The following are equivalent:
\begin{itemize}              
\item[(i)] there exists a groupoid homomorphism $\sigma \colon
\Gamma\to R$ 
such that $\pi\smalcirc\sigma\equal \mathrm{Id}$;
\item[(ii)] there exists an $S^1$-equivariant groupoid homomorphism
$\varphi\colon R\to S^1$;
\item[(iii)] $R\cong \Gamma\times S^1$ (as a product of groupoids).
\end{itemize}
\end{prop}

\begin{pf}
(i)$\implies$(ii): take $\varphi(r)\equal r(\sigma\smalcirc\pi(r))^{-1}$.

(ii)$\implies$(iii): the map $r\mapsto (\pi(r),\varphi(r))$
is a groupoid isomorphism from $R$ to $\Gamma\times S^1$.

(iii)$\implies$(i): obvious.
\end{pf}
\par\medskip
 
The set of $S^1$-central extensions of $\Gamma$ of the form
$R\equal (t^*\Lambda\times \overline{s^*\Lambda})/S^1 \toto M$, where
$\Lambda$ is an  $S^1$-principal bundle on $M$, is a subgroup
of $Tw^{sm}(\Gamma)$. The quotient of $Tw^{sm}(\Gamma)$ by this subgroup
is denoted by ${\mathcal{E}}^{sm}(\Gamma)$. 

We now  introduce the definition of Morita equivalence of $S^1$-central
 extensions, and define an abelian group structure on the set of Morita 
equivalence classes of extensions $S^1\to R'\toto \Gamma'$, with $\Gamma'$ Morita
equivalent to $\Gamma$.

\begin{defn}
Let $S^1\to R\to \Gamma\toto M$
 and $S^1\to R'\to \Gamma' \toto M'$ be $S^1$-central extensions.
We  say that a  generalized homomorphism $M\leftarrow Z
\rightarrow M'$ from $R$ to $R'$ is $S^1$-equivariant
if $Z$ is endowed with an action of $S^1$ such that
$$(\lambda r)\cdot z\cdot r'\equal r\cdot (\lambda z)\cdot r'\equal r\cdot
z \cdot
(\lambda r')$$
 whenever $(\lambda, r, r', z)\in S^1\times R\times R'\times Z$ and
the products make sense.
\end{defn}                        

\begin{lem}
\label{lem:eq}
Let $S^1\to R\to \Gamma \toto M$ and $S^1\to R'\to \Gamma' \toto M'$ 
be $S^1$-central extensions, and  $M\leftarrow Z
\rightarrow M'$  an $S^1$-equivariant  generalized homomorphism from $R$ to $R'$. Then the $S^1$-action on $Z$ is free and
$M\leftarrow  Z/S^1 \rightarrow M'$  defines a
generalized  homomorphism from $\gm$ to $\gm'$.
\end{lem}
\begin{pf}
Assume that $\lambda z \equal z$ for $\lambda \in S^1$ and $z\in Z$.
 From the compatibility condition $ (\lambda z)\cdot r'\equal  z \cdot (\lambda r')$
and the fact that the $R'$-action on $Z$ is free, we  obtain
that $\lambda r' \equal r'$. Hence  $\lambda \equal 1$ since the $S^1$-action on
$R'$ is free. The rest of the assertion follows
immediately from the compatibility condition again.
\end{pf}                 
\begin{defn}\label{def:morita equivalent extension}
Two $S^1$-central extensions  $S^1\to R\to \Gamma \toto M$
and $S^1\to R'\to \Gamma' \toto M'$ are called {\em Morita equivalent}
if there is a generalized $S^1$-equivariant
 isomorphism $M\leftarrow Z \rightarrow M'$.
In this case, $Z$ is called  an {\em equivalence bimodule}.
\end{defn}
 
As an immediate consequence of Lemma \ref{lem:eq},
in particular, if $S^1\to R\to \Gamma \toto M$ and
$S^1\to R'\to \Gamma' \toto M'$ are Morita equivalent $S^1$-central
 extensions,
then $\Gamma$ and $\Gamma'$ must be Morita equivalent groupoids.
 
The following result gives a useful construction
of $S^1$-equivariant generalized homomorphism,
 and in particular shows that for two Morita
equivalent  $S^1$-central extensions, one may
recover one from the other in terms of the equivalence bimodule.

Let $S^1\to R\to \Gamma \toto M$ be  an $S^1$-central extension,
and $ \tau: Z\to M$ a left principal $R$-bundle over $Z\to M':\equal Z/R $.
Then $Z$ admits an $S^1$-action defined as follows:
for all $\lambda\in S^1$ and $z\in Z$, denote by $\lambda_z\in R$
the element $\lambda \tau(z)$, where $\tau(z)\in R^{(0)}$ is considered
as an element of $R$. We let
$$\lambda \cdot z\equal  \lambda_z\cdot z.$$
It follows from the properties of $S^1$-central extensions
that this indeed defines an $S^1$-action. Moreover, by assumption, this
action is free and therefore $Z/S^1$ is a smooth manifold, which is denoted
by $X$. It is simple to see that the following identity holds: 
$$(\lambda r)\cdot z\equal r\cdot (\lambda z)
\quad\forall (\lambda,r,z)\in S^1\times R\times Z\mbox{ with }
s(r)\equal \tau(z).$$

\begin{prop}
As above, assume that  $S^1\to R\to \Gamma \toto M$ is an $S^1$-central
 extension, and $ \tau: Z\to M$  a principal $R$-bundle over
 $Z\to M':\equal Z/R \cong X/\gm$, where
  $X\equal Z/S^1$.  Let $\gm' \equal X\times_\gm X$  and
$R'\equal Z\times_R Z$. Then
\begin{itemize}
\item[(i)]
$R'\to \Gamma' \toto M'$ is an $S^1$-central extension of groupoids,
and $M'\leftarrow Z \rightarrow M$ with the natural actions
defines an $S^1$-equivariant generalized homomorphism from $R'$ to $R$.
\item[(ii)]
if moreover $\tau: Z\to M$ is a surjective submersion, then
$R'$ and $R$ are Morita equivalent $S^1$-central extensions.
\end{itemize}
As a consequence,  if $Z$ is an $S^1$-equivariant Morita equivalence
bimodule from $R$ to $R''$, then $ R''\cong R'$.
\end{prop}
\begin{pf}
The proof is a straightforward verification, and is similar to
\cite[Theorem 3.2-3.3]{Xu:90}.
\end{pf}

Given a Lie groupoid $\gm \toto M$,  there is a natural
abelian group structure on the set of 
Morita equivalence classes of $S^1$-central
extensions $S^1\to R'\to \Gamma' \toto M'$,
where $ \Gamma' \toto M'$ is a   Lie groupoid Morita
equivalent to $\gm \toto M$.
To see this, assume that
$S^1\to R_i\to \Gamma_i \toto M_i$, $i\equal 1, 2$, are two
such extensions. Since $\Gamma_i \toto M_i$, $i\equal 1, 2$, 
are Morita equivalent, there exists a generalized
isomorphism $M_1 \leftarrow X \rightarrow M_2$.
By pulling back using the above maps, one
obtains two $S^1$-central extensions over the
groupoid $\gm_1 [X] \cong \gm_2 [X]  \toto X$,
namely $R_1 [X] \toto X$ and $R_2 [X] \toto X$.
Thus one may define $[R_1 ]+[R_2 ]$ 
to be  the class of  the $S^1$-central extension $R_1 [X] \otimes R_2 [X]$
(see Eq. (\ref{eq:Rsum})).  It is simple to check that this
operation is well-defined. The inverse is defined by
$\overline{[R]}\equal [\bar{R}]$. 
Let us denote by ${\mathrm{Ext}}^{sm}(\Gamma,S^1)$ the
group thus obtained.
\par\bigskip

The zero element in ${\mathrm{Ext}}^{sm}(\Gamma,S^1)$
is characterized by the following
 
\begin{prop}\label{prop:trivial extensions}
Consider an $S^1$-central extension of Lie groupoids $S^1\to R\to\Gamma\toto M$.
The following are equivalent:
\begin{itemize}
\item[(i)] there exists an $S^1$-equivariant generalized homomorphism
$R\to S^1$;
\item[(ii)] there exists a cover $(U_i)$ of $M$ such that the extension
$S^1\to R[U_i]\to \Gamma[U_i]\toto \coprod U_i$ is strictly trivial;
\item[(iii)] the extension is Morita equivalent to a strictly trivial
$S^1$-central extension $0\to S^1\to \Gamma'\times S^1\to \Gamma'\toto M'$;
\item[(iii)'] the class of the extension in $Ext^{sm}(\Gamma,S^1)$ is 0;
\item[(iv)] the extension is Morita equivalent to the strictly trivial
$S^1$-central extension $0\to S^1\to \Gamma\times S^1\to \Gamma\toto M$;
\item[(iv)'] the class of the extension in $\eE^{sm}(\Gamma)$ is 0.
\end{itemize}
\end{prop}

\begin{pf}
(ii)$\implies$(iii)$\iff$(iii)'$\implies$(i) and
(iv)'$\iff$(iv)$\implies$(iii) are obvious.
(i)$\implies$(ii) is a consequence of Propositions~\ref{prop:strictly trivial}
and~\ref{prop:homomorphism cover}. To show (iii)$\implies$(iv),
let $Z$ be an equivalence bimodule between $\Gamma$ and $\Gamma'$,
then $Z\times S^1$ is obviously an equivalence bimodule between
the trivial central extensions $S^1\to \Gamma\times S^1\to \Gamma\toto M$
and $S^1\to \Gamma'\times S^1\to \Gamma'\toto M'$.
\end{pf}

$S^1$-extensions which satisfy any of the conditions in
 the previous proposition
are said to be trivial. Therefore, $\eE^{sm}(\Gamma)$ is the quotient
$Tw^{sm}(\Gamma)$ by trivial extensions.
Thus two $S^1$-central extensions
$S^1\to R_i\to\Gamma\toto M$
are equal in ${\mathcal{E}}^{sm}(\Gamma)$ if and only if
they are Morita equivalent.
\par\medskip

The groups ${\mathcal{E}}^{sm}(\Gamma')$, where $\Gamma'$ is a groupoid
Morita equivalent to $\Gamma$, form an inductive system. It
follows from Proposition~\ref{prop:trivial extensions}
 that
$$\mbox{Ext}^{sm}(\Gamma,S^1)\cong\lim_{{\Gamma'\sim\Gamma}}
{\mathcal{E}}^{sm}(\Gamma')\cong\lim_{{\mathcal{U}}}
{\mathcal{E}}^{sm}(\Gamma[{\mathcal{U}}]),$$
where $\mathcal{U}$ runs over open covers of $M$.
\par\medskip

\begin{rmk}
\label{rmk:gerbe}
 An $S^1$-central extension $S^1\to R\to \Gamma \toto M$
gives rise to an $S^1$-gerbe $\RR$ over the  differentiable
stack $\XX_\Gamma$ associated to
 the  groupoid $\Gamma \toto M$  \cite{BX, BX1},
  and Morita equivalent $S^1$-central extensions correspond to
isomorphic gerbes.

Conversely, given an $S^1$-gerbe  $\RR \stackrel{\pi}{\to}\XX$ over
a differential stack  $\XX$, if   $R\toto M$ and $R'\toto M'$ are
the Lie groupoids corresponding to
the  presentations $M\to \RR$ and
$M'\to \RR$ of $\RR$ respectively, and
 $\Gamma\toto M$  and $\Gamma'\toto M'$  are the
Lie groupoids corresponding to the induced presentations $M\to \XX$ and
$M'\to \XX$ of $\XX$ respectively,
then the $S^1$-central extensions $S^1\to R\to \Gamma \toto M$
and $S^1\to R'\to \Gamma' \toto M'$   are Morita equivalent  \cite{BX, BX1}.
The equivalence  bimodule is $Z\equal M\times_{\RR}M'$, which is a
principal $S^1$-bundle over $M\times_{\XX}M'$.

Therefore, given a Lie groupoid
$\gm \toto M$,   one may  identify an 
$S^1$-gerbe over   the stack $\XX_\gm$ as an element in
${\mathrm{Ext}}^{sm}(\Gamma,S^1)$, i.e.
a  Morita equivalence class of $S^1$-central
extensions $S^1\to R'\to \Gamma' \toto M'$,
where $ \Gamma' \toto M'$ is a  Lie groupoid Morita
equivalent to $\gm \toto M$.
 We will call such a Morita equivalence class 
 an isomorphism class of $S^1$-gerbes by abuse of notations.  Moreover,
the group structure on ${\mathrm{Ext}}^{sm}(\Gamma,S^1)$
corresponds to the abelian group structure
on  the $S^1$-gerbes over $\XX_\gm$. Therefore,
one may simply identify these two groups.
\end{rmk}

\subsection{Cohomology and characteristic classes}

In this subsection, we briefly review some basic cohomology
theory of groupoids, which will be needed later in the paper.

There exist many equivalent ways of introducing
cohomology  groups associated to a Lie
groupoid $\gm\toto \gm^{(0)}$ \cite{BX, BX1, CM}. 
A  simple and geometric  way is to consider
the   simplicial manifold canonically associated
to  the groupoid  and apply the usual cohomology theory. More precisely, 
let $\Gamma\toto \gm^{(0)}$ be a Lie groupoid. Define for all  $p\geq0$
$$\Gamma^{(p)}\equal \underbrace{\Gamma\times_{\gm^{(0)}}\ldots\times_{\gm^{(0)}}\Gamma}_{\text{$p$
times}}\,,$$
i.e., $\Gamma^{(p)}$ is the manifold of composable sequences of $p$ arrows
in the groupoid $\Gamma\toto \gm^{(0)}$.
We have $p+1$ canonical maps $\gm^{(p)}\to \gm^{(p-1)}$ 
giving rise to a              diagram
\begin{equation}\label{sim.ma}
\xymatrix{
\ldots \gm^{(2)}
\ar[r]\ar@<1ex>[r]\ar@<-1ex>[r] & \gm^{(1)}\ar@<-.5ex>[r]\ar@<.5ex>[r]
&\gm^{(0)}\,.}
\end{equation}

In fact,  $\gm \upcom$ is a simplicial manifold, so 
one can  introduce (singular) cohomology groups $H^k (\gm\upcom , \zz)$,
 $H^k(\Gamma\upcom,\rr)$ and $H^k(\Gamma\upcom,\rr/\zz)$.             
We refer the reader to \cite{Dupont} for the detailed study of
cohomology of simplicial manifolds.
In fact, for any abelian sheaf $F$ on the category of differentiable manifolds,
we have the cohomology groups $H^k(\Gamma\upcom,F)$ \cite{SGA4, BEFFGK, BX1}.
  One way to define
them is by choosing for every $p$ an injective resolution $F^p\to
{I^p}\upcom$ of sheaves on $\Gamma^{(p)}$, where $F^p$ is the small 
sheaf induced by
$F$ on $\Gamma^{(p)}$; then choosing homomorphisms $f^{*} {I^{p-1}}\upcom\to
{I^p}\upcom$ for every map $f:\Gamma^{(p)}\to \Gamma^{(p-1)}$
in~(\ref{sim.ma}). This gives rise to a double complex
$I\upcom(\Gamma\upcom)$, whose total cohomology groups are the
$H^k(\Gamma\upcom,F)$.
Examples of abelian sheaves on the category of manifolds are:
$\zz$,
$\rr$, $\rr/\zz$,  $\Omega^k$,
 ${\mathcal{R}}$ and ${\mathcal{S}}^1$.  The first three are
sheaves of locally constant functions, 
 ${\mathcal{R}}$ and ${\mathcal{S}}^1$ are the sheaves
 of differentiable $\rr$-valued and  $S^1$-valued functions,
respectively (see \cite{BEFFGK, BX, BX1})
With  respect to the first three, the notation $H^k(\Gamma\upcom ,F)$
does not conflict with the notation introduced before.    
Note that the cohomology groups
 $H^k(\Gamma\upcom,F)$ satisfy  the functorial property
with respect to generalized homomorphisms according to Proposition 
\ref{prop:functor Morita}.

Another  cohomology, which is  relevant to us, is
the De Rham cohomology.
Consider  the double complex $\Omega^\com(\Gamma\upcom)$:
\begin{equation}
\label{eq:DeRham}
\xymatrix{
\cdots&\cdots&\cdots&\\
\Omega^1(\gm^{(0)})\ar[u]^d\ar[r]^\partial &\Omega^1(\gm^{(1)})\ar[u]^d\ar[r]^\partial
&\Omega^1(\gm^{(2)})\ar[u]^d\ar[r]^\partial&\cdots\\
\Omega^0(\gm^{(0)})\ar[u]^d\ar[r]^\partial&\Omega^0(\gm^{(1)})\ar[u]^d\ar[r]^\partial
&\Omega^0(\gm^{(2)})\ar[u]^d\ar[r]^\partial&\cdots
}
\end{equation}

Its boundary maps are $d:
\Omega^{k}( \gm^{(p)} ) \to \Omega^{k+1}( \gm^{(p)} )$, the usual exterior
derivative of differentiable forms and $\partial
:\Omega^{k}( \gm^{(p)} ) \to \Omega^{k}( \gm^{(p+1)} )$,  the alternating
sum of the pull-back maps of (\ref{sim.ma}).
We denote the total differential by $\delta\equal (-1)^pd+\del$.
The  cohomology groups of the total complex $C^\com(\Gamma\upcom)$:
$$H_{DR}^k(\Gamma\upcom)\equal H^k\big(\Omega^\com(\Gamma\upcom)\big)$$
are called the \emph{De~Rham cohomology} groups of $\Gamma\toto \gm^{(0)}$.              

 The following proposition lists
some well-known properties regarding De~Rham cohomology groups of a
Lie groupoid.

\begin{prop} \cite{BEFFGK, BX, BX1, CM, hae84}
\begin{enumerate}
\item  For   any Lie groupoid $\gm\toto \gm^{(0)}$, we have
\begin{equation}\label{can.iso}
H_{DR}^k(\Gamma\upcom)\longiso H^k(\Gamma\upcom,\rr);
\end{equation}
\item if $\gm\toto \gm^{(0)} $ and $G\toto G^{(0)}$ 
are Morita equivalent, then
$$ H_{DR}^k(\Gamma\upcom) \longiso H_{DR}^k({G}\upcom), \ \ \ \mbox{and }
H^k (\Gamma\upcom , {\mathcal{S}}^1 )\longiso H^k ({G}\upcom ,  {\mathcal{S}}^1 ).$$
\end{enumerate}
\end{prop}

We call a De~Rham $k$-cocycle an \emph{integer cocycle}, if it
maps under~(\ref{can.iso}) into the image of the canonical map
$H^k(\Gamma\upcom, \zz)\to H^k(\Gamma\upcom,\rr)$.

\begin{numex}
\label{ex:coh}
\begin{enumerate}
\item When $\gm$ is a manifold $M$, it is clear that
$H^{k} (\Gamma\upcom, \zz)$ (or $H^{k} (\Gamma\upcom, \rr ))$
reduces to the usual cohomology $H^{k} (M, \zz)$ (or $H^{k} (M, \rr )$
respectively).
If $\{U_i \}$ is an open covering of $M$ and $ X\equal  \coprod_i U_i \to M$
is the \'etale map, then   $\gm: \equal X\times_M X\toto X$, which is
$\coprod_{ij} U_i \cap U_j \toto \coprod_i U_i$, is
Morita equivalent to $M\toto M$. Hence
$H^{k} (\Gamma\upcom, \zz)$ (or $H^{k} (\Gamma\upcom, \rr )$ respectively)
is isomorphic to  $H^{k} (M, \zz)$ (or $H^{k} (M, \rr )$ respectively).
The double complex (\ref{eq:DeRham}), when $\{U_i \}$ is a 
nice covering,  is the one used by Weil in his
proof of De Rham theorem \cite{Weil}.

\item When $\Gamma$ is a transformation groupoid $G\times
M\toto M$,  $H^{k} (\Gamma\upcom, \zz)$ (or $H^{k} (\Gamma\upcom, \rr )$
respectively) is  the  $G$-equivariant cohomology group $H^k_G (M, \zz)$
(or $H^k_G (M, \rr )$ respectively).
If $G$ is compact, $H^k_G (M, \rr )$ can be alternatively computed
by either Cartan model or Weil model  (see \cite{Guillemin} for
more details).

\item On the other hand, if $\Gamma\toto M$ is an \'etale groupoid
representing an orbifold \cite{Moer} and $\Lambda (\gm )\toto
\Gamma $ its associated inertia groupoid,
then $H^{k} ( \Lambda (\gm )\upcom, \rr )$
is the orbifold cohomology.
\end{enumerate}
\end{numex}

It is known that $H^2 (\gm\upcom, {\mathcal{S}}^1 )$ classifies
$S^1$-gerbes over the stack $\XX_{\gm}$ \cite{Giraud}. 
 As a consequence
(see Remark \ref{rmk:gerbe}), we have

\begin{prop}\label{prop:ext-H2}
For a  Lie groupoid $\gm \toto M$, we have
$$ {\mathrm{Ext}}^{sm}(\Gamma,S^1) \cong H^2 (\gm\upcom, {\mathcal{S}}^1 ).$$
\end{prop}

For instance, 
 when $\gm$ is a manifold $M$,
by Example~\ref{ex:gerbe}.1 below, ${\mathrm{Ext}}^{sm}(M,S^1)$
is isomorphic to the \v{C}ech cohomology group
$\check{H}^2(M,{\mathcal{S}}^1)\equal H^2(M,{\mathcal{S}}^1)$.

The exponential sequence $0\to \zz\to \rrr \to {\mathcal{S}}^1\to 0$
gives rise to a  long exact sequence:

\begin{equation}
\label{eq:exact}
\cdots \to H^2 (\Gamma\upcom, \zz )
\stackrel{\psi_2}{\to}
\cdots  \to H^2 (\Gamma\upcom, \rrr ) \to H^2(\Gamma\upcom,{\mathcal{S}}^1)
\stackrel{\phi}{\to}
 H^3(\Gamma\upcom,\zz)
\stackrel{\psi_3}{\to}   H^3(\Gamma\upcom, \rrr) \to \cdots
\end{equation}


\begin{lem}
$$ H^k (\Gamma\upcom, \rrr) \cong H^k (\gm, \rr),$$
where $H^k (\gm, \rr)$ denotes the (smooth)
groupoid cohomology with the trivial
coefficients $\rr$, i.e. the cohomology of the complex
$(C^\infty(\gm^{(n)},\rr))_{n\in\nn}$ with the differential
$$(dc)(g_1,\ldots,g_{n+1})
\equal c(g_2,\ldots,g_{n+1})+\sum_{k\equal 1}^n(-1)^kc(g_1,\ldots,g_kg_{k+1},
\ldots,g_{n+1})$$
$$\qquad\qquad+(-1)^{n+1}c(g_1,\ldots,g_n).$$
\end{lem}
\begin{pf}
There is a spectral sequence
$$ E_1^{p,q} \equal   H^q(\gm^{(p)}, \rrr) \implies H^{p+q}(\Gamma\upcom, \rrr)$$
Since $\gm^{(p)}$ is a manifold and
 the sheaf $\rrr|\gm^{(p)}$ is soft, $H^q(\gm^{(p)}, \rrr)\equal 0$ for $q>0$.
Therefore the spectral sequence degenerates.
It follows that  $H^* (\Gamma\upcom, \rrr)$  can be calculated using the
complex $H^0(\gm^{(p)}, \rrr)\equal C^{\infty}(\gm^{(p)})$. 
\end{pf}

By identifying the groups $H^k (\Gamma\upcom, \rrr)$
with $H^k (\gm, \rr)$, the homomorphism 
$$\psi_k : H^k(\Gamma\upcom,\zz) \to H^k(\Gamma\upcom, \rrr) $$
in the exact sequence (\ref{eq:exact}) is the composition
of the following sequences of morphisms:
$$H^k(\Gamma\upcom, \zz)\to H^k(\Gamma\upcom, \rr )
\longiso H^{k}_{DR}(\Gamma\upcom )\stackrel{\pr}{\to}
H^k (\Gamma , \rr )\longiso H^k(\Gamma\upcom, \rrr), $$
where
$\pr : H_{DR}^k (\Gamma\upcom )\lon H^k (\Gamma , \rr )$
 is given, on the cochain level, by the projection 
$\oplus_{i+j\equal k}\Omega^{i} (\gm_j )\to \Omega^0 (\gm_k )$.
 See \cite{C, WX} for details on (smooth) groupoid cohomology.

Note that in general $\phi: H^2(\Gamma\upcom,{\mathcal{S}}^1))
\to H^3 (\gm\upcom , \zz )$ is neither surjective nor
injective.
Write
$$  H^3_{gerbe} (\gm\upcom , \zz )
\equal \phi (H^2(\Gamma\upcom,{\mathcal{S}}^1)).$$

\begin{prop}
\begin{itemize}
\item[(i)] 
$  H^3_{gerbe} (\gm\upcom , \zz )$ is
a   subgroup of $H^3 (\gm\upcom , \zz )$
consisting of those elements whose image in
 $H^3_{DR} (\gm\upcom )$ projects to  zero under
$\pr : H_{DR}^3 (\Gamma\upcom )\lon H^3 (\Gamma , \rr )$.
\item[(ii)]   The kernel  of $\phi$ is  isomorphic to
$H^2 (\gm , \rr )/\psi_2 ( H^2 (\gm\upcom, \zz ))$.
\end{itemize}
\end{prop}

For an $S^1$-central extension $R\to \gm \toto M$,
let $[R]\in  H^2(\Gamma\upcom,{\mathcal{S}}^1)$ denote its
class.  The image of  $[R]$ in $H^3(\Gamma\upcom,\zz)$
 under the homomorphism $\phi$
 is called the \emph{Dixmier-Douady class} of $R$ \cite{BX, BX1}. 
  The Dixmier-Douady class behaves well
with respect to the pull-back and the tensor operation.
Unlike  the manifold case, in general the Dixmier-Douady class
does not completely determine an $S^1$-gerbe.
However this is true  when $\gm \toto M$ is a proper groupoid.
Let us recall its definition below.

\begin{defn}
Let $\Gamma\toto M$ be a locally compact groupoid. Then $\Gamma$ 
is said to be proper if any of the following  equivalent conditions 
is satisfied:
\begin{itemize}
\item[(i)] the map $(s,t)\colon \Gamma\to M\times M$
is   proper;
\item[(ii)] for every $K\subset M$ compact, $\Gamma_K^K$
is compact.
\end{itemize}
\end{defn}

For instance, compact groupoids are of course proper;  a transformation
groupoid $G\times M\toto M$ is proper if and only if the action
is proper.

\begin{lem}\label{lem:proper morita}
\begin{enumerate}
\item The notion of properness is invariant by Morita equivalence;
\item  for a proper groupoid $\gm\toto M$,
 the orbit space $M/\gm$ is  a Hausdorff topological space, and
 is invariant by Morita equivalence.
\end{enumerate}
\end{lem}

\begin{pf}
Suppose that $f\colon Y\to M$ is a surjective submersion.
If $\Gamma$ is proper, then for every $K\subset Y$ compact,
$(\Gamma[Y])_K^K$ is a closed subset of $K\times K\times
\Gamma_{f(K)}^{f(K)}$, and therefore it is compact. Hence,
$\Gamma[Y]$ is proper.

Conversely, if $\Gamma[Y]$ is proper, then for every $L\subset M$
compact, there exists $K\subset Y$ compact such that $f(K)\equal L$
(since $f$ is open surjective). Now,
$\Gamma_L^L$ is a continuous image of the compact set $(\Gamma[Y])_K^K$,
and thus  is compact. It follows that $\Gamma$ is proper. This proves (1).

The first assertion in (2) is proved for instance in
\cite[Proposition~6.3]{tu99}. For the second one, it is clear
that if $f\colon Y\to M$ is a surjective submersion, then
$f$ induces a homeomorphism $Y/(\Gamma[Y])\cong M/\Gamma$. 
\end{pf}

When $\gm \toto M$ is a proper Lie groupoid, since the
smooth groupoid cohomology $H^k (\gm , \rr )$ vanishes
when $k\geq 1$ according to Crainic \cite{C},
we see that $\phi$ is an isomorphism.

\begin{prop}\label{prop:H2 equal H3}
If $\gm \toto M$ is a proper Lie groupoid, then
$$\phi:\ \  \   H^2(\Gamma\upcom,{\mathcal{S}}^1) \lon H^3(\Gamma\upcom,\zz) $$
is an isomorphism.
\end{prop}
As a consequence, we have

\begin{cor}
\label{cor:eg}
\begin{enumerate}
\item If a Lie group $G$ acts  on a smooth manifold
 $M$ properly, then the equivariant cohomology $H^3_G (M, \zz)$ is isomorphic
to the abelian group of $S^1$-gerbes  over the stack
$M/G$   associated to the transformation groupoid $G\times M\toto M$.
\item If $\gm\toto M$ is an \'{e}tale groupoid corresponding to
an orbifold ${\mathbf{X}}$, then $H^3 (\gm\upcom , \zz )$
is isomorphic to the abelian group of $S^1$-gerbes  over the 
orbifold ${\mathbf{X}}$.
\end{enumerate}
\end{cor}

In other words, in both cases above, the third 
integer  cohomology classes can be geometrically described
by  Morita equivalent classes of  groupoid $S^1$-central extensions.
In particular,   for a smooth manifold $M$, since
$M\toto M$ is a special case when $G\equal 1$,  then 
 $H^3 (M, \zz )$ characterizes $S^1$-gerbes over the manifold $M$ \cite{brylinski1}.
 However the Dixmier-Douady  class does not completely characterize 
$S^1$-gerbes even when $\gm$ is a non-compact group as
we see below.

\begin{numex}
\cite{TW}
Consider  the abelian group $\rr^2$ as a groupoid
$\rr^2\toto \cdot$. It is clear that $H^k ({\rr^2}\upcom , \zz)
\cong H^k (B\rr^2 , \zz )\equal 0$
since $\rr^2$ is contractible. 
Therefore the kernel of $\phi$ is isomorphic
to the 2nd  group cohomology  of $\rr^2$,
which is in turn isomorphic to the
 2nd Lie algebra cohomology with trivial coefficients
since $\rr^2$ is simply connected. The latter
is isomorphic to the invariant De Rham cohomology
of $\rr^2$ under the translation, and therefore
is one-dimensional as a $\rr$-vector space. 
More explicitly, the group 2-cocycle is
given by
$$\sigma ((x, y), (x' , y' ))\equal \frac{1}{2} (x'y -xy' ) .$$

In other words, the group 2-cocycle $\exp({2\pi i \sigma})$
defines a non-trivial $S^1$-central extension of
$\rr^2$ (hence a non-trivial $S^1$-gerbe over $\XX_{\rr^2}$)
with the trivial Dixmier-Douady class. 
\end{numex}

It is often  useful to use differential forms
to describe the Dixmier-Douady class as in the manifold case.
Recall that a pseudo-connection
 is $\theta+B \in \Omega^1 (R)\oplus \Omega^2 (M)$
such that $\theta $ is a connection one-form of the principal
$S^1$-bundle $R\stackrel{\pi}{\to}\gm$ \cite{BX}.
Its pseudo-curvature
 $\eta+\omega+\Omega \in \Omega^1 (\gm^{(2)}) \oplus
\Omega^2 (\gm )\oplus \Omega^3 (M)
\subset C^3 (\gm\upcom )$
 is defined by
$$\delta(\theta+B)\equal \pi\upst(\eta+\omega+\Omega). $$

Then we have the following  \cite{BX}:

\begin{them}
1. $[\eta+\omega+\Omega] $ is independent of the  pseudo-connection
and defines an integer class in $H^3_{DR} (\gm\upcom )$.
Under the
canonical homomorphism $H^3(\Gamma\upcom,\zz)\to H^3_{DR}(\Gamma\upcom)$, the
Dixmier-Douady class of $R$ maps to $[\eta+\omega +\Omega]$.

2. Assume that 
$\gm \toto M$ is proper.
Given any integer 3-cocycle $\eta+\omega+\Omega$ as
above,
by passing to a Morita equivalent  groupoid $\gm'\toto M'$ if
necessary,
there is an $S^1$-central extension $R\to   \gm$ with a pseudo-connection whose
 pseudo-curvature equals $\eta+\omega+\Omega$.
\end{them}

In conclusion, 
for  a proper Lie groupoid $\gm\toto M$,
if   $H^3 (\Gamma\upcom,\zz)$ has no torsion, then
$H^3 (\Gamma\upcom,\zz)\to H^3 (\Gamma\upcom,\rr)$
is injective by the universal coefficient theorem.
Hence any  integer class in $H^3_{DR} (\Gamma\upcom )$ can
be represented uniquely by an $S^1$-gerbe over
$\XX_\Gamma$ and vice-versa. In this case, one
can define $K$-theory twisted by such a class $[\eta+\omega +\Omega]$.
However, in general, our  twisted $K$-theory
is  only defined for twisting a  class in
 $H^2 (\Gamma\upcom , {\mathcal{S}}^1 )$
not for an integer 3rd De-Rham class $[\eta+\omega +\Omega]$.
This is an essential difference when dealing
with general groupoids.

Let us end this subsection by some examples,  which have
been studied extensively in the literature.

\begin{numex}
\label{ex:gerbe}
\begin{enumerate}
\item Let $M$ be  a manifold and $\alpha \in H^3 (M, \zz )$, and
let  $\{U_i \}$ be  a  good  covering of $M$. Then  the groupoid
$\coprod_{ij} U_{ij}\toto \coprod_i U_i$, where $U_{ij}
\equal  U_i \cap U_j $ is Morita equivalent to $M\toto M$. See Example 
\ref{ex:coh}. Then the $S^1$-gerbe  corresponding to 
the class $\alpha$ can be  realized as an $S^1$-central
extension of groupoids $\coprod_{ij} R_{ij}\to \coprod_{ij} U_{ij}\toto \coprod_i U_i$, where $R_{ij}$ are  $S^1$-bundles over $U_{ij}$, and
the groupoid multiplication is defined as follows: taking  a trivialization
$R_{ij}\cong U_{ij}\times S^1$, then
\begin{equation}
\label{eq:ex1}
(x_{ij}, \lambda_1) (x_{jk}, \lambda_2 )\equal (x_{ik}, \lambda_1  \lambda_2 c_{ijk}),
\end{equation}
where $x_{ij}, \ x_{jk}, x_{ik}$ are the same point $x$ in the
three-intersection $U_{ijk}$ considered as elements in the
two-intersections, and $c_{ijk}: U_{ijk}\to S^1$
 is a  2-cocycle which represents the
Cech class in $H^2 (M, {\mathcal{S}}^1 )$ corresponding to $\alpha$.
Note that $c_{ijk}$ can also be considered as an $S^1$-valued groupoid
2-cocycle of the groupoid $\coprod_{ij} U_{ij}\toto \coprod_i U_i$,
and Eq.  (\ref{eq:ex1}) above is a special case of Eq. (\ref{eq:extension}).
See \cite{brylinski1, Hitchin} for details. 

\item  Let  $\Gamma$ be  a transformation groupoid $G\times M\toto M$,
where $G$ acts on $M$ properly. By Corollary  \ref{cor:eg},
we have $H^3_G (M, \zz)\cong H^2 (\Gamma\upcom, {\mathcal{S}}^1 )$.
Assume that there exists a $G$-invariant good cover  $\{U_i\}$,
then $\gm\toto M$ is Morita equivalent to $\coprod_{ij} G\times U_{ij}\toto 
\coprod_i U_i$, where  the groupoid structure is given by
$s(g, x_{ij})\equal x_{j}, \ \ t(g, x_{ij})\equal gx_{i}$, and 
$$(g, x_{ij} )\cdot (h, y_{jk})\equal (gh, z_{ik})$$
where $x\equal h y$ and $y\equal z$.     
Then the $S^1$-gerbe  corresponding to
the class $\alpha$ can be  realized as an $S^1$-central
extension of groupoids $S^1\to
\coprod_{ij} R_{ij}\to \coprod_{ij}  G\times U_{ij}
\toto \coprod_i U_i$, where $R_{ij}$ are  $S^1$-bundles over
 $G\times U_{ij}$.
For all $i$, $j$, take an open cover $(V_\alpha)_{\alpha\in I_{ij}}$
of  $G$
such that the restriction $R_{ij\alpha}$ of $R_{ij}$ over $V_\alpha
\times U_{ij}$ is isomorphic to the trivial bundle
$V_\alpha\times U_{ij}\times S^1$. The product
$$R_{ij\alpha}\times_{U_j} R_{jk\beta}\to R_{ik\gamma}$$
has the form
\begin{equation}\label{eqn:product}
(i,j,\alpha,g,x,\lambda)(j,k,\beta,h,y,\mu)\equal (i,k,\gamma,gh,y,\lambda
\mu c_{ijk;\alpha\beta,\gamma}(g,x,h,y)),
\end{equation}
where $c_{ijk;\alpha\beta,\gamma}\colon \{(g,x,h,y)\in V_\alpha\times
U_{ij}\times V_\beta\times U_{jk}\vert\; x\equal hy,\,gh\in V_\gamma\}\to S^1$
satisfies the following cocycle relation which expresses that
the product is associative:
\begin{eqnarray*}
\lefteqn{c_{ijk;\alpha_1\alpha_2,\alpha_{12}}(g_1,x,g_2,y)
c_{ikl;\alpha_{12}\alpha_3,\alpha_{123}}(g_1g_2,y,g_3,z)}\\
&\equal &c_{jkl;\alpha_2\alpha_3;\alpha_{23}}(g_2,y,g_3,z)
c_{ijl;\alpha_1\alpha_{23},\alpha_{123}}(g_1,x,g_2g_3,z).
\end{eqnarray*}
Conversely,  given a cocycle as above,
then one can associate to it an $S^1$-central extension
$$S^1\to R\to \coprod_{i,j}G\times U_{ij}\toto \coprod_iU_i.$$
The proof is elementary but tedious. We omit it here.
\end{enumerate}
\end{numex}

\begin{rmk}
There is a canonical map $H^3_G (M, \zz)\to H^3(M, \zz)$
induced by the inclusion   of $M$ to the unit space
of  $G\times M\toto M$.
This implies that an equivariant gerbe should
induce a gerbe over $M$. From the picture of $S^1$-central extensions,
such a gerbe over $M$ is simply the restriction of 
the $S^1$-central extension $R' \to \gm'\toto M'$ to the unit space,
where $R' \to \gm'\toto M'$ is an $S^1$-central extension representing
this equivariant gerbe. In some cases, we have an isomorphism
$H^3_G (M, \zz)\cong  H^3(M, \zz)$. It is interesting to investigate 
how an  $S^1$-gerbe over $M$ can be made an equivariant one under
this assumption. 
\end{rmk}

\subsection{Continuous case}

The purpose of this subsection is to clarify the relation with \cite{kum98}.
To relate their  constructions  to ours, let ${\mathcal{S}}_{cont}^1$
be the sheaf of continuous $S^1$-valued functions. We need
to determine whether the natural map
$H^2(\Gamma\upcom,{\mathcal{S}}^1)\to H^2(\Gamma\upcom,
{\mathcal{S}}^1_{cont})$ is an isomorphism. Unfortunately,
we don't know the answer in general, but we can prove that it is
an isomorphism in our main case of interest:

\begin{prop}\label{prop:smooth cocycle}
Let $\Gamma$ be a proper Lie groupoid. Then
the natural map
\begin{equation}\label{eqn:smooth cocycle}
H^2(\Gamma\upcom,{\mathcal{S}}^1)\to
H^2(\Gamma\upcom,{\mathcal{S}}^1_{cont})
\end{equation}
is an isomorphism.
\end{prop}

\begin{pf}
Recall from Proposition~\ref{prop:H2 equal H3} that $H^2(\Gamma\upcom,
{\mathcal{S}}^1)$ is isomorphic to $H^3(\Gamma\upcom,\zz)$.
We claim that $H^2(\Gamma\upcom,{\mathcal{S}}^1_{cont})$
is also isomorphic to $H^3(\Gamma\upcom,\zz)$. Indeed,
Crainic's proof that smooth groupoid cohomology
vanishes \cite{C} also works for continuous cohomology
since Crainic only uses integration and cutoff functions, and
never uses differentiation.
\end{pf}

Exactly the same constructions can be performed in the category
of locally compact groupoids: let us denote by $Tw^{lt}(\Gamma)$,
$\eE^{lt}(\Gamma)$ and ${\mbox{Ext}}^{lt}(\Gamma,S^1)$ the groups
thus obtained. The superscript ``lt'' stands for ``locally trivial'',
since central extensions are required to be locally trivial
$S^1$-principal bundles (in the continuous sense), and Morita
equivalences between groupoids are required to be locally trivial
principal bundles. An immediate consequence of
Proposition~\ref{prop:smooth cocycle} is the following

\begin{cor}
Let $\Gamma$ be a proper Lie groupoid. Then the natural map
$$\mbox{Ext}^{sm}(\Gamma,S^1)\to \mbox{Ext}^{lt}(\Gamma,S^1)$$
is an isomorphism.
\end{cor}
\par\medskip

However, in \cite{kum98}, $S^1$-central extensions $S^1\to R\to\Gamma$
are not required to be locally trivial: the homomorphism
$R\to \Gamma$ is only required to be open surjective. Moreover,
the notion of Morita equivalence in \cite{kum98} is weaker
since in their definition of an equivalence bimodule
$\Gamma_1^{(0)}\stackrel{\tau}{\leftarrow}
Z\stackrel{\sigma}{\to} \Gamma_2^{(0)}$,
the maps $\sigma$ and $\tau$ are just open surjective, and
the actions of $\Gamma_1$ and $\Gamma_2$ on $Z$ are free and proper,
but $Z$ is not necessarily a locally trivial $\Gamma_i$-principal bundle.
Let us denote by $Tw^{lc}(\Gamma)$, $\eE^{lc}(\Gamma)$
and $\mbox{Ext}^{lc}(\Gamma,S^1)$ the groups constructed in
\cite{kum98}. There are obvious natural morphisms
$$\xymatrix{
Tw^{sm}(\Gamma)\ar[r]\ar[d]
     & Tw^{lt}(\Gamma)\ar[r]\ar[d]
     & Tw^{lc}(\Gamma)\ar[d]\\
\eE^{sm}(\Gamma)\ar[r]\ar[d]
     & \eE^{lt}(\Gamma)\ar[r]\ar[d]
     & \eE^{lc}(\Gamma)\ar[d]\\
{Ext}^{sm}(\Gamma,S^1)\ar[r]
     & {Ext}^{lt}(\Gamma,S^1)\ar[r]
     & {Ext}^{lc}(\Gamma,S^1).
}$$
Since any $S^1$-central extension of Lie groupoids is the pull-back
of the central extension $$S^1\to U(\HH)\to PU(\HH)$$ which
is locally trivial (Section~\ref{sec:DD}), the map
$Tw^{lt}(\Gamma)\to Tw^{lc}(\Gamma)$ is an isomorphism. Therefore
$\eE^{lt}(\Gamma)\to \eE^{lc}(\Gamma)$ and
$\mathrm{Ext}^{lt}(\Gamma,S^1)\to {\mathrm{Ext}}^{lc}(\Gamma,S^1)$
are surjective.

From Proposition~\ref{prop:trivial extensions} and its analogue
for $\eE^{lt}$ and $\eE^{lc}$ instead of $\eE^{sm}$, an extension
is zero in $\eE^{lt}(\Gamma)$ if and only if there exists an
open cover $(U_i)$ such that its class is zero in $Tw^{lt}(\Gamma
[U_i])$, and similarly for $\eE^{lc}$. Therefore,

\begin{prop}\label{prop:sm equal lc}
Let $\Gamma$ be a Lie groupoid. Then
\begin{itemize}
\item[(a)] the natural maps $\eE^{lt}(\Gamma)\to\eE^{lc}(\Gamma)$
and $Ext^{lt}(\Gamma, S^1)\to Ext^{lc}(\Gamma, S^1)$ are isomorphisms.
\item[(b)] If $\Gamma$ is proper, then
$Ext^{sm}(\Gamma, S^1)\stackrel{\sim}{\to} Ext^{lt}(\Gamma,S^1)
\stackrel{\sim}{\to} Ext^{lc}(\Gamma,S^1)$.
\end{itemize}
\end{prop}

In \cite{kum98} is defined the Brauer group $Br(\Gamma)$ of $\Gamma$.
It is the group
of locally trivial bundles of $C^*$-algebras over $M$ endowed with an
action of $\Gamma$,
with fibers isomorphic to ${\mathcal{K}}$, divided by
Morita equivalence. Let $Br_0(\Gamma)$ be the subgroup of
$Br(\Gamma)$ consisting of those bundles whose Dixmier-Douady
class in $H^3(M,\zz)$ is zero. Then $Br_0(\Gamma)$ is
the group of bundles of the form $M\times {\mathcal{K}} $
with the diagonal action
$\gamma\cdot (s(\gamma),T)\equal (t(\gamma),\pi(\gamma)(T))$,
where
\begin{equation}\label{eqn:projective rep}
\pi\colon \Gamma\to \mbox{Aut}({\mathcal{K}} )\cong
PU(\HH)
\end{equation}
is a ``projective representation'' of $\Gamma$.
The group structure is given by tensor product:
$[\pi][\pi']\equal [\pi\otimes\pi']$,
 where $(\pi\otimes\pi' )(\gamma)
\in \mbox{Aut}({\mathcal{K}}(\HH\otimes\HH))\cong
\mbox{Aut}({\mathcal{K}} )$.
\par\medskip

Recall \cite{kum98} that
\begin{equation}\label{eqn:BR0}
Br_0(\Gamma)\cong \eE^{lc}(\Gamma).
\end{equation}
Indeed,
from the data $(M\times {\mathcal{K}} \to M ,\pi )$,
one obtains an $S^1$-central extension as follows:
$$S^1\to
\{(\gamma,U)\in\Gamma\times U(\HH)|\;
\pi(\gamma)\equal \mbox{Ad}(U)\}\to
\Gamma.$$
For the construction of a bundle of $C^*$-algebras obtained from a
central extension, see \cite{kum98} or Section~\ref{sec:DD}.
\par\medskip

If $\{U_i\}$ is a cover of $M$ by contractible open subspaces
and if $\Gamma'$ denotes $\Gamma[U_i]$, then
$\mbox{Ext}^{lc}(\Gamma,S^1)\cong
Br(\Gamma)\cong Br(\Gamma')\cong {\mathcal{E}}^{lc}(\Gamma')$.
To summarize,

\begin{prop}
If $\gm$ is a proper Lie groupoid, then we have
 $$Br(\gm)\cong
{\mathrm{Ext}}^{sm}(\gm,S^1)\cong H^2(\gm\upcom,{\mathcal{S}}^1)
\cong H^3(\gm\upcom,\zz). $$
\end{prop}

\subsection{$S^1$-gerbes via principal
 $G$-bundles over  groupoids}

The purpose of this subsection is to present
 another construction of $S^1$-gerbes using principal 
$G$-bundles over  groupoids
together with an $S^1$-central extension of $G$.
In fact, we show, in the next subsection,
 that every $S^1$-gerbe arises in this
way when $G$ is taken the projective  unitary group $PU(\HH)$
of a separable Hilbert space $\HH$.
Let us recall the definition of principal $G$-bundles.

\begin{defn}
\label{def:principal}
Let $\Gamma\rightrightarrows M$ be a Lie groupoid.       
A $\gm$-space consists of a smooth manifold $P$ together with
a smooth map $J: P\to M$ such that

(i) there is a map $\sigma : Q\to P$, where $Q$ is the fibered
product $Q\equal \gm \times_{s, M, J} P$.
We write $\sigma(\gamma, x)\equal \gamma \cdot x$.

\noindent This map is subject to the constraints

(ii) for all $x\in P$ we have
$$J(x) \cdot x\equal x$$

(iii)  for all $x\in P$ and all $\gamma,\delta\in\Gamma$ such that
$J(x)\equal s(\gamma)$ and $t(\gamma)\equal s(\delta)$ we have
$$(\delta \cdot \gamma )\cdot x\equal  \delta \cdot (\gamma  \cdot x).$$
\end{defn}

Note that, as a consequence of the above definition,
  for any $r\in \gm$, the map
\begin{equation}
\label{eq:lr}
l_{r}: J^{-1}(u)\to J^{-1}(v), \ \ \ x\to r\cdot x
\end{equation}
     must be  a diffeomorphism, where $u\equal s(r)$ and $v\equal t(r)$.

Associated to any $\gm$-space $J: P\to M$, there is a
natural groupoid $Q\toto P$, called the 
{\em transformation groupoid},
which is defined as follows  $Q\equal \gm \times_{s, M,J} P$,
the source and target maps are, respectively,
$s(\gamma , x)\equal x$, $t(\gamma , x)\equal \gamma \cdot x$,  and the multiplication
\begin{equation}
\label{eq:transformation}
(\gamma ,y )\cdot (\delta , x)\equal (\gamma \cdot\delta , x), \ \ \ \ \mbox{where }
y\equal \delta \cdot x.
\end{equation}
It is simple to check that  the first projection
defines a (strict) homomorphism of groupoids from $Q\toto P$
to $\gm \toto M$.

\begin{defn}\label{def:G bundle}
A  principal $G$-bundle over $\gm\toto M$ is a principal  right $G$-bundle
$P\stackrel{J}{\to} M$, which, at the same time, is also a $\gm$-space
 such that  the following compatibility condition is satisfied:
for all $x\in P$ and $\gamma \in \gm$, $s(\gamma )\equal J(x)$
\begin{equation}
(\gamma \cdot x) \cdot g\equal \gamma \cdot (x \cdot g ) .
\end{equation}
\end{defn}
In this case  $Q \to \gm $ also  becomes a principal (right)
$G$-bundle.

\begin{numex}
\label{ex:equivariant}
Let $\gm$ be the  transformation groupoid $H\times M\toto M$.
Then a principal $G$-bundle over $\gm$ corresponds exactly
to an $H$-equivariant principal (right) $G$-bundle over $M$.
\end{numex}

A principal $G$-bundle over a groupoid
$\gm \toto M$  can also be equivalently considered as a generalized
homomorphism from  $\gm \toto M$ to $G\toto \cdot$.
As a consequence of
Proposition~\ref{prop:composition homomorphisms},
we see that
principal bundles behave well under the ``generalized homomorphisms'' in
the following sense.

\begin{prop}\label{prop:pullback principal bundle}
Let $f$ be a generalized homomorphism from  $ \Gamma_1 \toto M_1$
to $ \Gamma_2 \toto M_2$ given by
$$M_1 \stackrel{\tau}{\leftarrow} X
\stackrel{\sigma}{\rightarrow} M_2 .$$
 Then for any principal $G$-bundle
$P\to M_2$ over $ \Gamma_2 \toto M_2$,
$$f^*P\stackrel{def}{\equal }X\times_{M_2} P\to M_1$$
is  a principal $G$-bundle over $ \Gamma_1 \toto M_1$. As a
consequence, if $ \Gamma_1 \toto M_1$  and $ \Gamma_2 \toto M_2$
are Morita equivalent groupoids, then there is a bijection
between their principal $G$-bundles.
\end{prop}

%
%
%
%
%

Given a principal $G$-bundle $J: P\to M$ 
over $\gm \toto M$, let
 $\frac{P \times P}{G}\toto M$ be the  gauge  groupoid.
 We denote by $(p_1, p_2)$ an element of $P \times P$
and by $\overline{(p_1, p_2 )}$
the class of this element in $\frac{P \times P}{G}$.
A map from $\Gamma$ to $ \frac{P \times P}{G}$ is defined by $\gamma \mapsto
\overline{(\gamma p, p)}$ where $p$ is any element that satisfies
$J(p)\equal s(\gamma)$. Thus we obtain the following groupoid
homomorphism:

\begin{equation}
\label{eq:Q1}
\xymatrix{
\Gamma \ar@<-.5ex>[d]\ar@<.5ex>[d]\ar[r] &
\frac{P \times P}{G} \ar@<-.5ex>[d]\ar@<.5ex>[d]\\
M\ar[r] & M}
\end{equation}

Since any transitive groupoid is Morita equivalent  to its isotropy group,
$\frac{P\times P}{G}\toto M $ is Morita equivalent to $G \toto \cdot$.
It is not hard to check that the homomorphism (\ref{eq:Q1}) and  
the $G$-principal bundle $P$ define  the isomorphic
 generalized homomorphisms from $\gm$ to $G$.


%


From  Proposition~\ref{prop:functor Morita},
it follows that a generalized homomorphism $f$ from  $ \Gamma_1 \toto M_1$
to $ \Gamma_2 \toto M_2$ induces a natural 
 homomorphism, called the {\em pull back map}: 
$$f^*: H^2(\Gamma_2\upcom, {\mathcal{S}}^1 )
\lon H^2(\Gamma_1\upcom, {\mathcal{S}}^1 ). $$

In what follows, we describe a construction of gerbes over
a  stack which is similar to the construction in \cite {brylinski1, BCMMS}.
 Assume that $S^1 \to \tilde{G}\to G$ is an
$S^1$-central extension of Lie groups, $P$ is a $G$-bundle over $\gm\toto M$.
Then  $P$ defines a  generalized homomorphism from $\gm\toto M$
to $G\toto \cdot$, and therefore induces 
a pull back map $H^2 (G\upcom,  {\mathcal{S}}^1 )\to H^2 (\gm \upcom,  {\mathcal{S}}^1 )$.
 By pulling back the class of $S^1 \to \tilde{G}\to G$
in $H^2 (G\upcom,  {\mathcal{S}}^1 )$ via this
map, one  obtains
 an element in $H^2(\Gamma\upcom, {\mathcal{S}}^1 )$, i.e., an 
$S^1$-gerbe over the stack $\XX_\gm$
associated to $\Gamma$.


Since  the equivalence classes
of $G$-bundles $P\stackrel{J}{\to} M$ over $\gm\toto M$
are classified by  $H^1 (\gm\upcom , G)$,  
we have a  map

\begin{equation}
\label{eq:Phi2}
\Phi : H^1 (\gm\upcom , G)\times H^2 (G\upcom , {\mathcal{S}}^1 )
\lon  H^2 (\gm\upcom , {\mathcal{S}}^1 ).
\end{equation}

Below we describe an explicit construction of the map 
$\Phi$ in a special case  more relevant to
us.

Besides the above assumption, we furthermore assume that,
as  a $G$-principal bundle,  $P\to M$ can be lifted to
a $\tilde{G}$-principal bundle $\tilde{P}\to M$.
Note that if $\tilde{P}\to M$ is a principal $\tilde{G}$-bundle,
there is  a natural $S^1$-action on $\tilde{P}$ defined
as follows: $\forall \lambda\in S^1 , \tp \in
\tilde{P}, \  \lambda \cdot \tp \equal  (\lambda \cdot 1_{\tilde{G}}) \tp$,
where $\lambda \cdot 1_{\tilde{G}}$ is considered as
an element in $\tilde{G}$.
Then $\tilde{P}/S^1$ is a principal $G$-bundle over $M$, which is
isomorphic to the reduced principal bundle induced by the group homomorphism
$\tilde{G}\to G$. We require that as a principal
$G$ bundle $P\cong \tilde{P}/S^1 $.
In this case, it is simple to see that 
$$\frac{\tilde{P} \times \tilde{P}}{ \tilde{G}} \to 
\frac{ P\times P}{G}\toto M$$ is an $S^1$-central
extension, which is Morita equivalent to
$$\tilde{G}\to G\toto \cdot$$
Here the $S^1$-equivalent Morita equivalence bimodule  
$M\stackrel{\tilde{\tau}}{\leftarrow} \tilde{P}
\stackrel{\tilde{\sigma}}{\rightarrow} \cdot $
 is  given by the composition of the projection
$\tilde{P}\to P$  with
$M\stackrel{\tau}{\leftarrow} P
\stackrel{\sigma}{\rightarrow} \cdot $,
the left action of $\frac{\tilde{P} \times \tilde{P}}{ \tilde{G}}$
 on $\tilde{P}$
is  $[(\tilde{p}_1,\tilde{p}_2)]\cdot \tilde{p}_3\equal \tilde{p}_1
\tilde{g}$ where $\tilde{g}$ is the unique element in $\tilde{G}$
such that $\tilde{p}_3\equal \tilde{p}_2\tilde{g}$, and the
right $\tilde{G}$-action is the usual one.

Let $R\to \gm$ denote the pull-back $S^1$-bundle of
$\frac{\tilde{P} \times \tilde{P}}{ \tilde{G}} \to
\frac{ P\times P}{G}$ via the map $\gm \to \frac{ P\times P}{G}$
 as in Eq. (\ref{eq:Q1}).

\begin{prop}
\label{pro:pull-back}
Under the same hypothesis as above, 
 $R\to \gm$ is a groupoid $S^1$-central extension,  whose corresponding class
in $H^2 (\gm\upcom , {\mathcal{S}}^1 )$ is equal to
$\Phi (\alpha , \beta)$. 
Here $\alpha\in H^1(\gm\upcom , G)$ is
the class defined by  $P\to M$ and $\beta \in H^2 (G\upcom, {\mathcal{S}}^1 )$
is the class corresponding to  the central extension
$S^1\to \tilde{G}\to G$.
\end{prop}

\subsection{Hilbert bundle and Dixmier-Douady class}
\label{sec:DD}

\label{sec:gerbe-principal}
The purpose of this subsection is to show that 
every $S^1$-gerbe over a differential stack always arises
from a principal  $PU(\HH )$-bundle over the stack
as in the case of manifolds \cite{Atiyah}. However,
unlike the manifold case, such projective bundles may
not be  unique. Nevertheless, we show that 
there always  exists a canonical one.
 In the following, we describe an explicit construction of such a projective
 bundle.

We now  fix a  separable Hilbert space $\HH$ and consider the
canonical  $S^1$-central extension:

\begin{equation}
S^1\to U(\HH )\to PU(\HH ),
\end{equation}
which is a generator of $H^2 (PU(\HH )\upcom, {\mathcal{S}}^1 )$ \cite{brylinski1}.
Thus Eq. (\ref{eq:Phi2}) induces a group homomorphism 

\begin{equation}
\label{eq:Phi1}
\Phi' : H^1 (\gm\upcom , PU(\HH ) )
\to  H^2 (\gm\upcom , {\mathcal{S}}^1 ).
\end{equation}

Note that  $H^1(\gm\upcom,PU(\HH))$ can be
endowed with the following abelian group structure:
$[\pi][\pi']\equal [\pi\otimes\pi']$ if $\pi$, $\pi'\colon
\Gamma[U_i]\to PU(\HH)$ are groupoid homomorphisms and $(U_i)$
is an open cover of $\Gamma$. See \cite{brylinski1} 
for the case when $\gm$ is a manifold $M\toto M$.
 
In other words, any principal $PU(\HH )$-bundle over $\gm \toto M$
defines an element in $H^2 (\gm\upcom , {\mathcal{S}}^1 )$, or an $S^1$-gerbe
over the associated  stack $\XX_\gm$.

When $\gm $ is a manifold $M\toto M$, $\Phi'$ is
indeed an isomorphism \cite{Atiyah, brylinski1}.  However,
in general,  $\Phi'$ may not be
injective.\footnote{For instance, let $G$ be a compact Lie group and
$\pi$ any unitary representation of $G$ such that $\pi(g)$ is not
a scalar multiple of the identity for some $g\in G$. Then the
associated element $[\pi]\in H^1(G\upcom, PU(\HH))$ is nonzero,
but $\Phi'([\pi])\equal 0$ since the composition
$H^1(G\upcom, U(\HH))\to H^1(G\upcom, PU(\HH))
\to H^2(G\upcom, \sS^1)$ is zero.}
We will see below
  that $\Phi'$
 admits a canonical left inverse. Therefore it is
always surjective.

First of all, let us assume that 
$\alpha \in H^2 (\gm\upcom , {\mathcal{S}}^1 )$
is the class defined by a groupoid $S^1$-central extension
$R\to \gm \toto M$.





\begin{defn}\label{defi:equivariant}
A  complex-valued function $f$ on 
$R$  is said to be equivariant if $ f(\lambda \tilde{\gamma})\equal 
 \lambda^{-1} f(\tilde{\gamma}) $ for any $\lambda \in S^1$
     and any $\tilde{\gamma} \in R$. 
\end{defn}

 Let $\lambda\equal (\lambda^x)_{x\in M}$ be a Haar system on $ R $,
i.e. $ \lambda^x$  is a measure on $R^x$ such that
for any $\tilde{\gamma}\in R $ the map $ L_{\tilde{\gamma}} :
 R^{s(\tilde{\gamma})}
 \to R^{t(\tilde{\gamma})}$ defined by $ \tilde{\gamma}'
\mapsto \tilde{\gamma} \tilde{\gamma}' $
preserves the measure.

  By ${\mathcal L}^2_{x}$,  we denote the space
$L^2 (R^x)^{S^1} $ consisting  of  $S^1$-equivariant functions defined 
on $R^x$ which  are $ L^2$ with respect to the Haar measure.
Let
\begin{equation}
\label{eq:Hx}
 {\mathcal H}_x \equal {\mathcal L}^2_{x}\otimes \HH, \ \ \ \mbox{ and } \ \ 
\tilde{\mathcal H}\equal  \coprod_x {\mathcal H}_x .
\end{equation}
Then $\tilde{\mathcal H}\to M$ is a
countably generated continuous field of infinite dimensional
Hilbert spaces over the finite dimensional space $M$,
and therefore is a locally trivial Hilbert bundle (indeed 
globally trivial) according to  Dixmier-Douady theorem \cite{dix-dou63}.

          For $x\in M$, let ${\mathcal B}_x$ be the set of orthonormal basis of
${\mathcal H}_x$ and ${\mathcal{B}}\equal \coprod_{x\in M}{\mathcal{B}}_x$.
We endow ${\mathcal B}$ with the following topology: identify
${\mathcal{B}}_x$ with the space $U({\mathcal{H}}_x,\HH)$ of unitary maps
from ${\mathcal{H}}_x$ to $\HH$. Then a section $x\mapsto u_x$ is
continuous if and only if for every $\xi\in \HH$, $x\mapsto u_x^{-1}\xi$
is a continuous section of the field $\tilde{\mathcal{H}}\to M$.
    The fiber bundle ${\mathcal B} \to M$
 is a principal $U(\HH )$-bundle.
Now $S^1 $ naturally  acts on ${\mathcal B}$  by scalar multiplication.
 Let $P{\mathcal B}\equal {\mathcal B}/S^1 $  be its quotient.
         Then  $ P{\mathcal B}$ is a principal $PU (\HH )$-bundle
over $ M$.

  Let $U(x,y)$ be the set of unitary linear maps from ${\mathcal  H}_x$
  to ${\mathcal H}_y$, and
 $ U({\mathcal H},{\mathcal H})\equal \{U(x,y) |\; (x,y)\in M\times M\} $.
         Then $U({\mathcal H},{\mathcal H})$ is
naturally  a groupoid over $M$.
 
  Let $PU(x,y)$ be the set of unitary projective maps from $ {\mathcal
H}_x$
        to ${\mathcal H}_y$, and
let $ PU({\mathcal H},{\mathcal H})\equal \{PU(x,y) |\; (x,y)\in M\times M \}$.
Then $PU({\mathcal H},{\mathcal H})$ is a groupoid over $M$.

The groupoid $ R $ acts naturally on $ {\mathcal H}$:
 for any element $\tilde{\gamma} \in R $
 with $x\equal s(\tilde{\gamma})$ and  $y\equal  t(\tilde{\gamma})$,
and any equivariant function $f \in {\mathcal H}_x$, the action is given by:

$$ f\mapsto \tilde{\gamma}\cdot f, \qquad\mbox{ where }
(\tilde{\gamma}\cdot f)(r)\equal f(\tilde{\gamma}^{-1} r).$$

Since this action preserves the measure $\lambda$,
it induces a homomorphism of groupoids

$$ i : R \to U({\mathcal H},{\mathcal  H})  .$$

Since $i$ is equivariant under the $S^1$-actions,
it induces a homomorphism  of groupoids $j$:

$$ j: \Gamma \to PU({\mathcal H},{\mathcal H})  . $$

 In short, we have the following diagram  of groupoid homomorphisms:
           \begin{equation} \label{eq:dia}  \xymatrix{
        R \ar[d] \ar[r]^i & U({\mathcal H},{\mathcal H})\ar[d]   \\
        \Gamma \ar[r]^j & PU({\mathcal H},{\mathcal H})
}                      \end{equation}

It is obvious that $P{\mathcal B}\to M$ is a principal 
$PU(\HH)$-bundle over the groupoid $PU({\mathcal H},{\mathcal H})\toto M$. By
pushing forward the action using the above groupoid homomorphism,
$P{\mathcal B}\to M$ is naturally a principal $PU(\HH)$-bundle over
the groupoid $\gm \toto M$.

\begin{prop}
\label{pro:PB}
If $\alpha \in H^1(\gm\upcom , PU(\HH))$ denotes the
class defined by $P{\mathcal B}\to M$, then
$\Phi' (\alpha )$
is equal to the class in $H^2 (\gm\upcom , {\mathcal{S}}^1)$
 corresponding to the $S^1$-central extension $R\to \gm$.
\end{prop}
\begin{pf}
Note that as groupoids, $U({\mathcal H},{\mathcal  H})\toto M$
is isomorphic to $\frac{ {\mathcal B}\times {\mathcal B}}{U(\HH)}\toto M$,
and $PU({\mathcal H},{\mathcal  H})\toto M$
is isomorphic to $\frac{ P{\mathcal B}\times P{\mathcal B}}{PU(\HH)}\toto M$.
Thus, the conclusion follows from Proposition \ref{pro:pull-back}
and Diagram (\ref{eq:dia}).
\end{pf}

Now let us return to the general case.
Consider a  Lie groupoid $\Gamma\toto M$ and
an element $\alpha\in H^2(\Gamma\upcom, {\mathcal S^1})$.
 There exists a Lie  groupoid $\Gamma'\toto M'$ Morita equivalent
 to $\Gamma$ such that $\alpha$ is the class of an $S^1$-central extension
$$S^1\to R\to\Gamma'\toto M'.$$
Let $P{\mathcal{B}}\to M'$ be the
corresponding principal $PU(\HH)$-bundle over $\Gamma'\toto M'$
constructed above as in Proposition \ref{pro:PB}. Since 
$\Gamma'\toto M'$ and $\Gamma\toto M$ are Morita equivalent,
there is an associated principal $PU(\HH)$-bundle
$P_\alpha\to M$ over the groupoid $\Gamma\toto M$.
 In fact, $P_\alpha \equal  (Z\times_{M'}P{\mathcal{B}})/{\Gamma'}$,
where $M\leftarrow Z \to M'$ is an equivalence bimodule between
$\Gamma\toto M$ and $\Gamma'\toto M'$.
By construction, $P_\alpha$ and $P{\mathcal{B}}$ represent
the same generalized homomorphism, thus define
 the same element
in $H^1(\Gamma\upcom,PU(\HH))$.
 Moreover,
$P_\alpha$ does not depend on a particular choice of the $S^1$-central
extension $S^1\to R\to \Gamma'\toto M'$ realizing the class $\alpha$.
This follows from the following
\begin{lem}
Assume that $p\colon Y\to M$ is a surjective submersion. Let
$f\colon \gm[Y]\to \gm$ be the projection map. Assume that $S^1\to
R\to \gm\toto M$ is an $S^1$-central extension, and denote by $P$ (resp. $P'$)
the associated $PU(\HH)$-bundle over $\gm$ (resp. $\gm[Y]$).
Then $P' $ is isomorphic to $ P\smalcirc f$ as generalized morphisms $\gm[Y]\to PU(\HH)$.
\end{lem}

\begin{pf}
Let us first treat the case $Y\equal \amalg_{i\in I}U_i$. Let $I_x
\equal \{i\in I\vert\; x\in U_i\}$ and $\hH_I\equal \amalg_x \ell^2(I_x)$.
Then $\hH_I\to M$ is endowed with a structure of continuous field of
Hilbert spaces over $M$ (associated to the $C_0(M)$-Hilbert module
$\oplus_{i\in I} C_0(U_i)$, see Proposition~\ref{prop: module field}).

It is easy to see that $P'\equal PU(f^*(\tilde{\hH}\otimes\hH_I),\HH)$
(where ``$\otimes$'' denotes the tensor product of continuous fields over $M$).
Now, $\tilde{\hH}\otimes\hH_I\cong (\tilde{\hH}\otimes \HH)\otimes\hH_I
\cong \tilde{\hH}\otimes(\hH_I\otimes \HH)
\cong \tilde{\hH}\otimes\HH \cong \tilde{\hH}$
since $\hH_I\otimes \HH$ is the trivial continuous field
$M\times \HH \to M$ (see the argument below (\ref{eq:Hx})).
It follows that $P'\equal PU(f^*\tilde{\hH},\HH)\equal P\smalcirc f$.
\par\medskip
In the general case, i.e. for a general $Y$, consider a continuous
$p$-system $\mu\equal (\mu_x)_{x\in M}$, i.e. $\mu_x$ is a measure with
support $p^{-1}(x)$ such that
$$\forall \varphi\in C_c(Y),\quad [x\mapsto \int \varphi(y)\,
d\mu_x(y)]\in C_c(M).$$
A Haar system on $\gm[Y]$ is given by
$$\int_{\gm[Y]^y}\psi \equal 
\int_{\gamma\in\gm^{p(y)}} d\lambda^{p(y)}(\gamma)
\int_{z\in p^{-1}(s(\gamma))} d\mu_{s(\gamma)}(y')\psi(y,\gamma,y').  $$
Then, $\amalg_{x\in M} L^2(\mu_x)\to M$ is a continuous field
of Hilbert spaces (associated to the $C_0(M)$-module obtained by
the completion of $C_c(Y)$ with respect to the scalar product
$\langle\varphi,\varphi\rangle (x)
\equal \int |\varphi|^2\,d\mu_x$), such that
$\amalg_{x\in M} L^2(\mu_x)\otimes \HH$ is the trivial field
$M\times \HH\to M$. The proof is almost the same as above, except that
notations are more complicated. We omit details.
\end{pf}

Therefore we have proved the following

\begin{prop}
\label{pro:can}
Let $\Gamma\toto M$ be a Lie groupoid. 
Associated to any element   $\alpha\in H^2(\Gamma\upcom, {\mathcal S^1})$,
there is a canonical $PU(\HH)$-bundle over $\Gamma\toto M$,
 denoted by  $P_\alpha\to M$,  whose
corresponding class in $H^1 (\gm\upcom , PU (\HH ))$ goes to $\alpha$ under
the map $\Phi'$ in Eq. (\ref{eq:Phi1}).
\end{prop}

Clearly, if  $\alpha$ can be realized as an $S^1$-central extension
over the groupoid $\gm\toto M$ without the need of passing
to Morita equivariance, then $P_\alpha\equal P{\mathcal{B}}$.
As a consequence, when $\gm$ is a transformation groupoid,
 we obtain the following:

\begin{cor}
\label{cor:equivariant}
If  $G$ is a Lie  group acting on $M$ properly,
then there is a  group homomorphism:
\begin{equation}
\{\mbox{Isomorphism classes of $G$-equivariant $PU(\HH)$-bundles}\}
\to H^3_G (M, \zz ),
\end{equation}
which admits a canonical inverse. Namely, to any element in $H^3_G (M, \zz )$,
 there  associates a canonical $G$-equivariant $PU(\HH)$-bundle $P_\alpha\to M$.
\end{cor}

\section{Twisted $K$-theory and Fredholm bundles}

In this section, we introduce twisted $K$-theory groups
of  a  Lie groupoid (or more precisely, a differential
stack). In the case of proper Lie groupoids,
we describe these $K$-groups in terms of  homotopy classes
of certain $\Gamma$-invariant sections of Fredholm
operators associated to the projective
Hilbert bundle  as constructed  in Sect.~\ref{sec:DD}
 (Theorem~\ref{thm:fredholm proper}).

\subsection{The reduced $C^*$-algebra of an $S^1$-central extension}
Given an $S^1$-central extension of Lie groupoids
$S^1\to R\to \Gamma\toto M$,
let   $L\equal R\times_{S^1}\cc$ be its associated complex line bundle.
Then $L\to M$ can be considered as   a Fell bundle
of $C^*$-algebras over the groupoid $\Gamma\toto M$.
Therefore one can  construct a $C^*$-algebra
 out of it (see Appendix \ref{subsec:reduced}).

\begin{defn}\label{def:CGR}
Let $\Gamma$ be a Lie groupoid and $S^1\to R\to \Gamma\toto M$
an $S^1$-central extension. Then the reduced $C^*$-algebra
 of the central extension 
$C^*_r(\Gamma;R)$ 
 is defined   to be $C^*_r(\Gamma;L)$,
where $L\equal R\times_{S^1}\cc$ is the associated
complex line bundle considered as  a Fell bundle
of $C^*$-algebras  over $\Gamma \toto M$.
\end{defn}

There is another picture for this $C^*$-algebra.
Consider
$$C_c(R)^{S^1}\equal \{\xi\in C_c(R)|\; \xi(\lambda r)\equal
\lambda^{-1}\xi(r),\;\forall \lambda\in S^1,\;r\in R\}.$$
One easily checks that $C_c(R)^{S^1}$ is stable under both
the convolution and the adjoint, and that the map
\begin{eqnarray}\label{eqn:isoRE}
C_c(R)^{S^1}&\to& C_c(\Gamma;L),
\end{eqnarray}
$\xi\mapsto \eta$, where $\eta(g)\equal [(r,\xi(r))]\in L_g\equal R_g\times_{S^1}\cc$,
is well-defined and is indeed   an isomorphism
of convolution algebras. Let us define
\begin{equation}\label{eqn:CS1}
C^*_r(R)^{S^1}:\equal \overline{C_c(R)^{S^1}}\subset C^*_r(R),
\end{equation}
i.e. $C^*_r(R)^{S^1}$ is the norm-closure of $C_c(R)^{S^1}$ in $C^*_r(R)$
(see reference~\cite{ren80} for details on the construction of
the reduced $C^*$-algebra $C^*_r(\Gamma)$ of a groupoid $\Gamma$).

%

The algebra $C^*(S^1)\equal C^*_r(S^1)$ acts on $C^*_r(R)$
 by convolution operators. More precisely, there is a *-homomorphism

$$\Lambda\colon C^*(S^1)\to M(C^*_r(R))$$
such that for every $f\in C(S^1)$ and every $\xi\in C_c(R)$,
$$(\Lambda(f)\xi)(r)\equal 
\int_{S^1} f(\lambda)\xi(\lambda^{-1}r)\,d\lambda$$
where $d\lambda$ is the normalized Haar measure $\frac{d\theta}{2\pi}$
on $S^1$. Indeed, one only needs to check that 
$$(U_\lambda\xi)(r) \equal \xi(\lambda^{-1}r)$$
defines a unitary representation of $S^1$ into the unitary group
of $M(C^*_r(R))$.

The map $\Lambda$ is non-degenerate, for if $f_n$ is a sequence
in $C(S^1)$ converging to the delta function at 1,
then $\Lambda(f_n)a$ converges to $a$ for all $a\in C^*_r(R)$.
That is, $\Lambda(f_n)$ converges strictly to the identity. Therefore,
$\Lambda$ extends to a unital strictly continuous $*$-homomorphism
$M(C^*(S^1))\to M(C^*_r(R))$ \cite[paragraphs 3.12.10 and 3.12.12]{ped79}.

Let $P_n\in C_r^*(S^1)$ be the convolution by $z^n$, i.e. $P_n$
corresponds to the characteristic function of $\{n\}$ via
the Fourier transformation $C^*_r(S^1)\cong C_0(\zz)$. Let $Q_n\equal \Lambda(P_n)$.
Then the $Q_n$'s are pairwise orthogonal projections.
Since $\Lambda$ is non-degenerate, the sum $\sum Q_n$ is strictly
convergent to 1. Moreover, since $U_\lambda$ is in the center of
$M(C_r^*(R))$, the projections $Q_n$ also belong to the center
of $M(C^*_r(R))$.

Using the formula
$Q_n(\xi)(r)\equal \int_{S^1}\lambda^{n}\xi (\lambda^{-n} r)\,d\lambda$,
 one easily checks that the image of
$Q_n$ is the closure of the set of elements $\xi\in C_c(R)$
such that $\xi(\lambda r)\equal \lambda^n\xi(r)$ for all $(\lambda,r)\in
S^1\times R$. In particular, $C^*_r(R)^{S^1}$,  the closure
of $C_c(R)^{S^1}$ in $C^*_r(R)$,  is $Q_{-1}(C^*_r(R) )$.

Similarly as in Eq. (\ref{eqn:isoRE}), there is an isometric isomorphism
of Hilbert $C_0 (M)$-modules:
\begin{equation}\label{eqn:iso H}
L^2(R)^{S^1}\to L^2(\Gamma;L).
\end{equation}

If $\xi\in C_c(R)^{S^1}$ and $\eta(g)\equal [(r,\xi(r))]$, then
the norm of $\xi$,  as a convolution operator acting on
$L^2(R)$,  is equal to the norm of $\eta$,  as a convolution operator
acting on $L^2(\Gamma;L)$.

Noting that $L^2(R)^{S^1}$ is the image of the projection $Q_{-1}$,
we have

\be
&&\|\xi\|_{C^*_r(R)}\equal \|Q_{-1}\xi \|_{C^*_r(R)}\equal \|\xi Q_{-1}\|_{C^*_r(R)}
\equal \sup_{\|\varphi\|_{L^2(R)}\equal 1} \|\xi*Q_{-1}\varphi\|\\
&&\ \ \ \ \ \
\equal \sup_{\varphi\in L^2(R)^{S^1},\;\|\varphi\|_{L^2(R)}\equal 1} \|\xi*\varphi\|
\equal \sup_{\|\psi\|_{L^2(\Gamma;L)}\equal 1}\|\eta*\psi\|_{L^2(\Gamma;L)}
\equal \|\eta\|_{C^*_r(\Gamma;L)}.
\ee

It follows that
\begin{equation}
C^*_r(\Gamma;R)\cong C^*_r(R)^{S^1}.
\end{equation}

We summarize the above discussion in the following:

\begin{prop}
\label{prop:direct sum}
Let $S^1\to R\to\Gamma$ be an $S^1$-central extension of Lie groupoids.
Then there is a
canonical isomorphism
$$C^*_r(R)\cong\oplus_{n\in\zz} C^*_r(\Gamma;R^n), $$
where $C^*_r(\Gamma;R^n)$ is the $C^*$-algebra of the central extension
$$S^1\to R^n\equal R\otimes\cdots\otimes R\to\Gamma$$
for all $n\ne 0$, and $C^*_r(\Gamma;R^0)\equal C^*_r(\Gamma)$ by convention.

The image of $C^*_r(\Gamma;R)$ in $C^*_r(R)$ consists of the
closure of $C_c(\Gamma,R)^{S^1}$ defined in Eq.~(\ref{eqn:CS1}).

For $f\in C_c(R)\subset C^*_r(R)$, the image $f_n$
of $f$ in $C^*_r(\Gamma;R^n)$
is given by
$$f_n(r)\equal \int_{S^1}\lambda^{-n}f(\lambda^n r)\,d\lambda, $$
where $d\lambda$ is the normalized Haar measure on $S^1$.
\end{prop}

For the  $S^1$-central extension 
 $\coprod_{ij} R_{ij}\to \coprod_{ij} U_{ij}\toto \coprod_i U_i$
  in Example \ref{ex:gerbe} (1), we refer to \cite{RT, rosenberg1, kum96}
for a detailed discussion on the  $C^*$-algebra $C^*_r (\gm , R)$.

\subsection{Definition of twisted $K$-theory and first properties}
\begin{prop}\label{prop:morita invariance}
Let $R_i\to\Gamma_i\toto M_i$ ($i\equal 1,2$) be Morita equivalent 
$S^1$-central extensions. Then $C^*_r(R_1)^{S^1}$ and $C^*_r(R_2)^{S^1}$
are Morita equivalent $C^*$-algebras.
\end{prop}

\begin{pf}
This follows from \cite{ren87} or from
\cite[Theorem 11]{muh01}.
\end{pf}

We are now ready to define twisted $K$-theory.

\begin{defn}\label{def:twisted K theory}
Let $\Gamma$ be a Lie groupoid and
$\alpha\in H^2(\Gamma\upcom ,{\mathcal{S}}^1)$.
We define the twisted $K$-theory as
$$K^i_\alpha(\Gamma\upcom )\equal K_{-i}(C^*_r(R)^{S^1} ),  $$
where $S^1\to R\to \Gamma'\toto M'$ is any central extension
realizing the class $\alpha$ and $\Gamma'$ is Morita equivalent
to $\Gamma$.
\end{defn}

From Proposition~\ref{prop:morita invariance}, it  follows that
if two $S^1$-central extensions are Morita equivalent,  their
twisted $K$-theory groups are isomorphic and therefore only
depend on the corresponding stack and the $S^1$-gerbe over the
stack. Consequently, twisted $K$-theory is well-defined.  

\begin{numex}
(1). When $\gm$ is a manifold $M\toto M$  and  
$\alpha \in H^3 (M, \zz)\cong H^2 (\Gamma\upcom ,{\mathcal{S}}^1)$,
the above definition reduces to the one
introduced by Rosenberg \cite{rosenberg1}.
\par\medskip            

(2).  Assume that a    Lie group $G$ acts  on a smooth manifold
 $M$ properly. According to Corollary \ref{cor:eg},
the equivariant cohomology  $H^3_G (M, \zz)$ is isomorphic
to $H^2 (\gm\upcom , {\mathcal{S}}^1)$, where $\gm$ denotes the transformation
groupoid  $G\times M\toto M$. Let $\alpha \in H^3_G (M, \zz)$.
We define the twisted equivariant $K$-theory
$$K^i_{G, \alpha} (M): \equal K_{-i}(C^*_r(R)^{S^1} ) , $$
where $S^1\to R\to \Gamma'\toto M'$ is any  $S^1$-central extension
realizing the class $\alpha$ and $\Gamma'$ is Morita equivalent
to $\Gamma$. According to   the observation following
Definition \ref{def:twisted K theory}, 
  we have    the following

\begin{prop}
If $G$ acts  on a smooth manifold $M$ properly and freely so that $M/G$
 is a manifold, then 
$$K^i_{G, \alpha} (M)\cong K^i_{\alpha'} (M/G),$$
where $\alpha' $ is the image of $\alpha$ under the
isomorphism  $H^3_G (M, \zz )\stackrel{\sim}{\to} H^3 (M/G, \zz)$.
More generally, if $H$ is a normal subgroup of $G$ which acts on
$M$ properly and freely, then 
$$ K^i_{G, \alpha} (M)\cong K^i_{G/H, \alpha'} (M/H), $$
where $\alpha' $ is the image of $\alpha$ under the isomorphism
  $H^3_G (M, \zz )\stackrel{\sim}{\to} H^3_{G/H} (M/H, \zz)$. 
\end{prop} 

Note that the proposition above  is a non-trivial theorem even in the
non-twisted case, i.e., $\alpha \equal 0$, in the  ordinary
equivariant $K$-theory of Segal \cite{Segal}.
The advantage of our approach is that these facts  are encoded
as a part of the definition since they
are obvious consequences of the Morita equivalence between  the
 transformation groupoids $G\times M\toto M$ and
$G/H\times M/H\toto M/H$. The hard part
is to prove that this definition coincides
with the topological one which is more often used by geometers.
\par\medskip

(3). Given an orbifold ${\mathbf{X}}$, let
  $\gm\toto M$ be  an \'{e}tale groupoid  representing this
orbifold. Now given $\alpha \in  H^3 (\mathbf{X} , \zz ) \cong
H^2 (\gm\upcom ,  {\mathcal{S}}^1 ) $, we define the twisted
orbifold $K$-theory
$$K^i_{\alpha} ({\mathbf{X}}): \equal K_{-i}(C^*_r(R)^{S^1} ),  $$
where $S^1\to R\to \Gamma'\toto M'$ is any $S^1$-central extension
realizing the class $\alpha$ and $\Gamma'$ is Morita equivalent
to $\Gamma$.  It would be interesting to investigate
the relation between our definition with the
one given by Lupercio and Uribe \cite{LU}.
\end{numex} 

Next let  us deduce some properties that are immediate from the definition.

\begin{prop}\label{prop:Bott}[Bott periodicity]
Let $S^1\to R\to\Gamma\toto M$ be an
$S^1$-central extension of Lie
groupoids. Then for all $i$,
\begin{eqnarray*}
K^i_\alpha(\Gamma\upcom )&\cong& K^{i+2}_\alpha(\Gamma\upcom ),\\
K^{i+n}_\alpha(\Gamma\upcom )&\cong& K^i_{\alpha_n}((\Gamma \times\rr^n )\upcom ), 
\end{eqnarray*}
where $\alpha_n$ is the class of the extension $R\times \rr^n
\to\Gamma\times \rr^n\toto M\times \rr^n$.
\end{prop}

Note also that $K^i_{\alpha_n} ((\Gamma \times\rr^n )\upcom )$
is the kernel of the morphism $K^i_{\alpha'_n}((\Gamma\times S^n)\upcom )
\to K^i_{\alpha}(\Gamma\upcom)$ induced by the inclusion
$\Gamma\times\{pt\}\subset \Gamma\times S^n$,
where $\alpha'_n$ is the class
of the extension $R\times S^n\to\Gamma\times S^n\toto M\times S^n$.

We say that a subgroupoid $\Gamma_1 \toto M_1$
of $\Gamma\toto M$ is saturated if  $M_1$
is an invariant subset of $M$ (i.e. $\Gamma_{M_1}^{M_1}\equal \Gamma_{M_1}$)
such that  $\Gamma_1\equal \Gamma_{M_1}^{M_1}$.

\begin{prop}
Let $S^1\to R\to\Gamma\toto M$ be an $S^1$-central extension of Lie
groupoids and denote by $\alpha$ its class in
$H^2(\Gamma\upcom ,{\mathcal{S}}^1)$.
Suppose that $\Gamma_1$ is an open saturated subgroupoid
of $\Gamma$ and let $\alpha_1$ be the class of the corresponding
$S^1$-central extension of $\Gamma_1$. Then the inclusion $i\colon
\Gamma_1\to \Gamma$ induces a canonical map
$$i_*\colon K^n_{\alpha_1}(\Gamma_1\upcom )\to
K^n_{\alpha}({\Gamma}\upcom ).$$
\end{prop}

\begin{pf}
Using the obvious notation,
$C^*_r(R_1)^{S^1}$ is an ideal of the $C^*$-algebra
$C^*_r(R)^{S^1}$. Indeed, it is not hard to check that
$C_c (R_1)^{S^1}\subset C_c (R)^{S^1}$ is stable
under the convolution and the adjoint. Since $R_1$ is a saturated
subgroupoid of $R$, we have
\begin{eqnarray*}
\|f\|_{C^*_r(R_1)}&\equal &\sup_{x\in R_1^{(0)}}\sup_{\xi\in C_c (R_x)}
\|f*\xi\|_{L^2(R_x)}\\
&\equal &\sup_{x\in R^{(0)}}\sup_{\xi\in C_c (R_x)}
\|f*\xi\|_{L^2(R_x)}\\
&\equal &\|f\|_{C^*_r(R)}, 
\end{eqnarray*}
and thus $C^*_r(R_1)^{S^1}$ is a sub-$C^*$-algebra of
$C^*_r(R)^{S^1}$. Moreover, for all
$f\in C_c (R_1)^{S^1}$ and $f'\in C_c (R)^{S^1}$
we have $f*f'\in C_c(R_1)^{S^1}$. Therefore $C^*_r(R_1)^{S^1}$ is an
ideal in $C^*_r(R)^{S^1}$.
\end{pf}

Recall (see for instance \cite[Section 3]{hry93})
that if $I_1$ and $I_2$ are two closed ideals in a $C^*$-algebra
$A$ such that $A\equal I_1+I_2$, then there is a six-term exact sequence
$$\xymatrix{
K_0(I_1\cap I_2)\ar[rr]^{(j_1)_*\oplus (j_2)_*}&&
K_0(I_1)\oplus K_0(I_2)
\ar[rr]^{(i_1)_*- (i_2)_*}&&
K_0(A)\ar[d]^\partial\\
K_1(A)\ar[u]^\partial&&
K_1(I_1)\oplus K_1(I_2)
\ar[ll]_{(i_1)_*- (i_2)_*}&&
K_1(I_1\cap I_2)\ar[ll]_{(j_1)_*\oplus (j_2)_*}
}$$
where $j_k\colon I_1\cap I_2\to I_k$ and $i_k\colon I_k\to A$ are the
inclusions ($k\equal 1,2$). Therefore, we get

\begin{prop}\label{prop:Mayer1}[Mayer-Vietoris sequence 1]
Let $S^1\to R\to\Gamma\toto M$ be an $S^1$-central extension of Lie
groupoids and denote by $\alpha$ its class in
$H^2(\Gamma\upcom ,{\mathcal{S}}^1)$.
Suppose that $\Gamma$ is the union of two open
saturated subgroupoids $\Gamma_1$ and $\Gamma_2$. Let
$\Gamma_{12}\equal \Gamma_1\cap \Gamma_2$, and let $\alpha_1$, $\alpha_2$
and $\alpha_{12}$ be the classes of the induced $S^1$-central extensions
and denote by
$\Gamma_{12} \stackrel{j_k}{\to}
\Gamma_k\stackrel{i_k}{\to}\Gamma$
($k\equal 1,2$) the inclusions.
Then we have an hexagonal exact sequence
$$\xymatrix{
K^0_{\alpha_{12}}(\Gamma_{12}\upcom )\ar[rr]^{(j_1)_*\oplus (j_2)_*}&&
K^0_{\alpha_{1}}(\Gamma_{1}\upcom )\oplus K^0_{\alpha_{2}}(\Gamma_{2}\upcom )
\ar[rr]^{(i_1)_*- (i_2)_*}&&
K^0_\alpha(\Gamma\upcom )\ar[d]^\partial\\
K^1_\alpha(\Gamma\upcom )\ar[u]^\partial&&
K^1_{\alpha_{1}}(\Gamma_{1}\upcom )\oplus K^1_{\alpha_{2}}(\Gamma_{2}\upcom )
\ar[ll]_{(i_1)_*- (i_2)_*}&&
K^1_{\alpha_{12}}(\Gamma_{12}\upcom )\ar[ll]_{(j_1)_*\oplus (j_2)_*}
}$$
\end{prop}

\begin{pf}
It is clear that both  $I_1\equal C^*_r(R_1)^{S^1}$ and
$I_2\equal C^*_r(R_2)^{S^1}$ are ideals of $A\equal C^*_r(R)^{S^1}$.

To check that $I_1\cap I_2\equal C^*_r(R_{12})^{S^1}$,
note that $I_1\cap I_2\equal I_1I_2$ (this is a standard result in $C^*$-algebras)
and that $f_1*f_2\in C_c(R_{12})^{S^1}$
if $f_1\in C_c(R_1)^{S^1}$ and $f_2\in C_c(R_2)^{S^1}$.

To check that $I_1+I_2\equal A$, take a partition of unity $(\varphi_1,\varphi_2)$
associated to the cover $(\Gamma_i^{(0)}/\Gamma)_{i\equal 1,2}$
of $M/\Gamma$. Let $f\in C_c(R)^{S^1}$. Then, considering $\varphi_i$
as $R$-invariant functions on $M$, we have
$f\equal (\varphi_1 f)+(\varphi_2 f)\in C_c(R_1)^{S^1}+C_c(R_2)^{S^1}$.
\end{pf}

\begin{prop}
Let $S^1\to R\to\Gamma\toto M$ be an $S^1$-central extension of Lie
groupoids and denote by $\alpha$ its class in
$H^2(\Gamma\upcom ,{\mathcal{S}}^1)$.
Assume that $\Gamma_1$ is a closed saturated subgroupoid
of $\Gamma$. Let $\alpha_1 \in H^2(\Gamma\upcom_1 ,{\mathcal{S}}^1)$
 be the class of the   corresponding $S^1$-central extension. 
Then the inclusion $i\colon
\Gamma_1\to \Gamma$ induces a canonical map
$$i^*\colon K^n_{\alpha}({\Gamma}\upcom )\to
K^n_{\alpha_1}(\Gamma_1\upcom )$$
\end{prop}

\begin{pf}
Using the obvious notation,
$C^*_r(R_1)^{S^1}$ is a quotient of the $C^*$-algebra
$C^*_r(R)^{S^1}$. Indeed, the restriction map $C_c(R)^{S^1}\to
C_c(R_1)^{S^1}$ is a surjective $*$-homomorphism of convolution algebras
and is norm-decreasing, and therefore induces a surjective $*$-homomorphism
$C^*_r(R)^{S^1}\to C^*_r(R_1)^{S^1}$.
\end{pf}

Suppose that $A_1\equal A/I_1$ and $A_2\equal A/I_2$ are two quotients of
a $C^*$-algebra such that $I_1\cap I_2\equal \{0\}$, and let
$A_{12}\equal A/(I_1+I_2)$. Denote by $p_k\colon A\to A_k$
and by $q_k\colon A_k\to A_{12}$ the quotient maps. Then there is a six-term
exact sequence
$$\xymatrix{
K_0(A_{12})\ar[d]^\partial&&
K_0(A_{1})\oplus K_0(A_{2})
\ar[ll]_{(q_1)_*\oplus (q_2)_*}&&
K_0(A)\ar[ll]_{(p_1)_*- (p_2)_*}\\
K_1(A)\ar[rr]^{(p_1)_*- (p_2)_*}&&
K_1(A_{1})\oplus K_1(A_{2})
\ar[rr]^{(q_1)_*\oplus (q_2)_*}&&
K_1(A_{12})
\ar[u]^\partial}$$
Since we cannot  locate this standard fact in the literature,
here is a sketch of the proof. For every locally compact space $X$,
we will denote by $A(X)$ the $C^*$-algebra $C_0(X,A)$.
Consider the $C^*$-algebra
$$D\equal \{(f_{-},a,f_+)\in A_1(-1,0]\oplus A\oplus A_2[0,1)
|\; f_{-}(0)\equal p_1(a)\mbox{ and }p_2(a)\equal f_+(0)\}.$$
There is an obvious exact sequence
$$0\to I_2(-1,0]\oplus I_1[0,1)\to D\to A_{12}(-1,1)\to 0.$$
Since $I_2(-1,0]\oplus I_1[0,1)$ is contractible, the six-term
exact sequence in $K$-theory yields an isomorphism $K_i(D)
\cong K_i(A_{12}(-1,1))$, hence
\begin{equation}\label{eqn:isoDA}
K_i(D)\cong K_{i+1}(A_{12}).
\end{equation}
Now, the obvious exact sequence
$$0\to A_1(-1,0)\oplus A_2(0,1)\to D\to A\to 0$$
gives a six-term exact sequence in $K$-theory, which yields
the result 
via Bott periodicity and Eq. (\ref{eqn:isoDA}).

\begin{prop}\label{prop:Mayer2}[Mayer Vietoris sequence 2]
Let $S^1\to R\to\Gamma\toto M$ be an $S^1$-central
 extension of \emph{proper} Lie
groupoids and denote by $\alpha$ its class in $H^2(\Gamma\upcom ,{\mathcal{S}}^1)$.
Suppose that $\Gamma$ is the union of two closed
saturated groupoids $\Gamma_1$ and $\Gamma_2$. Let
$\Gamma_{12}\equal \Gamma_1\cap \Gamma_2$. Let $\alpha_1$, $\alpha_2$
and $\alpha_{12}$ be the classes of their induced 
$S^1$-central extensions
and denote by
$\Gamma_{12}\stackrel{j_k}{\to}
\Gamma_k\stackrel{i_k}{\to}\Gamma$
($k\equal 1,2$) the inclusions.
Then we have an hexagonal exact sequence
$$\xymatrix{
K^0_{\alpha_{12}}(\Gamma_{12})\ar[d]^\partial&&
K^0_{\alpha_{1}}(\Gamma_{1})\oplus K^0_{\alpha_{2}}(\Gamma_{2})
\ar[ll]_{(j_1)^*\oplus (j_2)^*}&&
K^0_\alpha(\Gamma)\ar[ll]_{(i_1)^*- (i_2)^*}\\
K^1_\alpha(\Gamma)\ar[rr]^{(i_1)^*- (i_2)^*}&&
K^1_{\alpha_{1}}(\Gamma_{1})\oplus K^1_{\alpha_{2}}(\Gamma_{2})
\ar[rr]^{(j_1)^*\oplus (j_2)^*}&&
K^1_{\alpha_{12}}(\Gamma_{12})
\ar[u]^\partial}$$
\end{prop}

\begin{pf}
Let $\Gamma'_k$ ($k\equal 1,2$ or $12$) be the complementary of
$\Gamma_k$. Since $\Gamma_k$ is closed and saturated, it follows that
$\Gamma'_k$ is an open saturated subgroupoid of $\Gamma$.
With the obvious notations, write  $I_k\equal C^*_r(R'_k)^{S^1}$,
$A\equal C^*_r(R)^{S^1}$ and $A_k\equal A/I_k$, where
 $R'_k$ denotes the complementary of
$R_k$. Since $\Gamma\equal \Gamma_1\cup\Gamma_2$,
we have $I_1\cap I_2\equal \{0\}$. It is also clear that $I_1+I_2\equal I_{12}$.

To obtain the Mayer-Vietoris sequence, it  suffices to show that
$C^*_r(R_k)^{S^1}\equal A/I_k$. This is not always true for every groupoid.
I.e. the sequence
$$0\to C^*_r(R'_k)^{S^1}\to C^*_r(R)^{S^1}\to C^*_r(R_k)^{S^1}\to 0$$
is not necessarily exact.  However, the analogous sequence with
$C^*$ instead of $C^*_r$ is always exact
by the universal property of the full $C^*$-algebra of a groupoid,
and we have $C^*_r\equal C^*$ for proper groupoids (or even for amenable
groupoids).
\end{pf}

\begin{numex}
Assume that $\Gamma$ is a transformation groupoid $G\times M \toto M$,
where $G$ is a Lie group acting on $M$ properly.  Assume that
$U_1$ and $U_2$ are $G$-invariant  open submanifolds of $M$ such that
$M\equal U_1\cup U_2$. Then Proposition~\ref{prop:Mayer1} yields that 
$$\xymatrix{
K^0_{G,\alpha_{12}}(U_1\cap U_2)\ar[rr]^{(j_1)_*\oplus (j_2)_*}&&
K^0_{G,\alpha_{1}}(U_{1})\oplus K^0_{G,\alpha_{2}}(U_{2})
\ar[rr]^{(i_1)_*- (i_2)_*}&&
K^0_{G,\alpha}(M)\ar[d]^\partial\\
K^1_{G,\alpha}(M)\ar[u]^\partial&&
K^1_{G,\alpha_{1}}(U_{1})\oplus K^1_{G,\alpha_{2}}(U_{2})
\ar[ll]_{(i_1)_*- (i_2)_*}&&
K^1_{G,\alpha_{12}}(U_1\cap U_2)\ar[ll]_{(j_1)_*\oplus (j_2)_*}
}$$

Similarly, if $F_1$ and $F_2$ are  $G$-invariant closed 
submanifolds of $M$ such that $F_1\cup F_2 \equal M$, then
Proposition~\ref{prop:Mayer2}  yields that
$$\xymatrix{
K^0_{G,\alpha_{12}}(F_1\cap F_2)\ar[d]^\partial&&
K^0_{G,\alpha_{1}}(F_{1})\oplus K^0_{G,\alpha_{2}}(F_{2})
\ar[ll]_{(j_1)^*\oplus (j_2)^*}&&
K^0_{G,\alpha}(M)\ar[ll]_{(i_1)^*- (i_2)^*}\\
K^1_{G,\alpha}(M)\ar[rr]^{(i_1)^*- (i_2)^*}&&
K^1_{G,\alpha_{1}}(F_{1})\oplus K^1_{G,\alpha_{2}}(F_{2})
\ar[rr]^{(j_1)^*\oplus (j_2)^*}&&
K^1_{G,\alpha_{12}}(F_1\cap F_2)
\ar[u]^\partial}$$
\end{numex}

\subsection{The main theorem}

Let $S^1\to R\to \gm\toto M$ be  an $S^1$-central extension of groupoids, and
$P{\mathcal B}\to M$ its associated principal $PU(\HH)$-bundle over
the groupoid $\gm \toto M$ as constructed in Proposition \ref{pro:PB}.
Let 
$$\lL(\tilde{\hH}) \equal P{\mathcal B}\times_{PU(\HH )}{\mathcal{L}}(\HH )\to M$$
and
$$\kK(\tilde{\hH}) \equal P{\mathcal B}\times_{PU(\HH )}{\mathcal{K}}(\HH )\to M$$
be its associated bundles of $C^*$-algebras, where ${\mathcal{L}}(\HH )$
denotes the algebra of bounded operators on $\HH$ endowed with the
$*$-strong topology, and
${\mathcal{K}}(\HH )$ denotes the $C^*$-algebra of
compact operators on  $\HH $ endowed with the norm-topology.
The group $PU(\HH)$ acts on $\lL(\HH)$ and
$\kK(\HH)$ by conjugation.
To justify the notation, we show that these bundles are
isomorphic to those
of bounded and compact operators associated to the Hilbert bundle
$ \cup_{x\in M} L^2 (R^x)^{S^1}\otimes \HH$ (see Eq. (\ref{eq:Hx}))
as in the appendix (Propositions~\ref{prop:LeE}
and~\ref{prop:KeE}).

Indeed it is simple  to see that the fiber of
$\lL(\tilde{\hH}) \to M$ and $\kK(\tilde{\hH}) \to M$  at each $x\in M$ 
are, respectively, ${\mathcal{L}}({\mathcal H}_x)$ and
${\mathcal{K}}({\mathcal H}_x)$. The map
\begin{eqnarray*}
P\bB\times_{PU(\HH)}\lL(\HH)&\to&\coprod_{x\in M} {\mathcal{L}}
(\hH_x)\\
(u,T)&\mapsto& uTu^{-1}
\end{eqnarray*}
is clearly a  bijection. To identify the topology
of $P\bB\times_{PU(\HH)}\lL(\HH)$, we can assume that $\hH$ is a trivial
bundle (since it is locally trivial). Then
$P\bB\times_{PU(\HH)}\lL(\HH)\cong M\times PU(\HH)\times_{PU(\HH)}
\lL(\HH)\cong M\times \lL(\HH)$ is obviously the bundle of bounded
operators associated to a continuous field of Hilbert spaces
(see Proposition~\ref{prop:LeE}).
The proof for $\kK(\tilde{\hH})$ is similar (see Proposition~\ref{prop:KeE}).

The groupoid $\gm $-action on $P{\mathcal B}\to M$
induces  an action on the $C^*$-algebra bundle
$\lL(\tilde{\hH}) \to M$
(and $ \kK(\tilde{\hH}) \to M$
 respectively).
On the other hand, the associated line bundle
$L\equal R\times_{S^1}\cc \to M$ can be considered as
a Fell bundle  over the groupoid $\gm\toto M$. Therefore
the general construction of Yamagami
 (see Appendix \ref{subsec:reduced}) gives rise to a continuous
action of the groupoid $\gm \toto M$ on
 the $C^*$-bundle $\mathcal{L}( \tilde{L}^2(\Gamma;L))\to M$
(and $\mathcal{K}( \tilde{L}^2(\Gamma;L))\to M$ as well), which
extends to an action on the $C^*$-bundle 
$\mathcal{L}( \tilde{L}^2(\Gamma;L)\otimes \HH)\to M$
(and $\mathcal{K}( \tilde{L}^2(\Gamma;L)\otimes \HH)\to M$).
From Eq. (\ref{eqn:iso H}), it follows that
the Hilbert bundles $\tilde{L}^2(\Gamma; L) \otimes \HH\to M$
and $\tilde{\hH}\to M$ are canonically isomorphic. 
In fact we have the following:

\begin{lem}
The $C^*$-algebra bundles $\lL(\tilde{\hH})\to M$ 
($\kK(\tilde{\hH}) \to M$,  respectively) and
$\mathcal{L}( \tilde{L}^2(\Gamma;L)\otimes \HH)\to M$
(and $\mathcal{K}( \tilde{L}^2(\Gamma;L)\otimes \HH)\to M$
respectively) are canonically    isomorphic, 
and the isomorphism respects  the $\gm$-action.
\end{lem}

Recall that a section $x\mapsto T_x
\in {\mathcal{L}}(\tilde{\mathcal{H}}_x), \ x\in M$,
is strongly continuous if $x\mapsto T_x\xi$ is
norm-continuous for all continuous sections
 $\xi\in C(M, \tilde{\hH} )$, and that $x\mapsto T_x$
is $*$-strongly continuous if $x\mapsto T_x$ and $x\mapsto T_x^*$
are strongly continuous. 

By $C_b(M,{\mathcal{L}}(\tilde{\hH}))$ we denote the space of norm-bounded,
\emph{$*$-strongly}-continuous sections of  bounded operators
on $\hH$, and by  $C_b(M,{\mathcal{L}}(\tilde{\hH}))^{\gm}$ we denote
the subalgebra of $\gm$-invariant sections.
Similarly,
 by $C_0(M,{\mathcal{K}}(\tilde{\hH}))$
we denote the space  of \emph{norm}-continuous
sections of compact operators vanishing at infinity.

Let ${\mathcal{K}}_\Gamma (\mathcal{H})$
be the space  of norm-continuous  $\gm$-invariant
sections $\{T_x|x\in M\}$  of the $C^*$-algebra bundle 
 $\kK(\tilde{\hH})  \to M$ 
satisfying the boundary condition
$\|T_x\|\to 0 \mbox{ when }x\to\infty\mbox{ in }
M/\Gamma$. Note that $\|T_x\|$
can be considered as a function on  the 
orbit space $M/\Gamma$ due to the invariance
assumption.

Denote by ${\mathcal{F}}^0_\alpha$ the space
of $T\in C_b(M,{\mathcal{L}}(\tilde{\hH}))^{\gm}  $ such that there exists
$S\in C_b(M,{\mathcal{L}}(\tilde{\hH}))^{\gm}$ satisfying
$1-TS$, $1-ST\in {\mathcal{K}}_\Gamma({\mathcal{H}})$.
In other words,
\begin{itemize}
\item[(i)] $T_x$ and $S_x$ are Fredholm for all $x$,
and the sections
$x\mapsto T_x$ and $x\mapsto S_x$ are $*$-strongly continuous
and $\Gamma$-invariant;
\item[(ii)] $1-T_xS_x$, $1-S_xT_x$ are compact
operators for all $x$ and the sections
$x\mapsto 1-T_xS_x$, $x\mapsto 1-S_xT_x$
are norm-continuous and vanish at $\infty$ in $M/\Gamma$.
\end{itemize}

Denote by ${\mathcal{F}}^1_\alpha$ the space of self-adjoint elements
in ${\mathcal{F}}^0_\alpha$. Our main theorem is

\begin{them}\label{thm:fredholm proper}
Let $\Gamma\toto M$ be a proper Lie groupoid, $S^1\to R\to \Gamma$
an $S^1$-central extension and denote by $\alpha$ its class
in 
$H^2(\Gamma\upcom ,{\mathcal{S}}^1) (\cong H^3 (\Gamma\upcom , \zz ))$.
 Then
$$K^i_\alpha(\Gamma\upcom )\equal 
\{[T]|\; T\in {\mathcal{F}}^i_\alpha\},$$
where $[T]$ denotes the homotopy class of $T$.
\end{them}

The proof of (a generalization of) Theorem~\ref{thm:fredholm proper}
is the content of the next section.

Another way to   formulate Theorem~\ref{thm:fredholm proper} is as
follows.

 Let $\Gamma\toto M$ be  a proper Lie groupoid  and
 $\alpha\in H^3(\Gamma\upcom,\zz)$. Let
$P_\alpha\to M$ be its corresponding
 canonical $PU(\HH)$-bundle over the groupoid 
$\Gamma\toto M$ as in Proposition \ref{pro:can}.
Consider its associated bundles:
\begin{eqnarray}
&& Fred^i_\alpha  (\HH ) :\ \ \equal   P_\alpha\times_{PU(\HH)} Fred^i(\HH)\to M
\label{eqn:fredi}\\
&& \kK_\alpha (\HH ): \ \ \equal  P_\alpha\times_{PU(\HH)} \kK(\HH) \to M, \label{eqn:compi}
\end{eqnarray}
where $Fred^i(\HH)$ is endowed with the $*$-strong topology while
$\kK(\HH)$ is endowed with the norm-topology.
By $\fF^i_\alpha$, we denote the the space of
norm-bounded, $\Gamma$-invariant,
continuous sections $x\mapsto T_x$ of the bundle
$ Fred^i_\alpha  (\HH )\to M$ such that there exists a norm-bounded,
$\Gamma$-invariant, continuous section $x\mapsto S_x$ of
$Fred^i_\alpha (\HH )\to M$ with the property  that
$1-T_xS_x$ and $1-S_xT_x$ are continuous sections of
$\kK_\alpha (\HH)$ vanishing
at infinity\footnote{This definition of $\fF^i_\alpha$ obviously agrees
with the previous one.}.
Then we have the following

\begin{them}\label{thm:fredholm proper1}
Let $\Gamma\toto M$ be a proper Lie groupoid, 
and  $\alpha \in H^3 (\Gamma\upcom , \zz )$.
 Then
$$K^i_\alpha(\Gamma\upcom )\equal
\{[T]|\; T\in {\mathcal{F}}^i_\alpha\},$$
where $[T]$ denotes the homotopy class of $T$.
\end{them}

\begin{rmk}
Note that there may exist different
 $PU(\HH)$-principal bundles over $\Gamma\toto M$
other than $P_\alpha$, which   also  map to
 $\alpha\in H^2(\Gamma\upcom,\sS^1)$ under the map $\Phi'$.
However,
only the construction using this particular
principal bundle $P_\alpha$ gives the right answer for the
Fredholm picture of twisted $K$-theory groups.
\end{rmk}

\begin{numex}
\begin{enumerate}
\item When $\gm$ is a  compact 
manifold $M$, the principal $PU(\HH)$-bundle
$P_\alpha\to M$  over $M$ is represented by a 1-cocycle
$g_{kl}:  U_{kl}\to PU(\HH)$. A class $[T]$ in  
$K^{i}_\alpha(M )$ corresponds to a section of the bundle
$P_\alpha\times_{PU(\HH)} Fred^i(\HH)$, thus to 
a family of  $*$-strongly continuous
maps $T_k: U_k \to   Fred^{i} (\HH )$ satisfying
$T_l \equal g_{kl}^{-1}T_k g_{kl}$ on $U_{kl}$ \cite{Atiyah, F}.

\item When $\Gamma$ is a transformation groupoid $G\times M\toto M$,
for any $\alpha \in H^3_G (M, \zz )$, the associated
$PU(\HH)$-bundle $P_\alpha\to M$ over $\Gamma$
is  a  $G$-equivariant $PU(\HH)$-bundle over $M$. Therefore,
its associated bundles $ Fred^i_\alpha  (\HH ) \to M$
and $\kK_\alpha (\HH )\to M$ are $G$-equivariant.
Thus $K_{G, \alpha}^{i}(M)$ can be   represented as the group
of homotopy classes of $G$-invariant $*$-strongly  continuous sections of
the Fredholm bundle $ Fred^i_\alpha  (\HH ) \to M$.

In terms of  local charts,
the principal $PU(\HH)$-bundle $P_\alpha$ is represented by a $1$-cocycle
$\varphi_{kl}\colon (G\times M)_{U_k}^{U_l}\to PU(\HH)$ (see
Eqs~(\ref{eqn:projective rep}) and
(\ref{eqn:BR0})).
Then a class $[T]$ in $K_{G, \alpha}^{i}(M)$ corresponds to
a family of $*$-strongly continuous maps $T_k\colon U_k\to Fred^{i}(\HH)$
satisfying $T_l(y)\equal \varphi_{lk}(g, x)T_k(x)\varphi_{kl}
(g^{-1}, y)$ if $x\in U_k$, $y\in U_l$ and $y\equal gx$.
\end{enumerate}   
\end{numex}


\section{$C^*$-algebras of Fell bundles over proper groupoids}

The goal of this section is to prove a general result
in $C^*$-algebras,  which includes
Theorem \ref{thm:fredholm proper} as a special case.
More precisely, we prove that if $\Gamma\toto M$ is a proper
groupoid and ${{E}}\equal \coprod_{g\in \Gamma} E_g$
is a u.s.c. Fell bundle over $\Gamma$ (see Appendix),
then its associated $C^*$-algebra $C^*_r(\Gamma;E)$  is
isomorphic to the space of $\Gamma$-equivariant continuous sections
of the $C^*$-algebra bundle $\coprod_{x\in M} {\mathcal{K}}
(L^2(\Gamma_x;E))\to M$ that vanish at infinity in $M/\Gamma$
(Proposition~\ref{prop:structure of algebra of action}),
and to deduce that the $K$-theory groups of $C^*_r(\Gamma;E)$
are isomorphic to homotopy classes of certain generalized
Fredholm operators (Theorem~\ref{thm:K theory Fredholm}).

Let us recall
\begin{defn}\label{defi:cutoff}
Let $\Gamma\toto M$ be a proper groupoid with Haar system $\{\lambda^x |
x\in M\}$.
A continuous function $c\colon M \to \rr_+$ is
called a cutoff function if
\begin{itemize}
\item[(i)] for all $x\in M$, $\int_{\gamma\in \Gamma^x}
c(s(\gamma))\,\lambda^x(d\gamma)\equal 1$; and
\item[(ii)] for all $K\subset M$ compact, the support of
$(c\smalcirc s)_{|\Gamma^K}$ is compact.
\end{itemize}
\end{defn}

The condition  (ii) means that if $F\subset M$ is the saturate of a
compact set, then $F\cap{\mathrm{supp}}(c)$ is compact.
It is known that  a cutoff  function exists if and only if $\Gamma$
is proper \cite[Proposition~6.7]{tu99}.

Cutoff functions allow us to make ``averages''. Namely, let
\begin{equation}\label{eqn:TGamma}
T^\Gamma_x\equal \int_{\gamma\in\Gamma^x}\alpha_\gamma(T_{s(\gamma)})
c(s(\gamma))\,\lambda^x(d\gamma) \in \lL(L^2(\Gamma_x;E)).
\end{equation}

Then $T\mapsto T^\Gamma$ is a linear projection of norm one from
${\mathcal{L}}(L^2(\Gamma;E))$ onto ${\mathcal{L}}(L^2(\Gamma;E))^\Gamma$.
(If $\Gamma\equal G$ is a compact group then $T^G$ is the average
$\int_G \alpha_g(T)\,dg$.)

More generally,
let $\Gamma \toto M$ be  a proper groupoid with Haar system
acting continuously on a u.s.c. field  of  $C^*$-algebras $\aA \to M$.
By $A\equal C_0(M,\aA)$ we denote its corresponding $C^*$-algebra
of continuous sections vanishing at infinity.
As in Proposition~\ref{prop:LeE}, there exists a (not u.s.c.)
field of $C^*$-algebras $M (\aA )\to M$
with the fiber at $x\in M$ being $M (\aA_x )$, such that
\begin{itemize}
\item[(a)] a section $x\mapsto T_x\in M(\aA_x)$ is a continuous section
of $M(\aA)\to M$ if and only if for every continuous section
$x\mapsto a_x$ of $\aA\to M$, $x\mapsto T_xa_x$ and $x\mapsto T_x^*a_x$
are continuous sections of $\aA\to M$.
\item[(b)] The algebra $C_b(M,M(\aA))$ of continuous, norm-bounded sections
is isomorphic to the multiplier algebra $M(A)$.
\end{itemize}

It is clear that the $\gm$-action on $\aA \to M$ induces
a natural $\gm$-action on $M (\aA )\to M$.  
By $ M(A)^\Gamma$, we denote the $C^*$-subalgebra of $M(A)$
consisting of $\gm$-invariant sections.
For any $T\in A$, let $T^\Gamma\in M(A)$ be the element such that
$$(T^\Gamma)_x\equal \int_{\gamma\in \Gamma^x} \alpha_\gamma(T_{s(\gamma)})
c(s(\gamma))\,\lambda^x(d\gamma)\quad
\in M({\mathcal{A}}_x). $$
Of course, the averaging map
$T\mapsto T^\Gamma$ depends on the choice of the cutoff function.
\par\medskip

Let us introduce some notations. If ${\mathcal{E}}$
is a $C_0(M)$-Hilbert module, let
\begin{equation}\label{eqn:C(E)}
{\mathcal{C}}({\mathcal{E}})\equal \{ T\in {\mathcal{L}}({\mathcal{E}})|\quad
\varphi T\in
 {\mathcal{K}}({\mathcal{E}})\quad\forall \varphi\in C_0(M)\}.
\end{equation}

If,  moreover,  ${\mathcal{E}}$
is a $(\Gamma,E)$-equivariant Hilbert module,
let
\begin{equation}\label{eqn:K(E)}
{\mathcal{K}}_\Gamma({\mathcal{E}})\equal \{
T\in {\mathcal{C}}({\mathcal{E}})^\Gamma|\quad
\|T_x\|\to 0\mbox{ when }x\to \infty\mbox{ in }
M/\Gamma\}.
\end{equation}

More generally, if $\Gamma$ acts on a $C_0(M)$-algebra $A$, let
\begin{eqnarray}\nonumber
\lefteqn{A_\Gamma\equal \{T\in M(A)^\Gamma|\;
\forall \varphi\in C_0(M), \varphi T\in A}\\
\label{eqn:AG}
&&\qquad\mbox{and }
\|T_x\|\to 0 \mbox{ when } x\to\infty \mbox{ in } M/\Gamma\}.
\end{eqnarray}

For example, if $A\equal C_0(\rr)$ and $\Gamma\equal \zz$ acts on $A$ by
translations, then $A_\Gamma$ is the space of $1$-periodic
continuous functions on $\rr$, while the algebra $A^\Gamma$, consisting
of $\Gamma$-invariant elements in $A$,  is $\{0\}$.

\begin{lem}\label{lem:A_Gamma}
With the above assumptions, $A_\Gamma$ is a  $C^*$-subalgebra
of $M(A)$ and is equal to $\{T^\Gamma|\;T\in A\}$.
\end{lem}

\begin{pf}
The first assertion is easy.
\par\medskip

To show $\{T^\Gamma|\;T\in A\}\subset A_\Gamma$, let $T\in A$.
By a density argument, we may assume that $T$ is supported
on a compact subset of $M$.  That is,  $T_x\equal 0$
for $x$ outside a compact set. Let $\varphi\in C_c(M)$.
Then $(\varphi \cdot T^\Gamma)_x\equal \int_{\gamma\in \Gamma^x}\varphi(x)
\alpha_\gamma(T_{s(\gamma)})c(s(\gamma))\,d\gamma$ is the integral
on $\Gamma^x$
of a compactly supported element of $t^*A$. Thus
$\varphi \cdot T^\Gamma$
belongs to $A$. Moreover, it is clear that $T^\Gamma$
is zero outside a compact subset of $M/\Gamma$.
\par\medskip

To show that $\{T^\Gamma|\; T\in A\}
\supset A_\Gamma$,
let $T\in A_\Gamma$. Assume  first that
$T_x\equal 0$ outside a compact set $K$ of $M/\Gamma$. Let
$L\equal ({\mathrm{supp}}(c))\cap \pi^{-1}(K)$ where $\pi\colon M
\to M/\Gamma$ is the projection. Then, for all
$\varphi\in C_0(M)$ such that $\varphi\equal 1$ on $L$, one has
\begin{eqnarray*}
T_x&\equal &T_x\int_{g\in \Gamma^x} c(s(g))\varphi(s(g))\,\lambda^x(dg)\\
&\equal &\int_{g\in \Gamma^x} T_xc(s(g))\varphi(s(g))\,\lambda^x(dg)\\
&\equal &\int_{g\in \Gamma^x} \alpha_g(T_{s(g)})c(s(g))\varphi(s(g))\,\lambda^x(dg).
\end{eqnarray*}
Thus $T\equal (\varphi T)^\Gamma$.

In the general case, one has  $T\equal \sum_nT_n$ where $T_n$ is zero outside
a compact subset of $M/\Gamma$ and $\|T_n\|\le 2^{-n}$
for $n$ large enough. From the previous paragraph, we have
$T_n\equal (\varphi_n T_n)^\Gamma$,  and therefore
$T\equal (\sum_n\varphi_nT_n)^\Gamma$.
\end{pf}

\begin{prop}\label{prop:structure of algebra of action}
If $\Gamma\toto M$ is a proper groupoid with Haar system and $E$
is a u.s.c. Fell bundle over $\Gamma$, then
$$C^*_r(\Gamma;E)\equal \{T^\Gamma|\; T\in{\mathcal{K}}(L^2(\Gamma;E))\}
\equal {\mathcal{K}}_\Gamma(L^2(\Gamma;E)).$$
\end{prop}
\begin{pf}
Let us explain the idea of the proof in the case that $\Gamma$ is
a compact group $G$. In this case $C^*_r(G)$ is the closure of the
space of convolution operators on $L^2(G)$. Since these operators
have a $G$-invariant kernel $K\in C(G\times G)$, they are
compact $G$-invariant. Conversely, any compact invariant operator
$T$ is the limit of operators with $G$-invariant kernel, and such kernels
are of the the form $a(gh^{-1})$ where $a\in C(G)$.
It follows that $T$ is the limit of convolution operators. 
\par\medskip

Now we consider the general case of a proper groupoid.
Let us first show that
$C^*_r(\Gamma;E)\supset
\{T^\Gamma|\; T\in{\mathcal{K}}(L^2(\Gamma;E))\}$.
Let $T\in {\mathcal{K}}(L^2(\Gamma;E))$. 
We need to show that $T^\Gamma$
lies in the image of $C^*_r(\Gamma;E)$.
  We may assume that
$T\equal T_{b,b'}:\equal (\xi\mapsto b\langle b',\xi\rangle)$, where
$b$, $b'\in C_c(\Gamma;E)$, i.e.
$$(T\xi)(g)\equal b(g)\int_{h\in \Gamma_{s(g)}} b'(h)^*\xi(h)\,\lambda_{s(g)}(dh).$$
Then,
$$(T^\Gamma\xi)(g)\equal \int_{h\in \Gamma_{s(g)}}\int_{\gamma\in
\Gamma^{s(h)}} b(g\gamma)b'(h\gamma)^*\xi(h)c(s(\gamma))\,
\lambda^{s(h)}(d\gamma)\lambda_{s(g)}(dh).$$
Set $f(g)\equal \int_{\gamma\in \Gamma^{s(g)}} b(g\gamma)b'(\gamma)^*c(s(\gamma))
\,\lambda^{s(g)}(d\gamma)$. Let us check that $f\in C_c(\Gamma;E)$.
By Proposition~\ref{prop:approximation}, $(g,h)\mapsto
b(g)b'(h^{-1})^*$ can be approximated by sums of the form
$\sum_i f_i(g,h)\zeta_i(gh)$,  where $f_i\in C_c(\Gamma)$ and $\zeta_i
\in C_c(\Gamma;E)$. Therefore $(g,\gamma)\mapsto b(g\gamma)d(\gamma)^*$
is approximated by $\sum_i f_i(g\gamma,\gamma^{-1})\zeta_i(g)$.
Then, approximating $f_i(g\gamma,\gamma^{-1})$ by $\sum_j f_{1,i,j}(g)
f_{2,i,j}(\gamma)$, we see that $f(g)$ is approximated by
$\sum_if'_i(g)\zeta_i(g)$ where $f'_i\in C_c(\Gamma)$.  Hence
$f$ is a continuous section. Moreover, since $f$ is obviously
compactly supported, we have $f\in C_c(\Gamma;E)$.
Now,

$$f(gh^{-1})\equal \int_{\gamma\in \Gamma^{s(g)}}
b(g\gamma)b'(h\gamma)^*c(s(\gamma))\,\lambda^{s(g)}(d\gamma),$$
and $\pi_l(f)\xi(g)\equal \int_{h\in \Gamma_{s(g)}}f(gh^{-1})\xi(h)\,
\lambda_{s(g)}(dh)$, where $\pi_l\colon C_c(\Gamma;E)\to \lL(L^2(\Gamma;E))$
denotes the left regular representation.
Therefore, $T_{b,b'}^\Gamma\equal \pi_l(f)$.
\par\medskip

Next we  show that $C^*_r(\Gamma;E)\subset
\overline{\{T^\Gamma|\; T\in{\mathcal{K}}(L^2(\Gamma;E))\}}$.
Assume that  $a\in C_c(\Gamma;E)$. Let $K$ be a compact subset of $M$
such that $\Gamma^K_K$ contains the support of $a$, and 
$K_1$ a compact subset  such that the interior of ${K}_1$ contains $K$.
By the definition of $c$, the set
$L\equal ({\mathrm{supp}}(c))\cap\pi^{-1}(\pi(K_1) )$, where $\pi\colon
M\to M/\Gamma$ is the projection,
is compact.
By Proposition~\ref{prop:approximation}, one may approximate
$a(gh^{-1})$ on the compact set
$\Gamma_L^{K_1}\times \Gamma_L^{K_1}$ uniformly by
elements of the form
$$\sum_i b_i(g)b'_i(h)^*, $$
where $b_i$, $b'_i\in E$.
Therefore $(g,h)\mapsto a(gh^{-1})$ is approximated
uniformly on $\{(g,h)\in\Gamma^{K_1}\times\Gamma^{K_1}|\;
s(g)\equal s(h)\}$ by elements of the form
$\sum_i\int_{\gamma\in\Gamma^{s(g)}} b_i(g\gamma)b'_i(h\gamma)^*c(s(\gamma))
\,d\gamma$.
Replacing $b_i(g)$ by $\varphi(t(g))b_i(g)$ and
$b'_i(g)$ by $\varphi(t(g))b'_i(g)$,  where $\varphi\in C_c(M)$,
$0\le\varphi\le 1$, $\varphi\equal 1$ on $K$ and $\varphi\equal 0$ on $M
-K_1$, we define an operator

$$T\equal \sum T_{b_i,b'_i}.$$
Then $\pi_l(a)\in {\mathcal{L}}(L^2(\Gamma;E))$
is approximated by the operator $T^\Gamma$.

The other inclusions are proved using Lemma~\ref{lem:A_Gamma} with
$A\equal {\mathcal{K}}(L^2(\Gamma;E))$.
\end{pf}

\par\bigskip
To continue, let us introduce the following convention.
 For any $C^*$-algebras $A$ and $B$ such that
$B\subset M(A)$, we say that $B$ contains an approximate unit
for $A$ if there exist $u_i\in B$ such that $0\le u_i\le 1$
and $u_ia\to a$ for all $a\in A$. This terminology is
slightly abusive since $u_i$ may not belong to $A$.

\begin{lem}\label{lem:AGammamultiplier}
Let $\Gamma\toto M$ be a proper groupoid with Haar system
acting on a  u.s.c. field  of  $C^*$-algebras $\aA \to M$, and
 $A\equal C_0(M,\aA)$. Then 
\begin{itemize}
\item[(a)] $A_\Gamma$
contains an approximate unit for $A$;
\item[(b)] $\overline{A_\Gamma A}\equal A$;
\item[(c)] $M(A_\Gamma)\equal M(A)^\Gamma$;
\end{itemize}
\end{lem}

\begin{pf}
(a) Let $(u_i)_{i\in I}$ be an approximate unit in $A$
(it is standard that this always exists, see~\cite{dix77}). Let
$\tilde{u}_i\equal (u_i)^{\Gamma}$.
Then $\tilde{u}_i\in A_\Gamma$ by Lemma~\ref{lem:A_Gamma}.
It suffices to show that for all $a\in A_U$, where $U$ is a 
 relatively compact open subset of $M$, we have $\tilde{u}_ia\to a$.\par\medskip

Let $b(g)\equal a_{t(g)}c(s(g))$. Then $b\equal (t^*a)(c\smalcirc s)\in t^*A$.
Set  $v_i(g)\equal \alpha_g((u_i)_{s(g)})$. Since $\alpha\colon s^*A
\stackrel{\sim}{\to} t^*A$ is an isomorphism and 
$v_i\equal \alpha (s^*u_i)$, it follows that $v_i$ is an approximate unit
for $t^*A$. Thus $v_i b\to b$, i.e.,
$$\sup_{g\in \Gamma}\|a_{t(g)}c(s(g))-\alpha_g(u_i)_{s(g)}a_{t(g)}c(s(g))\|
\to 0.$$
By integration on $\Gamma^x$ ($x\in M$), it follows easily that
$\|a-\tilde{u}_ia\|\to 0$.
\par\medskip

(b) clearly follows  from (a).
\par\medskip

(c) The map $ M(A)^\Gamma\to M(A_\Gamma), \ \ a\mapsto  \mu(a)$,
where $\mu(a)b\equal ab$, $\forall a\in M(A)^\Gamma, \  b\in A_\Gamma$, 
is well-defined and $*$-linear.
To obtain  its inverse, by identifying  $M(A)$ with
${\mathcal{L}}(A)$ (the space of $A$-linear adjointable
operators on the $A$-Hilbert module $A$),  the map
$\nu\colon T\in M(A_\Gamma)\equal {\mathcal{L}}(A_\Gamma)\mapsto T\otimes 1
\in {\mathcal{L}}(A_\Gamma\otimes_AA)\equal {\mathcal{L}}(\overline{A_\Gamma A})
\equal {\mathcal{L}}(A)\equal M(A)$ takes its value in $M(A)^\Gamma$.
 It is clear that $\nu\colon M(A_\Gamma)\to M(A)^\Gamma$ and $\mu$
are inverse of  each other.
\end{pf}


\begin{cor}\label{coro:multiplier}
If $\Gamma$ is a proper groupoid with Haar system 
and $E$ is a u.s.c. Fell bundle over $\Gamma$, then
$$M({\mathcal{K}}_\Gamma(L^2(\Gamma;E )\otimes\HH))
\equal {\mathcal{L}}(L^2(\Gamma;E )\otimes \HH)^\Gamma.$$ 
\end{cor}

\begin{pf}
For a $C^*$-algebra $A$ and an  $A$-Hilbert module ${\mathcal{E}}$,
we have ${\mathcal{L}}({\mathcal{E}})\equal M({\mathcal{K}}({\mathcal{E}}))$
\cite[Thm 13.4.1]{black98}. Hence the result follows from Lemma \ref{lem:AGammamultiplier} (c).
\end{pf}

Let ${\mathcal{F}}^0(\Gamma,E)$
be the set consisting of all  $T\in {\mathcal{L}}(L^2(\Gamma;E)\otimes \HH)^\Gamma$
which are invertible modulo
${\mathcal{K}}_\Gamma(L^2(\Gamma;E )\otimes\HH)$,
and ${\mathcal{F}}^1(\Gamma,E)$ the subset of ${\mathcal{F}}^0(\Gamma,E)$
consisting of  self-adjoint
elements.  We denote by $[T]$ the
homotopy class of $T$.

\begin{them}\label{thm:K theory Fredholm}
Let $\Gamma$ be a proper groupoid with a Haar system.
Suppose that $E\equal (E_g)_{g\in\Gamma}$
is a u.s.c. Fell bundle over $\Gamma$. Then
$$K_i(C^*_r(\Gamma;E))\equal \{[T]|\; T\in {\mathcal{F}}^i(\Gamma,E)\}.$$
\end{them}

\begin{pf}
 Recall that if $B$ is a $C^*$-algebra then $K_0(B)$
is the set of homotopy classes of elements
$T\in M(B\otimes{\mathcal{K}}(\HH ))$
which are invertible modulo $B\otimes {\mathcal{K} (\HH )}$, and $K_1(B)$
is the set of homotopy classes of elements
$T\in M(B\otimes{\mathcal{K}}(\HH ))$
which are self-adjoint and invertible
modulo $B\otimes {\mathcal{K}}(\HH )$ (\cite[Cor 12.2.3]{black98},
\cite[Thm. 17.3.11]{wegge93}).

The theorem is thus a consequence of
Proposition~\ref{prop:structure of algebra of action}
and  Corollary~\ref{coro:multiplier},  by taking
$B\equal {\mathcal{K}}_\Gamma(L^2(\Gamma,{{E}}))$.
\end{pf}

\section{Twisted vector bundles}

In many situations, it is desirable to describe the $K$-theory 
groups in terms of geometrical objects such as vector bundles.
For  the twisted $K$-theory group $K^0_\alpha(\Gamma\upcom )$,
a natural candidate will be twisted vector bundles. However,
these vector bundles do not always exist. In fact, a necessary
condition is that the twisted class 
$\alpha \in H^2 (\gm\upcom, {\mathcal{S}}^1 ) $ must be
a torsion class. The main purpose of this section is
to explore the conditions under which $K^0_\alpha(\Gamma\upcom )$
can be expressed by twisted vector bundles. More precisely,
we prove that given an $S^1$-central extension
$S^1\to R\to\Gamma\toto M$ of a proper
Lie groupoid $\Gamma$ such that $M/\Gamma$ is compact,
the $K$-theory group $K^0_\alpha(\Gamma\upcom )$
twisted by the class $\alpha$ of the above central extension
is the Grothendieck group of twisted vector bundles
$K^{vb}_\alpha(\Gamma\upcom)$, provided
some conditions are fulfilled (see Theorem~\ref{thm:vector bundles}).
\par\medskip

The proof is divided into five steps outlined  as follows.
Let $L\equal R\times_{S^1}\cc$ be the associated line bundle over
$\gm$.

\begin{itemize}
\item[Step 1:] From the previous
 section,   it is known that $K^0_\alpha(\Gamma\upcom)$ is
  isomorphic to
$K_0({\mathcal{K}}_\Gamma(L^2(\Gamma; L)\otimes\hh ))$. Therefore,
if ${\mathcal{K}}_\Gamma(L^2(\Gamma;L)\otimes\hh)$ 
has an approximate unit consisting of projections, then
$K^0_\alpha(\Gamma\upcom)$ is the Grothendieck group of projections in
${\mathcal{K}}_\Gamma(L^2(\Gamma;L)\otimes\hh)$
\cite[Prop. 5.5.5]{black98};
\item[Step 2:] ${\mathcal{K}}_\Gamma(L^2(\Gamma;L)\otimes\hh)$ 
has an approximate unit consisting of projections if and only if
the $(\Gamma,L)$-equivariant Hilbert module $L^2(\Gamma;L)\otimes\hh$ 
satisfies a certain property  that we denote by \emph{AFGP};
\item[Step 3:] If $L^2(\Gamma)\otimes\hh$ is AFGP and if there exists
a twisted vector bundle, then $L^2(\Gamma;L)\otimes\hh$ 
is AFGP;
\item[Step 4:] Projections in ${\mathcal{K}}_\Gamma(L^2(\Gamma;L)\otimes\hh)$
correspond to $(\Gamma,L)$-equivariant Hilbert modules ${\mathcal{E}}$
such that ${\mathrm{Id}}_{\mathcal{E}}\in{\mathcal{K}}_\Gamma({\mathcal{E}})$
(see notation~(\ref{eqn:K(E)}));
\item[Step 5:] $(\Gamma,L)$-equivariant Hilbert modules ${\mathcal{E}}$
such that ${\mathrm{Id}}_{\mathcal{E}}\in{\mathcal{K}}_\Gamma({\mathcal{E}})$
correspond to twisted
vector bundles, which can be considered as  a generalization of Swan's theorem.
\end{itemize}

%
%
%

\subsection{Definition of  twisted vector bundles}

In this subsection, we give the definition of a twisted vector bundle
and show that if such a vector bundle exists, then the $S^1$-central
extension must be a  torsion.

Let us first recall the definition of a $\gm$-vector bundle.

\begin {defn}
Let $\gm\toto M$ be a groupoid.
A $\gm$-vector bundle is a vector bundle
$J: E\to M$ such that $E$ is a $\gm$-space in the sense of
Definition \ref{def:principal}, and the map (\ref{eq:lr})
is  a linear map.
\end{defn}

Note that in this case for any $r\in \gm$, the map
\begin{equation}
l_{r}: J^{-1}(u)\to J^{-1}(v), \ \ \ x\to r\cdot x,
\end{equation}
where $u\equal s(r)$ and $v\equal t(r)$, must be  a linear isomorphism.

For example, given a $G$-bundle $P$
over $\gm\toto M$ and  a representation
$G\to \mbox{End} V$, the associated vector bundle
$E: \equal (P\times V)/G\to M$
 naturally becomes a  $\gm$-vector bundle.

\begin{defn}\label{def:twisted vb}
Let $S^1\to R\to\Gamma\toto M$ be an $S^1$-central extension of Lie
groupoids.  By a $(\gm,R)$-twisted vector bundle, we mean an $R$-vector bundle
satisfying the compatibility condition:
$$ (\lambda \cdot r)\cdot x\equal \lambda (r \cdot x ), \ \ \lambda \in S^1,
r\in R \ \mbox{ and } \  x\in E  \mbox{ such that } s(r)\equal J(x).$$
Here $S^1$ is considered as the unit circle in $\cc$.
\end{defn}

The following gives an equivalent definition of twisted vector bundles.

\begin{lem}
Let $S^1\to R\stackrel{\pi}{\to}\Gamma\toto M$ be  an $S^1$-central
 extension of Lie
groupoids.  An $R$-vector bundle $E\to M$
is a $(\gm, R)$-twisted vector bundle if and only
$\ker \pi \cong M\times S^{1}$ acts on $E$ by
scalar multiplication, where $S^1$ is identified with
the unit circle of $\cc$.
\end{lem}

When $M$ is a  point,  the definition above reduces to the
usual projective representations of a  group.

\begin{numex}
\begin{enumerate}
\item  Consider the  $S^1$-central extension as in Example \ref{ex:gerbe} (1).
A twisted vector bundle  $E\to\coprod_i U_i$  of rank $n$
 corresponds to vector bundles $E_i \cong U_{i} \times \cc^n$,
where the transition functions $a_{ij}: U_{ij}\to GL(n, \cc)$
satisfy the twisted cocycle condition
$$a_{ij}a_{jk}a_{ki}\equal c_{ijk}.$$
Note that when the central extension is trivial, i.e., $c_{ijk}\equal 1$,
the transition functions $(a_{ij})$ define  an ordinary vector bundle over $M$.
In other words, a vector bundle over the groupoid $\coprod_{ij}
U_{ij}\toto \coprod_{i}U_i$ corresponds exactly 
to a vector bundle over $M$ in the usual sense.

\item Consider the    $S^1$-central extension as in Example 
\ref{ex:gerbe} (2).
Let $E\to \coprod_i U_i$ be a twisted vector bundle of rank $n$.
Then $E|_{U_i}\cong U_i\times \cc^n$.
For all $x\in U_i$ and $\xi\in \cc^n$,
denote by $[(i,x,\xi)]$ the corresponding element of $E|_{U_i}$.
Write 
$$(\alpha,g,x,\lambda)\cdot[(j,x,\xi)]\equal [(i,\lambda
a_{ij;\alpha}(g,x)\xi)], $$
 where $a_{ij;\alpha}\colon G\times U_{ij}\to GL_n(\cc)$.
Then we have the cocycle relation
$$a_{ij;\alpha}(g,x)a_{jk;\beta}(h,y)
\equal c_{ijk;\alpha\beta,\gamma}(g,x,h,y)a_{ik;\gamma}(gh,y).$$

\item Consider the case that $R$ is topologically trivial and 
the $S^1$-central extension is given by a  groupoid $S^1$-valued 2-cocycle
$c(x, y)$ as in Eq. (\ref{eq:extension}).  Let $E\to M$ be  a
(non-equivariant) trivial vector bundle over $M$, i.e., $E\cong M\times \cc^n$.
Then $E\to M$ defines a twisted vector bundle of  $R$
if and only if there is a smooth map $\phi: \gm \to GL(n, \cc)$
satisfying the condition:
$$\phi (x) \phi (y)\equal c(x, y) \phi (xy ), \ \ \ \forall (x, y )\in
\gm^{(2)}.$$
\end{enumerate} 
\end{numex}


In the proposition below, we show that twisted vector bundles exist
only when the $S^1$-central extension defines
 a torsion class in $H^2(\Gamma\upcom ,{\mathcal{S}}^1)$.

\begin{prop}\label{prop:exists twisted}
Let $S^1\to R\to \Gamma\toto M$ be an $S^1$-central extension of Lie groupoids.
Consider the following properties:
\begin{itemize}
\item[(i)] there exists a rank $n$ twisted vector bundle;
\item[(ii)] there exists an $S^1$-equivariant generalized
homomorphism $R\to GL_n(\cc)$, where  $GL_n(\cc)$ is naturally
considered as an  $S^1$-central extension of
$PGL_n(\cc)$: $S^1\to GL_n(\cc) \to PGL_n(\cc)$;
\item[(ii)'] there exists a generalized homomorphism 
$\gm \to PGL_n(\cc)$ such that $R$ is the pull-back of $GL_n(\cc) \to PGL_n(\cc)$;
\item[(iii)] there exists an open cover $(U_i)$ and
$\psi\colon R[U_i]\to \cc^*$ such that
$\psi(\lambda r)\equal \lambda^n \psi(r)$ for all $\lambda\in S^1$ and $r\in R$;
\item[(iv)] $R^n$ is a trivial extension;
\item[(v)] there exists an open cover $(U_i)$ and  $\zz_n$-central extensions
$\zz_n\to R'[U_i]\to \Gamma[U_i]\toto \coprod U_i$
such that $R[U_i]\equal R'\times_{\zz_n}S^1$, where $\zz_n$ is identified
with  the group of $n$-th roots of unity in $S^1\subset\cc$.
\end{itemize}
Then $(i)\iff (ii)\implies (iii)\iff (iv)\iff (v)$.
\end{prop}
\begin{pf}
(i)$\implies$(ii): let $E$ be a rank $n$ twisted vector bundle.
Since $E$ is locally trivial, replacing
$R$ by $R[U_i]$ one may assume that $E\cong M\times \cc^n$ as
a (non-equivariant) vector bundle. Hence the action of $R$ on 
$M\times \cc^n$ defines an $S^1$-equivariant homomorphism.

(ii)$\implies$(i): let $(U_i)$ be an open cover of $M$ such that
there exists an $S^1$-equivariant strict homomorphism
$R[U_i]\to GL_n(\cc)$ (see Proposition~\ref{prop:homomorphism cover}).
Let $Z\equal \coprod R_{U_i} \equal  \{(r,i)|\; s(r)\in U_i\}$.
Then $Z$ is naturally endowed with a right $R[U_i]$-action.
Let $E\equal Z\times_{R[U_i]}\cc^n$. Then $E$ is a twisted
vector bundle of rank $n$ (where the map $E\to M$
is $(r,i,\xi)\mapsto t(r)$).

(ii)$\implies$(iii): compose with the determinant $GL_n(\cc)\to \cc^*$.

(iii)$\implies$(iv): replacing $\psi$ by $\psi/|\psi|$, we may
assume  that the image of $\psi$ lies in $S^1$. Define
$\varphi(\lambda[r,\ldots,r])\equal \lambda\psi(r)$.
Then $\varphi$ is a well-defined $S^1$-equivariant
homomorphism from $R^n[U_i]$ to $S^1$.
Hence $R^n$ is a trivial extension
(see Proposition~\ref{prop:trivial extensions}).

(iv)$\implies$(iii): If $\varphi\colon R^n[U_i]\to S^1$ is $S^1$-equivariant,
then  $\psi(r): \equal \varphi[(r,\ldots,r)]: R[U_i]\to S^1 $ is the function
satisfying the desired property.

(iii)$\implies$(v): take $R'\equal \psi^{-1}(1)\subset R[U_i]$. Then
$(r,\lambda)\in R'\times_{\zz_n}S^1\mapsto \lambda r\in R[U_i]$
is an isomorphism.

(v)$\implies$(iii): the map $[(r,\lambda )]\in R'\times_{\zz_n}S^1
\mapsto \lambda^n \in S^1$ is well-defined and satisfies (iii).
\end{pf}

\begin{rmk}
It is worth noting that (v) means that the class in
 $H^2 (\gm\upcom , \sS^1 )$
defined by  the $S^1$-central extension $R\to \gm$ lies
 in the image of the homomorphism
$H^2 (\gm\upcom , \zz_n )\to H^2 (\gm\upcom , \sS^1 )$.
\end{rmk} 

By $\tilde{K}^0(\Gamma, R)$, we denote  the Grothendieck group 
of twisted vector bundles. As an immediate consequence of
$(i)\iff (ii)$, we have

\begin{cor}
Assume that $S^1\to R_i\to\Gamma_i\toto M_i$, $i\equal 1, 2$,
are 
Morita equivalent $S^1$-central extensions of groupoids.
Then 
$$\tilde{K}^0(\Gamma_1, R_1)\cong \tilde{K}^0(\Gamma_2, R_2).$$
\end{cor}

This allows us to introduce the following

\begin{defn}
Let $\Gamma\toto M$ be a Lie groupoid.
For any  $\alpha\in H^2(\Gamma\upcom,\sS^1)$,
 denote by $K^{vb}_\alpha(\Gamma\upcom)$ the Grothendieck
group of $(\Gamma' ,R' )$-twisted vector bundles, where
$S^1\to R'\to\Gamma'\toto M'$ is any $S^1$-central extension realizing
the class $\alpha$.
\end{defn}

This definition coincides with the definition of twisted
orbifold $K$-theory in the special case considered
by Adem and Ruan \cite{AR}.

Similarly, one can also work in the category of locally compact groupoids
and introduce the $K$-theory group $K^{vb,cont}_\alpha(\Gamma\upcom)$.
We will later find conditions which guarantee that the canonical
morphism $K^{vb}_\alpha(\Gamma\upcom)\to K^{vb,cont}_\alpha
(\Gamma\upcom)$ is an isomorphism (see Theorem~\ref{thm:vb smooth}).

\begin{rmk}
\begin{enumerate}
\item In general,  (iv) does not imply (i).
Even more, there does not exist a function $f\colon \nn\to\nn$
such that every $S^1$-central extension satisfying
(iv) has a twisted vector bundle of rank $\le f(n)$.

Let us assume the contrary. Choose any integer $N>f(n)$
such that $N$ and $n$ are not mutually prime, for instance
$N\equal nf(n)$. We show that
there exists an $S^1$-central extension of Lie groups
$S^1\to R\to\Gamma\toto\cdot$ such that any twisted vector bundle
has rank $\ge N$ and  $R^n$ is trivial.

Let $R'\equal U(N)$,
$\Gamma\equal R'/\zz_n$ where $\zz_n$ is identified with
 the group of $n$-th roots of unity. We consider the central extension
(of order dividing $n$)
$$S^1\to R'\times_{\zz_n}S^1\to \Gamma\toto \cdot.$$

Suppose that there exists a rank $n'$ twisted vector bundle
with $n'\le f(n)$. Then
by (ii) there exists a $\zz_n$-equivariant
group morphism $\pi\colon R'\to GL_{n'}(\cc)$.
Since $R'\equal U(N)$ is compact,
we may assume that $\pi$ is an irreducible unitary representation,
and since $\mbox{dim}\,\pi<N$, we have $\pi\equal (\mbox{det})^p$ for some $p$.


Since $\pi$ is $\zz_n$-equivariant, for $\omega\equal e^{2i\pi/n}$
we get $\pi(\omega r)\equal \omega\pi(r)$ and thus
$\omega^{Np}\equal \omega$. This is impossible since $N$ and $n$ have a
common factor.
%
%

\item   Consider the Lie  group $SL_2({\mathbb R})$.
 Its  fundamental group  is ${\mathbb Z}$.
 Let $H$ be its connected double covering.  Then
 $H$ is  a ${\mathbb Z}_2$-central extension  over $SL_2({\mathbb R}) $.
 Let $R\equal S^1\times_{{\mathbb Z}_2} H $ be its associated
 $ S^1$-central extension over $SL_2({\mathbb R}) $.
Then clearly $R$ defines a torsion class of degree $2$. 

 Let us show that $R\to SL_2({\mathbb R})$ does not admit
any   finite dimensional  twisted vector bundle, i.e.,  a projective
representation.
 It is known (see \cite[p.13]{PS}) that any
 group homomorphism $ \phi$ from  the universal extension $\hat{G}$ 
of $ SL_2({\mathbb R} )$
 to $GL_n({\mathbb C}) $ satisfies the following
property:

 \begin{equation}\label{eq:equivariancedemandee}
 \phi(z\cdot g)\equal \phi(g) , \ \ \ \forall z \in {\mathbb Z}, g \in \hat{G}
 \end{equation}

Assume that $ \psi: R\to GL_n({\mathbb C})$ is an $S^1$-equivariant group 
homomorphism.
Let $ \psi' : H\to GL_n({\mathbb C})$ be its restriction to
$H$. Then $ \psi'$ is a ${\mathbb Z}_2$-equivariant map.
Let $ \pi:\hat{G}\to H$ and $p:{\mathbb Z}\to {\mathbb Z}_2 $
 be the canonical projections. Since $\psi' \smalcirc \pi :
\hat{G}\to GL_n({\mathbb C})$  is  a group homomorphism,
 according to Eq. (\ref{eq:equivariancedemandee}), we have, 
 for any $z \in {\mathbb Z} $ and  $ g \in \hat{G}$, 

$$   \psi' \big( \pi(z\cdot g) \big)\equal  \psi'  \big( p(z) \cdot \pi(g) \big)
 \equal \psi' \big( \pi (g) \big).    $$

Since both $p$ and $ \pi$ are onto, it follows that
for any $a \in {\mathbb Z}_2 $ and any $ g \in H$

\begin{equation}\label{eq:deuxieme} 
 \psi' (a\cdot g)\equal \psi' (g).
  \end{equation}
 This contradicts to the assumption that $ \psi'$ is 
 ${\mathbb Z}_2$-equivariant.
\end{enumerate}
\end{rmk}

\subsection{Proof of step 2}
Recall that a positive element $a$ in a $C^*$-algebra $A$
is said to be strictly positive if
$A\equal \overline {aA}$, which is also equivalent to
 $ A\equal \overline{aAa}$, and that every
separable $C^*$-algebra has a strictly positive element, i.e.,
 $A$ is $\sigma$-unital \cite[3.10.6]{ped79}.

\begin{lem}\label{lem:approx unit for both}
Let $\Gamma$ be a proper groupoid with Haar system
acting on a  u.s.c. field  of  $C^*$-algebras $\aA \to M$, and
$A\equal C_0(M,\aA)$.
 Let $(u_i)\in A_\Gamma$ such that $0\le u_i\le 1$
(see notation~(\ref{eqn:AG})).
The following are equivalent:
\begin{itemize}
\item[(i)] $(u_i)$ is an approximate unit for $A_\Gamma$;
\item[(ii)] $(u_i)$ is an approximate unit for $A$.
\end{itemize}
\end{lem}

(Recall that in (ii) above, we mean that $u_ia\to a$ for all $a\in A$,
but $u_i$ does not necessarily belong to $A$.)

\begin{pf}
(i)$\implies$(ii): it is clear since $\overline{A_\Gamma A}\equal A$
according to  Lemma~\ref{lem:AGammamultiplier}.

(ii)$\implies$(i): by assumption, $u_ia\to a$ for all $a\in A$.
Since $a\mapsto a^\Gamma$ is linear and norm-decreasing,
we have $(u_ia)^\Gamma\to a^\Gamma$. 
On the other hand, it is simple to see that
$(u_ia)^\Gamma\equal u_i a^\Gamma$. Thus,
 from Lemma~\ref{lem:A_Gamma}, $u_ib\to b$ for all $b\in A_\Gamma$.
\end{pf}

\begin{prop}\label{prop:Punital}
Let $\Gamma$ be a proper groupoid with a Haar system
acting on a  u.s.c. field  of  $C^*$-algebras $\aA \to M$, and
$A\equal C_0(M,\aA)$.
 Then (i)--(iii) are equivalent, and (i) $\iff$ (iv) if $A$
is $\sigma$-unital.
\begin{itemize}
\item[(i)] $\exists P_i\in A_\Gamma$ approximate unit of $A_\Gamma$
consisting of projections;
\item[(ii)] $\exists P_i\in A_\Gamma$ approximate unit of $A$
consisting of projections;
\item[(iii)] for all $x\in A_+$ and $\varepsilon>0$,
there exists $a_{\varepsilon,x}\in (A_\Gamma)_+$ such that
$\mbox{sp}(a_{\varepsilon,x})$ does not contain any interval
$[0,\delta]$ ($\delta>0$) and $x
\le \varepsilon+ a_{\varepsilon,x}$;
\item[(iv)] there exists $x\in A_+$ strictly positive with the property
that for all $\varepsilon>0$,
there exists $a_{\varepsilon,x}\in (A_\Gamma)_+$ such that
$\mbox{sp}(a_{\varepsilon,x})$ does not contain any interval
$[0,\delta]$ ($\delta>0$) and $x
\le \varepsilon+ a_{\varepsilon,x}$.
\end{itemize}
\end{prop}

\begin{pf}
(i)$\iff$(ii) follows from Lemma~\ref{lem:approx unit for both}.

(ii)$\implies$(iii):  $ \forall \varepsilon>0$, by (ii), there exists
$i$ such that $\|x-P_ixP_i\|<\varepsilon$. Since
\begin{eqnarray*}
x\equal  (x-P_ixP_i)+P_i xP_i
&\le& \|x-P_ixP_i\|+P_i \|x\|P_i\\
&\equal & \|x-P_ixP_i\|+\|x\| P_i,
\end{eqnarray*}
we see that $a_{\varepsilon,x}\equal \|x\| P_i$ satisfies (iii).

(iii)$\implies$(iv): obvious

(iii)$\implies$(ii): let $x_1,\ldots,x_n\in A$ and $\varepsilon>0$.
We want to find a projection $P\in A_\Gamma$ such that
$\|(1-P)x_i\|\le\varepsilon$ for all $i\equal 1,\ldots,n$.

Let $x\equal \sum x_ix_i^*$. Choose a real number
$\eta$ such that $0<\eta<\varepsilon^2/2$ and $\eta$ does not belong
to the spectrum of $a_{\varepsilon^2/2,x}$. Then the spectral projection
$P\equal 1_{[\eta,\infty)}(a_{\varepsilon^2/2,x})$
of $a_{\varepsilon^2/2,x}$ on $[\eta,\infty)$ is an element of
$A_\Gamma$. Since $1-P$ is the spectral projection of
$a_{\varepsilon^2/2,x}$ on $[0,\eta]$, we have
$(1-P)a_{\varepsilon^2/2,x}(1-P)\le \eta (1-P)$.
Now, for any $i$, we have
\begin{eqnarray*}
(1-P)x_ix_i^*(1-P)&\le& (1-P)x(1-P)\\
&\le& (1-P)(\varepsilon^2/2+a_{\varepsilon^2/2,x})(1-P)\\
&\equal & (1-P)(\varepsilon^2/2+\eta)(1-P)\\
&\le& \varepsilon^2(1-P)\le\varepsilon^2,
\end{eqnarray*}
so $\|(1-P)x_i\|\equal \|(1-P)x_ix_i^*(1-P)\|^{1/2}
\le\varepsilon$ for all $i\equal 1,\ldots,n$.

(iv)$\implies$(ii): the same proof shows that if $x$ satisfies (iv),
there exist projections $P_i\in A_\Gamma$ such that $(1-P_i)x^{1/2}\to 0$,
and therefore $P_iy\to y$ for all $y\in \overline{x^{1/2}A}\equal A$.
\end{pf}

We note  that the approximate unit is not necessarily
increasing. In fact, we have the following:

\begin{prop}\label{prop:countable spectrum}
Let $\Gamma$ be a proper groupoid with a Haar system
acting on a countably generated u.s.c. field of
$C^*$-algebras $\aA\to M$. Let $A\equal C_0(M,\aA)$.
The following are equivalent:

\begin{itemize}
\item[(i)] there exist  projections $P_1\le P_2\le\cdots \le P_n$ 
in $A_\Gamma$ such that $P_i a\to a$ for all $a\in A_\Gamma$;
\item[(ii)] there exist  projections $P_1\le P_2\le\cdots \le P_n$ 
in $A_\Gamma$ such that $P_i a\to a$ for all $a\in A$;
\item[(iii)] there exists $a\in A_\Gamma$ strictly positive with
countable spectrum.
\end{itemize}
\end{prop}

\begin{pf}
(i)$\iff$(ii): follows from Lemma~\ref{lem:approx unit for both}.

(ii)$\implies$(iii): take $a\equal \sum_{n\equal 0}^\infty 2^{-n}(P_{n+1}-P_n)$ (with
$P_0\equal 0$ by convention). Then $1_{[2^{-n},1]}(a)\equal P_{n+1}$, and hence
$P_{n+1}A\subset aA$. It follows that $A\equal \overline{\cup_n P_nA}
\subset \overline{aA}$, so $a$ is strictly positive with spectrum
in $\{0\}\cup\{2^{-n}|\;n\in\nn\}$.

(iii)$\implies$(i): take $P_n\equal 1_{[\alpha_n,\infty)}(a)$,
 where $\alpha_n$
is a sequence decreasing to $0$ and $\alpha_n\notin \mathrm{sp}(a)$.
Then clearly $P_na\to a$. Thus $P_nb\to b$ for all $b\in \overline{aA}
\equal A$.
\end{pf}

\begin{defn}\label{def:GEmodule}
Let $E$ be a u.s.c. Fell bundle over a locally compact groupoid $\Gamma$,
and let $A\equal C_0(M;E)$.
A $(\Gamma,E)$-equivariant Hilbert module is an $A$-Hilbert module
${\mathcal{E}}$ with isomorphisms of $A_{s(g)}$-Hilbert modules
$${\mathcal{E}}_{t(g)}\otimes_{A_{t(g)}} E_g\to {\mathcal{E}}_{s(g)}$$
such that $(\xi\eta)\zeta \equal  \xi(\eta\zeta)$ whenever
$(g,h)\in \Gamma^{(2)}$ and $(\xi,\eta,\zeta)\in {\mathcal{E}}_{t(g)}
\times E_g\times E_h$. The product is required to be continuous in the
following sense: for all $\xi\in {\mathcal{E}}$ and
$\eta\in C_0(\Gamma;E)$,
$g\mapsto \xi(t(g))\eta(g)$ belongs to $s^* \eE$.
\end{defn}

Note that $\eE$ can be canonically identified  with a 
field of Banach   spaces over $M$ such that for any  $x\in M$
the fiber $\eE_x$ at $x$ is an $A_x$-Hilbert module
(see Proposition~\ref{prop: module field}). $\eE$
being $(\Gamma,E)$-equivariant  means, roughly speaking, 
  that this field  of Banach spaces  is equipped with an $E$-action.

\begin{defn}
Let $\Gamma$ be a locally compact groupoid with Haar system,
$E$ a u.s.c. Fell bundle over $\Gamma$,
$A\equal C_0(M;E)$ and ${\mathcal{E}}$ a $(\Gamma,E)$-equivariant
$A$-Hilbert module. Then ${\mathcal{E}}$ is said to be
approximately finitely generated projective (AFGP) if there exist
projections $P_i$ in ${\mathcal{K}}_\Gamma({\mathcal{E}})$
such that $P_i\xi\to\xi$ for all $\xi\in{\mathcal{E}}$.
\end{defn}

For the notation ${\mathcal{K}}_\Gamma({\mathcal{E}})$,
see Eq~(\ref{eqn:K(E)}).

\begin{lem}\label{rem:FGP compact}
If $A$ is a $C^*$-algebra and ${\mathcal{E}}$  an $A$-Hilbert module,
then ${\mathrm{Id}}_{\mathcal{E}}\in{\mathcal{K}} ({\mathcal{E}} )$
implies that ${\mathcal{E}}$
is finitely generated projective, and the converse holds if $A$
is unital. This explains the terminology.
\end{lem}

\begin{pf}
This is proved in the unital case in \cite[Thm 15.4.2,
Rem 15.4.3]{wegge93}. Below we outline a proof for the
direction ``$\implies$''.

If ${\mathrm{Id}}_{\mathcal{E}}$ is compact,
then ${\mathrm{Id}}_{\mathcal{E}}$ can be approximated by finite rank
operators, i.e. there exist $\xi_i$, $\eta_i$ such that
$S\equal \sum_{i\equal 1}^n
T_{\xi_i,\eta_i}$ satisfies $\|{\mathrm{Id}}-S\|<1$. In particular,
$S$ is invertible,  and therefore
$${\mathrm{Id}}_{\mathcal{E}}
\equal S^{-1}S\equal \sum_{i\equal 1}^n T_{S^{-1}\xi_i,\eta_i}.$$
Replacing $\xi_i$ by $S^{-1}\xi_i$, we may assume
that ${\mathrm{Id}}_{\mathcal{E}}\equal \sum_{i\equal 1}^nT_{\xi_i,\eta_i}$.
Now, define
\begin{eqnarray*}
U\colon {\mathcal{E}}&\to& A^n\\
\xi&\mapsto&(\langle\eta_1,\xi\rangle,\ldots,
\langle\eta_n,\xi\rangle)\\
V\colon A^n&\to&{\mathcal{E}}\\
(a_1,\ldots,a_n)&\mapsto&\xi_1a_1+\cdots+\xi_n a_n.
\end{eqnarray*}
 Then $VU\equal {\mathrm{Id}}_{\mathcal{E}}$. Hence $P\equal UV$ is an
idempotent in ${\mathcal{L}}(A^n)\equal M_n(A)$ and
${\mathcal{E}}\cong PA^n$ as right $A$-Hilbert modules.
%
\end{pf}

\begin{prop}\label{prop:AFGP iff stably unital}
Let $\Gamma$ be a proper groupoid with a Haar system,
$E$ a u.s.c. Fell bundle over $\Gamma$,
$A\equal C_0(M;E)$ and ${\mathcal{E}}$ a $(\Gamma,E)$-equivariant
$A$-Hilbert module. Then ${\mathcal{E}}$ is AFGP if and only if
${\mathcal{K}}_\Gamma(\mathcal{E})$ has an approximate unit consisting of
projections.
\end{prop}

\begin{pf}
\begin{eqnarray*}
\forall\xi\in{\mathcal{E}},\;P_n\xi\to \xi&\iff&
\forall\xi\in{\mathcal{E}},
(1-P_n)T_{\xi,\xi}(1-P_n)\equal T_{(1-P_n)\xi,(1-P_n)\xi}\to 0\\
&\iff& \forall\xi,\eta \in{\mathcal{E}},
(1-P_n)T_{\xi,\eta}(1-P_n)\to 0\\
&&\qquad\mbox{ since } T_{\xi,\eta}\equal (1/4)\sum_{\omega^4\equal 1}\omega T_{\xi+\omega\eta,
\xi+\omega\eta}\\
&\iff& \forall T\in{\mathcal{K}}({\mathcal{E}}),\,
(1-P_n)T(1-P_n)\to 0 \\
&\iff& (P_n) \mbox{ is an approximate unit for }
{\mathcal{K}}({\mathcal{E}})\\
&\iff& (P_n) \mbox{ is an approximate unit for }
{\mathcal{K}}_\Gamma({\mathcal{E}})\\
&&\qquad
{\mbox{by  Lemma }}\ref{lem:approx unit for both}.
\end{eqnarray*}
In the second from the last equivalence,
 we used the fact that $\|T(1-P_n)\|
\equal \|(T(1-P_n))^*T(1-P_n)\|^{1/2}
\equal \|(1-P_n)T^*T(1-P_n)\|^{1/2}$.
\end{pf}

An immediate consequence is the following:

\begin{cor}\label{coro:AFGP iff stably unital}
Let $S^1 \to R\to\Gamma\toto M$ be an $S^1$-central extension of proper Lie
groupoids, and let $L\equal R\times_{S^1}\cc$ be the associated
complex  line bundle.  Then $L^2(\Gamma;L)\otimes\HH$ is AFGP if and only if
$C^*_r(\Gamma;R)\otimes{\mathcal{K}} (\HH )$
 has an approximate unit consisting
of projections.
\end{cor}

\begin{pf}
Apply Proposition~\ref{prop:AFGP iff stably unital}
to ${\mathcal{E}}\equal L^2(\Gamma;L)\otimes\HH$ and use
Proposition~\ref{prop:structure of algebra of action}.
\end{pf}

\begin{cor}\label{coro:AFGP iff stably unital2}
Let $\Gamma$ be a proper Lie groupoid, then
$L^2(\Gamma)\otimes\HH$ is AFGP if and only if
$C^*_r(\Gamma)\otimes{\mathcal{K}} (\HH )$ has an approximate unit consisting
of projections.
\end{cor}

\begin{pf}
Apply Corollary~\ref{coro:AFGP iff stably unital}
to the trivial $S^1$-central extension.
\end{pf}

We end this subsection by listing some examples of AFGP modules.

\begin{prop}\label{prop:compact group AFGP}
If $G$ is a compact group and $\pi$ a unitary representation of
$G$ on a separable Hilbert space $\HH$ (considered
as a $G$-equivariant $\cc$-Hilbert module), then
$\HH$ is AFGP. 
\end{prop}

\begin{pf}
Choose a strictly positive element $a\in {\mathcal{K}}({\HH })^G$.
Since $a$ is a compact operator  on the Hilbert space
$\HH $, its spectrum is countable.
By Propositions~\ref{prop:AFGP iff stably unital}
and~\ref{prop:countable spectrum},
$\HH $ is AFGP.
\end{pf}

The following well known result
(see \cite[cor 15.1.4]{dix77} and \cite{koo57})
is a direct consequence of the above
Proposition~\ref{prop:compact group AFGP}.
\begin{cor}\label{coro:finite dimensional}
If $G$ is a compact group,  then
every irreducible unitary representation of $G$ is finite dimensional.
\end{cor}

\begin{pf}
Assume that  $\pi$ is  an irreducible representation on $H_\pi$.
Let $P_n$ be a sequence of compact, $G$-invariant projections
in $H_\pi$ such that
$P_n\xi\to\xi$ for all $\xi$. Since
the representation is irreducible, we have
 either $P_n\equal 0$ or $P_n\equal {\mathrm{Id}}$.
Therefore  $P_n\equal {\mathrm{Id}}$ for $n$ large enough. Since $P_n$
is a compact projection on a Hilbert space, its range
$H_\pi$ is finite dimensional.
\end{pf}

\begin{cor}
\label{coro:M times G}
If $\Gamma$ is  a transformation groupoid $G\times M\toto M$, 
where $M$ is a compact space and $G$ is a
compact group, then $L^2(\Gamma)\otimes\HH$ is AFGP. 
\end{cor}

\begin{pf}
Since $L^2(\Gamma)\otimes\HH\cong C(M)\otimes L^2(G)\otimes\HH$,
the $C^*$-algebra
${\mathcal{K}}(L^2(\Gamma,\HH))\cong
C(M)\otimes {\mathcal{K}}(L^2(G))\otimes{\mathcal{K}} (\HH ) $
is the tensor product of three $C^*$-algebras having approximate units
consisting of invariant projections.
\end{pf}

\subsection{Proof of step 3}

We need a sequence of lemmas.

\begin{lem}\label{lem:Id in C}
Let $M$ be a locally compact space, $F$ a Hermitian vector bundle and
${\mathcal{F}}\equal C_0(M,F)$ its space of continuous sections vanishing
at infinity considered as a $C_0(M)$-Hilbert module. Then
${\mathrm{Id}}_{\mathcal{F}}\in{\mathcal{C}}({\mathcal{F}})$
(see notation~(\ref{eqn:C(E)})).
\end{lem}

\begin{pf}
For every compact subspace $K$ of $M$,
the restriction of ${\mathcal{F}}$ to $K$, i.e. the $C(K)$-Hilbert
module ${\mathcal{F}}_K\equal {\mathcal{F}}\otimes_{C_0(M)}C(K)$,
is the space of sections of $F_{|K}$, thus by Swan theorem,
is a projective finitely generated $C(K)$-module. Therefore,
from Remark~\ref{rem:FGP compact},
the identity map on ${\mathcal{F}}_K$ is compact.

Let us show that this implies
${\mathrm{Id}}_{\mathcal{F}}\in{\mathcal{C}}({\mathcal{F}})$.
Let $\varphi$, $\psi\in C_c(M)$. Choose an open set $U$
and a compact set $K$ such that $U\subset K\subset M$ and
$U$ contains the supports of both  $\varphi$ and  $\psi$.
Since $\mathrm{Id}_{{\mathcal{F}}_K}$ is compact, there exist
$\xi_i$, $\eta_i\in {\mathcal{F}}_K$ such that
$\mathrm{Id}_{{\mathcal{F}}_K}\equal \sum_i T_{\xi_i,\eta_i}$.
Let $\xi'_i\equal \overline{\varphi}\xi_i$ and $\eta'_i\equal \psi\eta_i$. Then
$\varphi\psi\equal \sum_i T_{\xi'_i,\eta'_i}\in {\mathcal{L}}({\mathcal{F}}_K)$,
where $\varphi\psi$ denotes the multiplication
 operator acting on ${\mathcal{F}}$.
However since $\varphi$, $\psi$ and all $\xi'_i$, $\eta'_i$
are all supported in $U$, it is not hard to check that
the equality
$\varphi\psi\equal \sum_i T_{\xi'_i,\eta'_i}$
also holds in ${\mathcal{L}}({\mathcal{F}})$.
Therefore, $\varphi\psi$ is compact for all $\varphi$, $\psi\in
C_c(M)$. By a density argument, $\varphi$ is compact
for all $\varphi\in C_0(M)$, i.e.
$\mathrm{Id}_{{\mathcal{F}}}\in {\mathcal{C}}({\mathcal{F}})$.
\end{pf}

\begin{lem}\label{lem:tensor AFGP}
Let $S^1\to R_i\to\Gamma\toto M$,  $i\equal 1,2$,
 be $S^1$-central extensions of the Lie groupoid $\Gamma$, and 
${\mathcal{E}}_i$, $i\equal 1,2$,   
$(\Gamma,R_i)$-equivariant $C_0(M)$-Hilbert modules. Suppose that
${\mathcal{E}}_1$ is AFGP and  ${\mathrm{Id}}_{{\mathcal{E}}_2}
\in {\mathcal{C}}({\mathcal{E}}_2)$. Then
${\mathcal{E}}':\equal {\mathcal{E}}_1\otimes_{C_0(M)}{\mathcal{E}}_2$ is AFGP
as a $(\gm , R_1 \otimes R_2 )$-equivariant $C_0(M)$-Hilbert module.
\end{lem}

\begin{pf}
By assumption, there exists an approximate unit
$P_n\in{\mathcal{K}}_\Gamma({\mathcal{E}}_1)$ consisting of projections.
Let $P'_n\equal P_n\otimes_{C_0(M)}{\mathrm{Id}}_{\mathcal{E}_2}$.
It is clear that $P'_n$ is an invariant projection, and that
$\|P'_n(x)\|\to 0$ when $x\to 0$ in $M/\Gamma$. Let us show that
$P'_n\in{\mathcal{C}}({\mathcal{E}}')$.
For all $\varphi$, $\psi\in C_0(M)$, $(\varphi\psi)P'_n
\equal (\varphi P_n)\otimes_{C_0(M)}\psi\in {\mathcal{K}}
({\mathcal{E}}')$. It follows that $\phi P'_n
\in {\mathcal{K}}({\mathcal{E}}')$ for all $\phi\in C_0(M)$,
i.e. $P'_n\in {\mathcal{C}}({\mathcal{E}}')$.
Therefore $P'_n\in {\mathcal{K}}_\Gamma({\mathcal{E}}')$,
and it is clear that $P'_n$ is an approximate unit consisting of
projections.
\end{pf}

Before we proceed, we need to introduce some notation.

Let $E$ be a u.s.c. Fell bundle over the groupoid $\Gamma$,
$A\equal C_0(M;E)$, and let
${\mathcal{E}}$ be a (possibly non-equivariant) $A$-Hilbert module.
Consider
the field of Banach spaces over $\Gamma$
with fiber ${\mathcal{E}}_{t(g)}\otimes_{A_{t(g)}}E_g$, determined
by sections of the form $\eta(t(g))\otimes \zeta(g)$ where
$\eta\in {\mathcal{E}}$ and $\zeta\in C_0(\Gamma;E)$
(see Proposition~\ref{prop:section field}).
Denote by $C_c(\Gamma;E,{\mathcal{E}})$ the space of compactly supported
sections. 

\par\medskip
Endow $C_c(\Gamma;E,{\mathcal{E}})$ with an  $A$-valued scalar product
$$\langle\xi,\eta\rangle (x)\equal \int_{g\in \Gamma_x}
\langle \xi(g),\eta(g)\rangle\,\lambda_x(dg), \ \ \ \forall x\in M,$$
 and
denote by $L^2(\Gamma;E,{\mathcal{E}})$ its completion. Since
$L^2(\Gamma;E,{\mathcal{E}})$ is an $A$-Hilbert module, it can be
considered as a field of Banach spaces over $M$;
denote by $L^2(\Gamma_x;E,{\mathcal{E}})$ its fiber at $x$,
which  is an $A_x$-Hilbert module.

The usual action of $(\Gamma,E)$  on $L^2(\Gamma_x;E )$ (see Appendix)
extends naturally to an action on $L^2(\Gamma;E,{\mathcal{E}})$, which
is defined as follows:

$$L^2(\Gamma_x;E,{\mathcal{E}})\otimes_{A_x}E_{\gamma^{-1}}\cong
L^2(\Gamma_y;E,{\mathcal{E}})$$
$$\xi\otimes\eta\mapsto ( g\mapsto\xi(g\gamma)\otimes\eta )\in
{\mathcal{E}}_{t(g)}\otimes E_{g\gamma}\otimes E_{\gamma^{-1}}
\cong {\mathcal{E}}_{t(g)}\otimes E_g .$$

If moreover  ${\mathcal{E}}$ is a $(\Gamma,E)$-equivariant
module, then since ${\mathcal{E}}_{t(g)}\otimes_{A_{t(g)}}E_g
\cong {\mathcal{E}}_{s(g)}$, we get
\begin{equation}
\label{equ}
L^2(\Gamma;E, {\mathcal{E}})\cong L^2(\Gamma)\otimes_{C_0(M)}
{\mathcal{E}},
\end{equation}
where on the left-hand side we \emph{forget} the equivariant
structure on ${\mathcal{E}}$, while the right-hand side
is endowed with a ``diagonal'' action. The isomorphism above
is well-known in the case of a Hilbert space
$\HH$ endowed with a unitary representation 
$g\mapsto U_g$ of a locally compact group $G$. In this case,
the isomorphism $L^2(G, \HH )\to L^2(G)\otimes \HH $
(where $\HH $ is endowed with the \emph{trivial} action of $G$
on the left-hand side but $G$ acts \emph{diagonally} on the right-hand side)
is given by $$\xi\mapsto \eta(g)\equal U_g(\xi(g)).$$

\begin{prop}\label{prop:L2(Gamma,E) is AFGP}
Assume  that $\Gamma \toto M$ is a proper Lie groupoid such that
$M/\Gamma$ is compact and $S^1\to R\to
\Gamma$ is an $S^1$-central extension. Let $L\equal R\times_{S^1}\cc$.
If $L^2(\Gamma)\otimes\hh$ is AFGP,
and if there exists a topological (i.e. without
differentiable structure) $(\Gamma,R)$-twisted vector bundle (of finite
rank), then $L^2(\Gamma;L)\otimes\HH$ is AFGP.
\end{prop}

\begin{pf}
Let $F$ be a $(\Gamma, R)$-twisted vector bundle. Since $\Gamma$ is
proper, $F$ can be endowed with an invariant Hermitian metric, 
and therefore can be considered as a $(\Gamma,R)$-equivariant Hilbert 
module.  As a (non-equivariant) continuous field of Hilbert spaces
 over $M$,
$F\times\HH\to M$ is locally trivial with infinite dimensional
 fibers.
According to the  triviality theorem of Dixmier and Douady
\cite{dix-dou63}, $F\times\HH\to M$  is isomorphic to
$M\times\HH \to M$.
Moreover, the space of continuous sections $\fF\equal C_{0} (M, F)$ of $F\to M$
 can be considered as a $(\Gamma,R)$-equivariant
$C_0(M)$-module
such that ${\mathrm{Id}}_\fF\in {\mathcal{K}}_\Gamma(\fF)$
(see Lemma~\ref{lem:Id in C}).

Since $L^2(\Gamma)\otimes\HH$ is AFGP as a 
$\Gamma$-equivariant Hilbert module,
according to  Lemma~\ref{lem:tensor AFGP} and Lemma~\ref{lem:Id in C},
 we see that
$L^2(\Gamma)\otimes\HH\otimes_{C_0(M)} \fF$
is AFGP as a $(\Gamma,R)$-equivariant Hilbert module.
Using the isomorphism (\ref{equ}), we   deduce that 
$L^2(\Gamma;L,\HH\otimes \fF)$ is AFGP.
By  the triviality of the Hilbert bundle $F\times \HH$, we get
$L^2(\Gamma;L,\HH\otimes \fF)
\cong L^2(\Gamma;L,C_0(M)\otimes \HH)
\cong L^2(\Gamma;L)\otimes \HH$.
Therefore $L^2(\Gamma;L)\otimes \HH$ is AFGP.
\end{pf}

\subsection{Proof of step 4}
\begin{prop}\label{prop:stabilization}[Stabilization theorem]
Let $\Gamma\toto M$ be a proper groupoid with a  Haar system, and $E$
 a u.s.c. Fell bundle over $\Gamma$. Let   $A\equal C_0(M;E)$.
Assume that ${\mathcal{E}}$ is a $(\Gamma,E)$-equivariant countably generated
$A$-Hilbert module. Then we have the following
  isomorphism  of $(\Gamma, E)$-equivariant
Hilbert $C^*$-modules:
$${\mathcal{E}}\oplus L^2(\Gamma;E)\otimes\HH
\cong L^2(\Gamma;E)\otimes\HH.$$
\end{prop}

\begin{pf}
Since $\Gamma$ is proper, $C_0(M)$ is (as a $\Gamma$-equivariant
$C_0(M)$-Hilbert module) a direct factor of
$L^2(\Gamma)$ \cite{tu99}.
Hence
${\mathcal{E}}$ is a direct factor of $L^2(\Gamma)\otimes_{C_0(M)}
{\mathcal{E}}\cong L^2(\Gamma;E,{\mathcal{E}})$.
By Kasparov's stabilization theorem for non equivariant modules
(\cite[Thm 15.4.6]{wegge93}),
${\mathcal{E}}$ is a direct factor of $A\otimes\hh$, and  thus
${\mathcal{E}}$ is a direct factor of
$L^2(\Gamma;E,A\otimes\hh)\cong L^2(\Gamma;E)\otimes\hh$.  That is,
there exists ${\mathcal{E}}'$ such that
${\mathcal{E}}\oplus {\mathcal{E}}'\cong L^2(\Gamma;E)\otimes\HH$.
Therefore, we have
$$L^2(\Gamma;E)\otimes\HH\cong L^2(\Gamma;E)\otimes
(\HH\oplus\HH\oplus\cdots)
\cong {\mathcal{E}}\oplus {\mathcal{E}}'\oplus
{\mathcal{E}}\oplus {\mathcal{E}}'\oplus\cdots
\cong {\mathcal{E}}\oplus L^2(\Gamma;E)\otimes\HH .$$
\end{pf}

\begin{cor}\label{coro:grothprojections}
Let $\Gamma\toto M$ be a proper Lie groupoid,
 $S^1\to R\to\Gamma\toto M$  an $S^1$-central extension,
and $L\equal R\times_{S^1}\cc$.
Then there is an equivalence of categories between
the category of $(\Gamma,R)$-equivariant $C_0(M)$-Hilbert modules
${\mathcal{E}}$ such that ${\mathrm{Id}}_{\mathcal{E}}\in
{\mathcal{K}}_\Gamma({\mathcal{E}})$ and the category of projections in
$C^*_r(\Gamma;R)\otimes {\mathcal{K}} (\HH )$.
\end{cor}

\begin{pf}
As usual, let $L\equal R\times_{S^1}\cc$.
Recall from Proposition~\ref{prop:structure of algebra of action}
that $C^*_r(\Gamma;R)\otimes{\mathcal{K}} (\HH )$ is isomorphic
to ${\mathcal{K}}_\Gamma(L^2(\Gamma;L)\otimes\HH)$.
Given  a projection $P\in {\mathcal{K}}_\Gamma(L^2(\Gamma;L)\otimes\HH)$,
${\mathcal{E}}\equal P(L^2(\Gamma;L)\otimes\HH)$ is
 a $(\Gamma,R)$-equivariant Hilbert module.
It is clear that
${\mathrm{Id}}_{\mathcal{E}}\in {\mathcal{K}}_\Gamma({\mathcal{E}})$.

Conversely, if ${\mathcal{E}}$ is a countably generated
$(\Gamma,R)$-equivariant Hilbert module such that
${\mathrm{Id}}_{\mathcal{E}}\in {\mathcal{K}}_\Gamma({\mathcal{E}})$,
we know, from the stabilization theorem~\ref{prop:stabilization},
that there is an
invariant projection $P$ such that  
${\mathcal{E}}\equal P(L^2(\Gamma;L)\otimes\HH)$.
Since ${\mathrm{Id}}_{\mathcal{E}}\in
{\mathcal{K}}_\Gamma({\mathcal{E}})$, we have
$P\in {\mathcal{K}}_\Gamma(L^2(\Gamma;L)\otimes\HH)$.

A standard argument shows that two projections $P_1$
and $P_2\in {\mathcal{K}}_\Gamma(L^2(\Gamma;L)\otimes\HH)$
are Murray-Von Neumann equivalent if and only if the
associated Hilbert modules are isomorphic.
\end{pf}

\subsection{Proof of step 5}

The next proposition generalizes Serre-Swan theorem:
if $M$ is a compact space,  there is an equivalence of categories
between vector bundles on $M$ and finitely generated projective
$C(M)$-modules  (and thus $K^0(M)\equal K_0(C(M))$).

\begin{prop}\label{prop:swan}
Assume that $\Gamma \toto M$ is a proper Lie groupoid such that $M/\Gamma$
is compact, and $S^1\to R\to\Gamma\toto M$ is  an $S^1$-central extension.
\begin{itemize}
\item[(a)] 
The forgetful functor from the category of topological
(i.e. without differentiable structure)
$(\Gamma, R)$-twisted vector bundles
endowed with an $R$-invariant Hermitian  metric to the category of
$(\Gamma,R)$-twisted vector bundles is an equivalence of categories;
\item[(b)] The functor
from the category of topological $(\Gamma,R)$-twisted vector bundles
endowed with an $R$-invariant   metric
to the category of $(\Gamma,R)$-equivariant
$C_0(M)$-Hilbert modules ${\mathcal{E}}$
such that ${\mathrm{Id}}_{\mathcal{E}}
\in{\mathcal{K}}_\Gamma({\mathcal{E}})$, defined by
$$\Phi\colon
F\mapsto C_0(M,F),$$
is an equivalence of categories.
\end{itemize}
\end{prop}

\begin{pf}
To prove (a), note that by an averaging procedure using
cutoff functions (see Definition~\ref{defi:cutoff}),
every twisted vector bundle can be endowed with an invariant Hermitian  metric.
If two Hermitian  $R$-equivariant 
vector bundles $F_1$ and $F_2$ are  isomorphic as
$R$-equivariant vector bundles,
then by the polar decomposition,
they must be isometrically isomorphic.
Indeed, if $T_x\colon F_{1,x}\to F_{2,x}$ is an
$R$-equivariant isomorphism, then $U_x : \equal T_x(T_x^*T_x)^{-1/2}$,
 $\forall {x\in M}$, defines an $R$-equivariant isometric isomorphism.
\par\medskip

Let us prove (b). From Lemma~\ref{lem:Id in C} and the fact that
$M/\Gamma$ is compact, it is easy to see that
${\mathcal{E}}:\equal \Phi(F)$ satisfies ${\mathrm{Id}}_{\mathcal{E}}
\in{\mathcal{K}}_\Gamma({\mathcal{E}})$, and
can be endowed with a  $(\Gamma,R)$-action  so that
$\Phi$ is equivariant.
Therefore $\Phi$ is well-defined and functorial.
\par\medskip

Now, for every $(\Gamma,R)$-equivariant
$C_0(M)$-Hilbert module ${\mathcal{E}}$
satisfying
 ${\mathrm{Id}}_{\mathcal{E}}\in {\mathcal{K}}_\Gamma({\mathcal{E}})$,
${\mathcal{E}}$ is isomorphic to $C_0(M,F)$, 
where $F$ is a continuous field of Hilbert spaces on $M$
with fiber $F_x\equal {\mathcal{E}}\otimes_{ev_x}\cc$ (see Appendix).
For every compact $K\subset M$,  since ${\mathcal{E}}_K$ is a
 finitely generated projective module over $C(K)$, it follows from Swan
 theorem that $F_{|K}$ is a vector bundle
over $K$ (i.e.,  a locally trivial field of finite dimensional
Hilbert spaces). Since this is true for every compact $K$,  
it follows that
$F$ is a vector bundle.

Define $\Psi({\mathcal{E}})\equal F$. It is clear
that $\Phi$ and $\Psi$ are inverse from each other.
\end{pf}

\subsection{The main theorem: continuous case}\label{subsec:ccase}
\begin{them}\label{thm:vector bundles}
Let $\Gamma\toto M$ be a 
Lie groupoid, and $S^1\to R\to\Gamma\toto M$  an $S^1$-central extension.
Denote by $\alpha$ its  corresponding class in $H^2(\Gamma\upcom,\sS^1)$.
Assume  that
\begin{itemize}
\item[(a)] $\Gamma\toto M$ is proper;
\item[(b)] $M/\Gamma$ is compact;
\item[(c)] $L^2(\Gamma)\otimes\hh$ is AFGP. In other words,
there exists a sequence $(P_n)$ such that
\begin{itemize}
\item[(i)] $P_n\equal (P_n(x))_{x\in M}$ is a continuous section of
the field of compact operators
${\mathcal{K}}(\tilde{L}^2(\Gamma)\otimes\hh)\to M$;
\item[(ii)] $x\mapsto P_n(x)$ is $\Gamma$-equivariant;
\item[(iii)] $P_n(x)$ is a finite rank projection for all $x$;
\item[(iv)] for every compactly supported continuous section
$\xi$ of $\tilde{L}^2(\Gamma)\otimes\hh$, $(P_n\xi)(x)
{\to} \xi(x)$
uniformly on $M$ when $n\to\infty$.
\end{itemize}
\item[(d)] there exists a $(\Gamma,R)$-twisted vector bundle
(of finite rank).
\end{itemize}
Then $K^0_\alpha(\Gamma\upcom)$ is isomorphic to
$K^{vb,cont}_\alpha(\Gamma\upcom)$.
\end{them}

\begin{pf}
It is known that if $B$ is a stably unital $C^*$-algebra,
i.e., $B\otimes{\mathcal{K}} (\HH )$ has an approximate 
identity consisting of projections, then $K_0(B)$ is
 the Grothendieck group
of projections in $B\otimes {\mathcal{K}} (\HH )$
\cite[Prop 5.5.5]{black98}. We want to apply this fact to
$B\equal C^*_r(\Gamma;L)$.

Since $L^2(\Gamma;L)\otimes\HH$ is AFGP
according to Proposition~\ref{prop:L2(Gamma,E) is AFGP}, it follows
from Proposition~\ref{prop:structure of algebra of action}
and Proposition~\ref{prop:AFGP iff stably unital} that 
$C^*_r(\Gamma;R)$ is stably unital. Hence $K_0(C^*_r(\Gamma;R))$
is the Grothendieck
group of projections in $C^*_r(\Gamma;R)\otimes{\mathcal{K}} (\HH )$.
Therefore it is  the Grothendieck group  of
$(\Gamma,R)$-twisted vector bundles  according to
Corollary~\ref{coro:grothprojections} and Proposition~\ref{prop:swan}.
\end{pf}

\begin{rmk}\label{rem:abcd}
Note that  Conditions (a), (b), (c) and (d) are invariant 
under   Morita equivalence according to Lemma~\ref{lem:proper morita},
 Corollary~\ref{coro:AFGP iff stably unital2} and
Proposition~\ref{prop:exists twisted}.
\end{rmk}

We end this subsection by describing  an
explicit isomorphism $K^0_\alpha(\Gamma\upcom)\to K^{vb,cont}_\alpha(\Gamma\upcom)$.
We will use the Fredholm picture for $K^0_\alpha(\Gamma\upcom)$
(Theorem~\ref{thm:K theory Fredholm}).

Let $T\in{\mathcal{F}}^0_\alpha$.  By definition, there exists
$S\in{\mathcal{F}}^0_\alpha$ such that
$ST\equal 1+K$ where $K\in {\mathcal{K}}_\Gamma(
L^2(\Gamma;L)\otimes\HH)$.
By  Lemma~\ref{lem:A_Gamma}, we have
$K\equal K_0^\Gamma+K_1^\Gamma$
where $K_0\equal \sum_{i\equal 1}^n T_{\xi_i,\eta_i}$ and
$\|K_1\|<1$  is compact.
Since $L^2(\Gamma;L)\otimes\HH$ is
AFGP, we may assume that $\eta_i\in P(L^2(\Gamma;L)\otimes\HH)$,
where $P\in {\mathcal{K}}_\Gamma(L^2(\Gamma;L)\otimes\HH)$
is a projection. Then
\begin{eqnarray*}
(1+K_1^\Gamma)^{-1}ST(1-P)
&\equal &(1+K_1^\Gamma)^{-1}(1+K_1^\Gamma+K_0^\Gamma)(1-P)\\
&\equal &1-P+((1+K_1^\Gamma)^{-1}K_0(1-P))^\Gamma\\
&\equal &1-P+(\sum_i T_{(1+K_1^\Gamma)^{-1}\xi_i,(1-P)\eta_i})^\Gamma\\
&\equal &1-P
\end{eqnarray*}
Replacing $T$ by $T(1-P)$ and $S$ by $(1+K_1^\Gamma)^{-1}S$,
we may assume that $ST$  equals to the projection $1-P$.
Thus $TS$ is also a projection. Let $Q\equal 1-TS$.
Then the image of $[T]\in K^0_\alpha(\Gamma\upcom)$ is
$[P]-[Q]\in K_\alpha^{vb,cont}(\Gamma\upcom)$.
\par\medskip

Conversely, assume that $P$ is a projection
in ${\mathcal{K}}_\Gamma(L^2(\Gamma;L)\otimes\HH)$.
Let ${\mathcal{E}}\equal P(L^2(\Gamma;L)\otimes\HH)$, and 
$$T\colon L^2(\Gamma;L)\otimes\HH\cong {\mathcal{E}}
\oplus L^2(\Gamma;L)\otimes \HH\to L^2(\Gamma;L)\otimes \HH$$
be the projection.
More explicitly,
$$T\colon L^2(\Gamma;L)\otimes\HH\otimes \ell^2(\nn)
\to L^2(\Gamma;L)\otimes\HH\otimes
\ell^2(\nn)$$
$$(\xi_n)_{n\ge 0}\mapsto (P\xi_{n+1}+(1-P)\xi_n)_{n\ge 0}.$$
Then the map $K^{vb,cont}_\alpha(\Gamma\upcom)\to K^0_\alpha(\Gamma\upcom)$
 is given by $[P]\to [T]$.

\subsection{Discussion on the conditions in Theorem~\ref{thm:vector bundles}}

We would like to remark that Conditions (a)-(d) are all necessary
for Theorem~\ref{thm:vector bundles} to hold. Let us
go over these conditions one by one.

1) Condition (a) cannot be avoided even when
$\Gamma$ is a group $G$. Note that  the $K_0$-group of $C^*_r(G)$
in general  is not equal to the  (finite dimensional)
 representation ring  of  $G$ when $G$  is not compact.
\par\medskip

2) When $\Gamma$ is the manifold $\rr$,
$K_0(C^*_r(\Gamma))\equal \{0\}$ while vector bundles
on $\rr$ are obviously classified by their rank. Thus, condition (b)
cannot be removed.
\par\medskip

3) Condition (c) is not always true for every proper Lie groupoid.
For instance, let $G\equal SU(2)$, and $\gm$ be
 the transformation groupoid $G\times G \toto G$, 
where $G$ acts on itself by conjugation.
 It is known that $H^3(\gm \upcom ,\zz)\equal H^3_G (G, \zz)\equal \zz$
\cite{Mein}.
Let $S^1\to R\to\Gamma\toto M$ be an $S^1$-central extension
corresponding to the generator of $H^3_G (G, \zz)$.
Then $R$ is clearly a proper Lie groupoid, and
$C^*_r(R)\equal \oplus_{n\in\zz} C^*_r(\Gamma;R^n)$
according to Proposition~\ref{prop:direct sum}. Assume
that  $C^*_r(R)$ is  stably unital. Then $C^*_r(\Gamma;R)$
is  stably unital since a quotient of a stably unital
$C^*$-algebra is obviously stably unital.
 Therefore it follows that  there  exists a projection in 
$C^*_r(\Gamma;R)\otimes{\mathcal{K}} (\HH)$, and hence a
$(\Gamma,L)$-twisted vector bundle by
Corollary~\ref{coro:grothprojections}. This contradicts
Proposition~\ref{prop:exists twisted}.
In fact the above argument shows that (c) fails for any non-torsion
$S^1$-central extension of a proper Lie  groupoid.

However, note that condition (c) is fulfilled when $\Gamma$ is
a transformation groupoid $G\times M\toto M$,
where $G$ is a compact Lie  group acting on a compact manifold $M$
(Corollary~\ref{coro:M times G}), or when $\gm$ is  a compact \'{e}tale
groupoid (since in this case $C^*_r(\Gamma)$ is unital).
\par\medskip

4) Condition (d) implies that the class $\alpha$ of the $S^1$-central extension in
$H^2(\Gamma\upcom,{\mathcal{S}}^1)$ must be a  torsion.
 We conjecture that the converse holds:

\begin{quote}
\label{conj:twisted vb}
{\bf Conjecture }\ Let $ \Gamma\toto M$ be a proper Lie  groupoid such that $M/\gm$ is compact. 
Assume that  $S^1\to R\to\Gamma\toto M$ is an $S^1$-central
 extension of Lie groupoids
which corresponds to a  torsion  class in $H^2(\Gamma\upcom ,{\mathcal{S}}^1)$.
 Then there exists a $(\gm, R)$-twisted vector bundle.
\end{quote}

It is known by  Serre-Grothendieck theorem~\cite{don-kar70, Groth}
that Conjecture~\ref{conj:twisted vb} holds when $\Gamma$ is
Morita equivalent to a compact manifold.  It also holds if
$\Gamma$ is a compact group.  In this case,  $R$ is also a  compact group, so
$C^*_r(R)$ is stably unital (see,  for instance, Corollary~\ref{coro:AFGP iff stably unital}
 and   Proposition~\ref{prop:compact group AFGP}). Therefore there always exists a
twisted vector bundle, i.e., a finite dimensional projective representation.          
However, the conjecture remains open even for orbifold groupoids (i.e. \'{e}tale
proper groupoids). 

One possibility to prove this conjecture is to generalize Grothendieck's
proof \cite[Theorem 1.6]{Groth} to the simplicial CW-complex $\gm\upcom$
corresponding to the groupoid $\gm$. This requires some sophisticated
study of homotopy theory of  simplicial manifolds. In particular,
the following question arises naturally:

\begin{quote}
{\bf Question}  Let $PU(\infty )$ be the inductive limit of $PU(n)$,
and $\gm\toto M$ a proper Lie groupoid such that $M/\gm$ is compact.
Let $\gamma $ be an element in   $H^2 (\gm\upcom,  {\mathbb Q}/\zz ) $.
Does  $\gamma$ always induce a map from  simplicial manifolds $\tilde{\gm}\upcom$ to 
$PU(\infty )\upcom$, where $\tilde{\gm}$ is some Lie groupoid Morita equivalent to
$\gm$?
\end{quote}

Finally, we list some consequences of   Theorem~\ref{thm:vector bundles}
in various special cases.

\begin{cor}
Let $M$ be a compact manifold and $\alpha$ a torsion class in
$H^3(M,\zz)$. Then $K^0_\alpha(M)$ is isomorphic to
$K^{vb,cont}_\alpha(M)$.
\end{cor}
 
\begin{pf}
Use Remark~\ref{rem:abcd} and the discussion  following the
conjecture above.
\end{pf}

\begin{cor}
\label{cor:twisted-eq}
Let $M$ be a compact manifold and  $G$  a compact Lie group.
Assume that  $\alpha\in H^3_G (M, \zz )$ is a torsion class 
which admits at least one twisted vector bundle.
 Then $K^0_{G,\alpha}(M)$ is isomorphic
to $K^{vb,cont}_\alpha(\Gamma\upcom)$, where $\Gamma$ is
the transformation groupoid $G\times M\toto M$.
\end{cor}

Note that in the above case when $\alpha\equal 0$,
 twisted vector bundles
 simply correspond to  $G$-equivariant vector bundles over $M$, which
always exist. 
 Corollary \ref{cor:twisted-eq} simply
implies that the original definition of equivariant
$K$-theory of Segal \cite{Segal} is equivalent to the
$K$-theory of the crossed product $C^*$-algebra $C_0(M)\rtimes G$.

\begin{cor}
Let ${\mathbf{X}}$ be a compact orbifold. Assume 
 that ${\mathbf{X}}$ is reduced, or that ${\mathbf{X}}$ can
 be represented by a compact \'etale groupoid. If
$\alpha \in H^3 ({\mathbf{X}}, \zz)$
is a torsion which admits a twisted vector bundle, then
$K_\alpha^0 ({\mathbf{X}} )$ is isomorphic to
$K_\alpha^{vb,cont} ({\mathbf{X}})$.
\end{cor}
\begin{pf}
Recall that if an orbifold is reduced, it can be represented by a
crossed-product of a manifold by a compact group, and therefore
the result follows from Corollary~\ref{cor:twisted-eq}.

If $\Gamma$ is a compact \'etale groupoid, then $C^*_r(\Gamma)$ is
unital (the unit being represented by the characteristic function
of $\Gamma^{(0)}$), and therefore  Condition (c) is fulfilled.
\end{pf}

\subsection{The main theorem: smooth case}

Our goal in this subsection  is to prove the analogue of
Theorem~\ref{thm:vector bundles} for smooth vector bundles.
The main result is the following

\begin{them}\label{thm:vb smooth}
Under the same hypothesis as in Theorem~\ref{thm:vector bundles},
we have the following  commutative diagram of isomorphisms:
\begin{equation}\label{eqn:square iso}
\xymatrix{
K_0(C_c^\infty(\Gamma,R))\ar[r]^{V}\ar[d]^i
    & K_\alpha^{vb}(\Gamma\upcom )\ar[d]^{i'}\\
K_0(C^*_r(\Gamma;R))\ar[r]^{V'}
    & K_\alpha^{vb,cont}(\Gamma\upcom),
}
\end{equation}
where $i$ and $i'$ are naturally defined;  $V$ and $V'$
are defined as follows. For every projection
$P\in C^*_r(\Gamma;R)\otimes\kK(\HH)$, $V'(P)\equal P(\tilde{L}^2(\Gamma;L)
\otimes\HH)$, and for every projection 
$P\in C_c^\infty(\Gamma;R) \otimes \kK( \HH_n )$,
 $V(P)\equal P(\tilde{L}^2(\Gamma;L) \otimes\HH_n )$.
Here $\HH_n$ denotes  the $n$-dimensional  Hilbert space
 $\cc^n\subset \HH$.
\end{them}

It follows from Theorem~\ref{thm:vector bundles} that $V'$
is well-defined and is an isomorphism. 
To prove that $i$ and $i'$ are isomorphisms, we will first show
that $C_c^\infty(\Gamma,R)$ is stable under holomorphic functional
calculus. 
Let us recall the definition below.

\begin{defn}
Assume that $\bB$ is a subalgebra of a Banach algebra $B$.
Let $\tilde{\bB}$ and $\tilde{B}$ be the unitization of $\bB$
and $B$ respectively.  $\bB$ is said to be stable under holomorphic
functional calculus if for any $b\in \tilde{\bB}$
and any  $f$ holomorphic on a neighborhood of $sp(b)$, we have
$f(b)\in \tilde{\bB}$.
\end{defn}                     

If furthermore $\bB$ is endowed with a structure of
Fr\'echet algebra such that the inclusion
$\bB\to B$ is continuous, then  the following are equivalent
(see \cite[Appendix]{bost} or
\cite[Lemma~1.2, Thm~2.1]{sch}):
\begin{itemize}
\item[(i)]
$\bB$ is stable
under holomorphic calculus;
\item[(ii)] for all $n$, $M_n(\bB)$ is
stable under holomorphic calculus;
\item[(iii)] every element in $\tilde{\bB}$ which
is invertible in $\tilde{B}$ is actually invertible in $\tilde{\bB}$.
\end{itemize}
When any of the  conditions above is satisfied,
 the inclusion $\bB\to B$ induces an isomorphism  of
$K$-theory.
\par\medskip

Assume now that $\Gamma$ and $M$ are manifolds and  $s\colon \Gamma\to M$
is a submersion ($\Gamma$ is not necessarily a groupoid).
Let $E \to \gm $ be a Hermitian vector bundle. 
Assume  that there exists a smooth $s$-system
$\mu\equal (\mu_x)_{x\in M}$, i.e. $\mu_x$ is a measure on $\Gamma$
whose support is $\Gamma_x\equal s^{-1}(x)$ such that for every $f\in
C_c^\infty(\Gamma)$ the function
$x\mapsto \int_{g\in \Gamma_x}f(g)\mu_x(dg)$ is smooth.

\begin{rmk}
We will be interested in the case that $S^1\to R\to\Gamma\toto M$
is an $S^1$-central extension of Lie groupoids, $s$ is the source map,
$E$ is the associated line bundle and $\mu$ is a smooth Haar system.
It is well-known that  such a Haar system exists on any Lie
groupoid.
\end{rmk}

Let $F\subset \Gamma$ be a closed subset such that the restriction
 $s_{|F}\colon F\to M$ is proper.
Let $\aA_F\equal \{a\in C^\infty(\Gamma\times_s \Gamma,
pr_1^*(E)\otimes pr_2^*(E^*))|\; \mbox{supp}(a)\subset F\times_s F\}$,
where $pr_1, pr_2: \Gamma\times_s \Gamma \to \gm$ are the projections.
We endow $\aA_F$ with the convolution product
$$(a*b)(g,h)\equal \int_{\Gamma_{s(g)}}a(g,k)\cdot b(k,h)\,\mu_{s(g)}(dk), $$
where $\cdot$ denotes the obvious product $E_g\otimes E_k^*\otimes
E_k\otimes E_h^*\to E_g\otimes E_h^*$,
and the adjoint
$$(a^*)(g,h)\equal a(h,g)^*.$$
For any $\xi\in C_c^\infty(\Gamma_x;E)$, let
$$(\pi_x(a)(\xi))(g)\equal 
(a*\xi)(g)\equal \int_{h\in\Gamma_{x}}a(g,h)\xi(h)\,\mu_{x}(dh).$$
Then $a\mapsto \pi_x(a)$ defines a $*$-representation of $\aA_F$
in ${\mathcal{L}}(L^2(\Gamma_x;E))$.
Assume now that we are given a directed system (ordered by inclusion)
of closed subsets $F_i\subset \Gamma$ such that $s_{|F_i}$ is proper
for all $i$. Let $\aA\equal \lim_i\aA_{F_i}$,  and let $A$ be the completion of
 $\aA$ under the norm
$$\|a\|\equal \sup_{x\in M}\|\pi_x(a)\|.$$
Denote by $\tilde A$ and  $\tilde \aA$ the unitization of $A$
and $\aA$, respectively.

\begin{lem}\label{lem:gen holom}
${\aA}$ is a subalgebra of $A$, and  is stable under holomorphic
functional calculus.
\end{lem}
\begin{pf}
Let $\tilde b\in \tilde{\aA}$ be invertible in $A$. We  need to show that
$\tilde b$ is invertible in $\tilde{\aA}$.
Since $\tilde{\aA}$ is dense in $\tilde{A}$, there exists $x\in
\tilde{\aA}$ such that $\|1-\tilde{b} x\|<1/3$. Since
${\tilde b}^{-1}\equal x(\tilde{b}x)^{-1}$, we may assume that
$\|1-\tilde b\|<1/3$. Let $\tilde a \equal  1-\tilde b$.
We have $\tilde{a}\equal \lambda+b$, where $\lambda\in\cc$,
$|\lambda|<1/3$ and $b\in\aA$. Thus
 $\|b\|\equal  \|\tilde a-\lambda\| < 2/3$.
Let $a\equal (1-\lambda)^{-1}b$.
Since $(1-\tilde{a})^{-1}\equal (1-\lambda)^{-1}(1-a)^{-1}$,
it suffices  to prove that $(1-a)^{-1}\in\tilde{\aA}$ whenever $a\in\aA$
and $\|a\|<1$.

Let $a_n\equal a*a*\cdots*a$ ($n$ times). We show that the sum
$\sum_{n\equal 1}^\infty a_n$, and as well as all its derivatives,
 converges uniformly on every compact set.

Since $a_{n}(g,h)\equal [a_{n-1}*a(\cdot,h)](g)$, we have
$$\|a_n(\cdot,h)\|_{L^2(\Gamma_{s(h)})}\equal \|a_{n-1}*a(\cdot,h)\|
\le \|a\|^{n-1}\|a(\cdot,h)\|_{L^2(\Gamma_{s(h)})}$$
and similarly,
$$\|a_m(g,\cdot)\|_{L^2(\Gamma_{s(g)})}\le \|a\|^{m-1}
\|a(g,\cdot)\|_{L^2(\Gamma_{s(g)})}.$$
From the Cauchy-Schwarz inequality,
\begin{eqnarray*}
|a_{m+n}(g,h)|&\le& \|a_m(g,\cdot)\|_{L^2(\Gamma_{s(g)})}
\|a_n(\cdot,h)\|_{L^2(\Gamma_{s(h)})}\\
&\le& \|a\|^{m+n-2} \|a(g,\cdot)\|_{L^2(\Gamma_{s(g)})}
\|a(\cdot,h)\|_{L^2(\Gamma_{s(h)})}.
\end{eqnarray*}

It follows that $\sum_n a_n$ converges uniformly on every compact subset
of $\Gamma\times_s\Gamma$.  Similarly, one
shows that all derivatives converge
uniformly on any compact subset.
\end{pf}

\begin{prop}\label{prop:holom}
Let $S^1\to R\to\Gamma\toto M$ be an $S^1$-central extension of Lie groupoids.
Assume  that $\Gamma$ is proper. Then the subalgebra
$C_c^\infty(\Gamma;R)$ of $C^*_r(\Gamma;R)$ is stable under holomorphic
functional calculus.
\end{prop}

\begin{pf}
We use the construction above, where
$s\colon \Gamma\to M$ is the source map and the fiber bundle $E$
is $L\equal R\times_{S^1}\cc$. Let   $K\subset M$ be a compact subset, and 
$F_K\equal \Gamma^K$.  By the  properness of $\Gamma$, $s_{|F_K}$ is a
proper map. As above, define $\aA\equal \lim_{K} \aA_{\Gamma^K}$.
Denote by $\aA^\Gamma$
the subspace of $\aA$ consisting of $\Gamma$-invariant elements,
i.e. elements satisfying $a(g\gamma,h\gamma)\equal a(g,h)$, where
$L_g\otimes L_h^*$ and $L_{g\gamma}\otimes L_{h\gamma}^*$ are both
identified with  $L_{gh^{-1}}$.
Consider the map
\begin{eqnarray}\label{eqn:fa}
C_c(\Gamma;R)&\to& \aA^\Gamma\\
\nonumber
f&\mapsto&a,
\end{eqnarray}
given by  $a(g,h)\equal f(gh^{-1})\in L_{gh^{-1}}\cong L_g\otimes L_h^*$.
This  map is well-defined. Indeed,
if $f$ is compactly supported,  then there exists a compact subset $K$
of $\Gamma$ such that $\mbox{supp}(f)\subset \Gamma_K^K$. Therefore
it follows that $a\in{\aA}_{\Gamma^K}\cap\aA^\Gamma$.
Conversely, if $F\equal \Gamma^K$ and $a\in \aA_F\cap \aA^\Gamma$,
then $f(g)\equal a(g,s(g))$ is supported on $\Gamma_K^K$,
 which is compact
by the properness assumption, and $a(g,h)\equal f(gh^{-1})$ since $a$
is $\gm$-invariant. Therefore, the map defined by Eq. (\ref{eqn:fa})
 is bijective.
It is not hard to check that it is an isometric $*$-isomorphism
which extends to an isomorphism $C^*_r(\Gamma;R)\stackrel{\sim}{\to}A$.
The conclusion thus  follows from Lemma~\ref{lem:gen holom}.
\end{pf}

\begin{rmk}
Proposition~\ref{prop:holom} was proved in \cite[Lemma~7.5]{bc88}
in the non-twisted case for the crossed-product of a discrete group
acting properly on a manifold.
\end{rmk}

As an immediate consequence, we have the following

\begin{cor}\label{coro:iso in kth}
The inclusion $i\colon C_c^\infty(\Gamma;R)\to C^*_r(\Gamma;R)$
induces an isomorphism of $K$-theory.
\end{cor}


We now return  to the diagram (\ref{eqn:square iso}),
and  show that $V$ is well-defined.  We first need two preliminary
lemmas.

\begin{lem}\label{lem:smooth projection}
Let $P\in C^*_r(\Gamma;R)\otimes \kK$ be a projection and $\varepsilon>0$.
Then there exists a projection $P'\in C_c^\infty(\Gamma;R)\otimes \kK_0$,
where $\kK_0$ denotes the algebra of finite rank operators 
on $\HH$ such that $\|P'-P\|<\varepsilon$.
\end{lem}

\begin{pf}
Let $a\in C_c^\infty(\Gamma;R)\otimes\kK_0$ such that
$\|a-P\|<\varepsilon/2$. Then the spectrum of $a$ is contained in
the open set $U\equal B(0,\varepsilon/2)\cup B(1,\varepsilon/2)\subset\cc$. Let 
$f\colon U\to \cc$ be the function which is  equal to 0 on $B(0,\varepsilon/2)$
and  is equal to 1 on $B(1,\varepsilon/2)$.
Then $P':\equal f(a)$ is a projection such that
$\|P'-a\|<\varepsilon/2$, and $P'\in C_c^\infty(\Gamma;R)\otimes\kK_0$
by Proposition~\ref{prop:holom}.
\end{pf}

\begin{lem}\label{lem:smooth structure}
Suppose that $M$ is a manifold and $\pi\colon E\to M$ a Hermitian vector
bundle in the topological sense. Assume  that we are given a subspace
$\sS\subset C(M,E)$ such that
\begin{itemize}
\item[(a)] $\sS$ is a $C^\infty(M)$-module;
\item[(b)] for all $\xi$, $\eta\in \sS$, $x\mapsto \langle
\xi(x),\eta(x)\rangle$ is a smooth function on $M$;
\item[(c)] $\{\xi(x)|\; \xi\in \sS\}$ is dense in $E_x$ for all $x$.
\end{itemize}
Then there exists a unique smooth structure on the vector bundle $E$
such that $\sS$ consists of smooth sections.
\end{lem}

\begin{pf}
By the Gram-Schmidt orthonormalization process, there exists
an open cover $(U_i)$ of $M$ and sections $\xi_{i,1},
\ldots,\xi_{i,n}\in \sS$ such that for all $x\in U_i$,
$(\xi_{i,1}(x),\ldots,\xi_{i,n}(x))$ is an orthonormal basis
of $E_x$. Thus, we get local trivializations $\varphi_i
\colon \pi^{-1}(U_i)\cong U_i\times \cc^n$. Since $\langle\xi_{i,k},
\xi_{j,l}\rangle$ is smooth for all $i,j,k,l$, the change of
coordinates $\varphi_j\smalcirc\varphi_i^{-1}\colon
(U_i\cap U_j)\times \cc^n \to (U_i\cap U_j)\times \cc^n$
is smooth, thus we get a smooth structure on $E$. From
(b), it is clear that all elements of $\sS$ are smooth sections.
\par\medskip
Conversely, it is clear that if $E$ has a second smooth structure
such that all elements of $\sS$ are smooth sections, then $\varphi_i$
must be smooth for all $i$, and therefore the two smooth structures
coincide. 
\end{pf}

Now we return to the proof of Theorem \ref{thm:vb smooth}.

\begin{pff}
By assumption (see proof of Theorem~\ref{thm:vector bundles}),
there exists  an approximate unit $(P_n)$ in
$C^*_r(\Gamma;R)\otimes\kK$ consisting of projections.
According to  Lemma~\ref{lem:smooth projection},  there is  a projection $P'_n
\in C_c^\infty(\Gamma;R)\otimes \kK_0$ such that $\|P'_n-P_n\|<1/n$.
It is clear that
$(P'_n)$ is an approximate unit of $C_c^\infty(\Gamma;R)\otimes\kK_0$
consisting of projections. Hence according to  \cite[Proposition~5.5.5]{bla96},
 $K_0(C_c^\infty(\Gamma;R)
\otimes\kK_0)$ is the Grothendieck group of projections in
$C_c^\infty(\Gamma;R)\otimes\kK_0$.

Assume now that  $P\in C_c^\infty(\Gamma;R)\otimes\kK_0 \subset
{\mathcal{L}}(L^2(\Gamma;L)\otimes\HH)$ is  a projection.
Let $E\equal P(\tilde{L}^2(\Gamma;R)\otimes \HH )$.
 Then $E$ is a twisted vector bundle in the topological sense.
We say that a section of $E$ is smooth if it is of the form
$x\mapsto P_x\xi_x$, where $\xi\in C_c^\infty(\Gamma;L)\otimes \HH_n$
for some $n$. Since for any two smooth sections $\eta$ and $\zeta$,
$x\mapsto \langle\eta(x),\zeta(x)\rangle$ is smooth, the space
of smooth sections defines a smooth structure on $E$
according to  Lemma~\ref{lem:smooth structure}. It
follows that the map $V$ in (\ref{eqn:square iso}) is well-defined.
Also it is clear that the diagram (\ref{eqn:square iso}) is
commutative. 

Finally, we  prove that all maps in (\ref{eqn:square iso}) are
isomorphisms. For $i$ and $V'$, this follows from
Theorem~\ref{thm:vector bundles} and Corollary~\ref{coro:iso in kth}.
It  remains to show that $i'$ is injective.
\par\medskip

Assume  that $E$ and $F$ are smooth twisted vector bundles such that
$[E]-[F]\in \ker{i'}$. Then there exists a topological twisted
vector bundle $G$ such that $E\oplus G\cong F\oplus G$.
From the proof of Theorem~\ref{thm:vector bundles}, we know that there
exists a projection $P\in C^*_r(\Gamma;R)\otimes\kK$ such that
$G\cong P(\tilde{L}^2(\Gamma;R)\otimes\HH )$.  According to
Lemma~\ref{lem:smooth projection}, there exists a projection
$P'\in C_c^\infty(\Gamma;R)\otimes \kK_0$ such that
$\|P'-P\|<1$.  This implies that  $G\cong
P'(\tilde{L}^2(\Gamma;R)\otimes\HH )$. Therefore we can
 assume that $G$ is a smooth vector bundle. Replacing $E$ by $E\oplus G$ and $F$ by
$F\oplus G$, we see  that $E$ and $F$ are isomorphic as topological twisted
vector bundles. Let $T\equal (T_x)_{x\in M}$ be an isomorphism from $E$
to $F$. As in the proof of Proposition~\ref{prop:swan}, we can assume
 that $T_x$ is isometric for all $x$.
Let $T'\colon E\to F$ be
a (fiberwise linear, non equivariant) smooth morphism of vector
bundles such that $\|T'_x-T_x\|\le 1/2$ for all $x$. Choose
a smooth cutoff function $c\colon M\to \rr_+$ for the proper groupoid $R$.
Let
$$T''_x\equal \int_{R^x}\alpha_r(T'_{s(r)})c(s(r))\,\lambda^x(dr).$$
Since $T$ is $R$-invariant, $\|T''_x-T_x\|
\le \int_{R^x}\|\alpha_r(T'_{s(r)}-T_{s(r)})\|c(s(r))\,\lambda^x(dr)
\le 1/2\int_{R^x}c(s(r))\,\lambda^x(dr)\equal 1/2$. Therefore $T''_x$
is an isomorphism for all $x\in M$. Moreover, it is clear that
$x\mapsto T''_x$ is equivariant. It follows that $E$ and $F$
are isomorphic as smooth twisted vector bundles, and thus
$[E]-[F]\equal 0$. This completes the proof of the theorem.
\end{pff}

%

\section{The product $K_\alpha^i(\Gamma\upcom)\otimes
K_\beta^j(\Gamma\upcom)\to K_{\alpha+\beta}^{i+j}(\Gamma\upcom)$}
\subsection{The main idea} 
Let $S^1\to R\to \Gamma\toto M$ and $S^1\to R'\to \Gamma\toto M$
be $S^1$-central extensions, and let $\alpha$ and $\beta$ be their
corresponding classes in $H^2(\Gamma\upcom,\sS^1)$. 
It is simple to see that  there is a bilinear product
$$K_\alpha^{vb}(\Gamma\upcom)\otimes K_\beta^{vb}(\Gamma\upcom)
\to K_{\alpha+\beta}^{vb}(\Gamma\upcom)$$ 
defined as follows: let $E$ be a $(\Gamma,R)$-twisted vector bundle
and $E'$ a $(\Gamma,R')$-twisted vector bundle, then the product
of $[E]$ and $[E']$ is $[E\otimes E']$.

The objective of this section is to prove the following

\begin{them}
Let $\Gamma\toto M$ be a proper Lie groupoid such that $M/\Gamma$
is compact, and   $\alpha$, $\beta\in H^2(\Gamma\upcom,\sS^1)$.
Then there exists a bilinear, associative product
\begin{equation}
\label{eqn:KaKb}
K_\alpha^i(\Gamma\upcom)\otimes
K_\alpha^j(\Gamma\upcom)\to
K_{\alpha+\beta}^{i+j}(\Gamma\upcom)
\end{equation}
($i$, $j\in \{0,1\}$), which is compatible with the canonical map
$K_\gamma^{vb}(\Gamma\upcom)\to K_\gamma^0(\Gamma\upcom)$
($\gamma\in H^2(\Gamma\upcom,\sS^1)$).
\end{them}

In the theorem above, the  canonical map $\iota\colon
K_\gamma^{vb}(\Gamma\upcom) \to K_\gamma^0(\Gamma\upcom)$
is  constructed as in Section~\ref{subsec:ccase}.
Note that the construction of $\iota$ only requires the groupoid
$\Gamma$ to be proper, while the construction of the inverse of
$\iota$ as described in subsection~\ref{subsec:ccase}
requires all the hypotheses in
Theorem~\ref{thm:vector bundles}.
\par\bigskip
Recall that   in the Fredholm picture of twisted $K$-theory (see
Theorem~\ref{thm:fredholm proper}), the difficulty in
constructing the product (\ref{eqn:KaKb}) is to obtain a
Fredholm operator $T$ out of two Fredholm operators $T_1$
and $T_2$. Exactly the same difficulty appears in the
construction of the Kasparov product \cite{black98}. The existence
of the product $KK(A,D)\times KK(D,B)\to KK(A,B)$ (with
$A$, $D$ and $B$ separable $C^*$-algebras) is proved using
non constructive methods (in particular the Hahn-Banach theorem),
although explicit computations are possible in particular cases.

As a matter of fact, one can show, using Theorem~\ref{thm:fredholm proper},
that for a  proper groupoid   with compact orbit space
$\Gamma\toto M$, the \emph{untwisted} $K$-theory groups $K^i(\Gamma\upcom)$
are isomorphic to the equivariant $KK$-groups
$KK_{\Gamma}^i(C_0(M),C_0(M))$ defined by Le Gall~\cite{legall99}.
The existence of the product thus follows from the product in
$KK_\Gamma$-theory.\par\medskip

From the above discussion, it is quite natural to generalize
the $KK$-bifunctor further, and then try to identify it with
the twisted $K$-theory groups: this is the object of this section.
We will assume that the reader has some basic knowledge
about $KK$-theory \cite[chapter 8]{black98}. Since most of the
theory is already done in \cite{black98} or \cite{legall99},
we will only give those definitions and proofs that need substantial
modification.

\subsection{The $KK$ bifunctor}
Let us first recall a definition:

\begin{defn}
Let $A$ and $B$ be $C^*$-algebras. A $C^*$-correspondence from $A$
to $B$ is a pair $(\eE,f)$ where $\eE$ is a $B$-Hilbert module
and $f$ is a non-degenerate $*$-homomorphism from $A$ to $\lL(\eE)$.
\end{defn}

Recall that a $*$-homomorphism $\pi\colon A\to \lL(\eE)$ is said to be non-degenerate
if the closed linear span of $\pi(A)\eE$ is equal to $\eE$;
by Cohen's theorem, this is equivalent to $\pi(A)\eE\equal \eE$.
In particular, if $f\colon A\to B$ is a $*$-homomorphism, then 
$f$ induces a $C^*$-correspondence $(B,f)$.

Correspondences can be composed using the internal tensor product of
Hilbert modules: if $(\eE,f)$ is a $C^*$-correspondence from $A$ to
$B$ and $(\eE',g)$ is a $C^*$-correspondence from $B$ to $D$,
then $(\eE',g)\smalcirc (\eE,f)\equal (\eE\otimes_g\eE',f\otimes 1)$
is a $C^*$-correspondence from $A$ to $D$. Therefore, there is a
category $\cC$ whose objects are $C^*$-algebras and morphisms are
$C^*$-correspondences. And also there is a functor from the usual category
of $C^*$-algebras $\mathbf{C^*}$ to $\cC$
(given by the map $f\mapsto (B,f)$  as above).
 Moreover, isomorphism in the category $\cC$ is Morita equivalence.
\par\medskip

Recall that given a locally compact group $G$, Kasparov constructed
a bifunctor from the category of $G$-$C^*$-algebras to abelian
groups $(A,B)\mapsto KK_G(A,B)$ which is covariant in $B$ and
contravariant in $A$, and which is endowed with an associative
product $KK_G(A,B)\otimes KK_G(B,D)\to KK_G(A,D)$.
This construction was generalized by Le Gall to locally compact groupoids
admitting Haar systems \cite{legall99}.
Our goal in this subsection is to generalize
this construction further, by allowing the groupoid to act by
Morita equivalences on the $C^*$-algebras instead of by $*$-automorphisms,
i.e. to work in the category $\cC$ instead of $\mathbf{C^*}$.
This idea was communicated to us by Le Gall.
\par\medskip
For convenience, let us introduce some terminology.
\begin{defn}
Let $\Gamma\toto M$ be a locally compact groupoid.
Let $A$ be a $C^*$-algebra. A \emph{generalized action} of $\Gamma$
on $A$ is given by
\begin{itemize}
\item[(i)] a u.s.c. Fell bundle $\aA$ over $\Gamma$;
\item[(ii)] an isomorphism $A\cong C_0(M;\aA)$.
\end{itemize}
\end{defn}

For instance, if $\Gamma$ acts on $A$ in the usual sense, then
there exists a u.s.c. field of $C^*$-algebras $\aA'$ with fiber
$\aA'_x \cong A_x$ at $x\in M$ such that $A \cong C_0(M;\aA')$, and
the action of $\Gamma$ on $A$ induces $*$-isomorphisms
$\alpha_g\colon \aA'_{s(g)}\to \aA'_{t(g)}$. Let $\aA\equal s^*\aA'$,
with the product
\begin{eqnarray*}
\aA_g\otimes \aA_h\cong A_{s(g)}\otimes A_{s(h)}
&\to&  A_{s(h)}\cong \aA_{gh}\\
(a,b)&\mapsto& \alpha_{h^{-1}}(a)b
 \end{eqnarray*}
and the involution
\begin{eqnarray*}
\aA_g\cong A_{s(g)}&\to& A_{t(g)}\cong\aA_{g^{-1}}\\
a&\mapsto&\alpha_g(a^*).
\end{eqnarray*}
Then $\aA$ is a u.s.c. Fell bundle over $\Gamma$, and thus defines
a generalized action of $\Gamma$ on $A$.
\par\bigskip

If $A$ and $B$ are $C^*$-algebras endowed with a $\gm$-action,
there is a notion of equivariant $*$-homomorphism $f\colon A\to B$.
More generally, we want to introduce the definition of an
equivariant correspondence (Definition~\ref{def:equiv correspondence}).

We first introduce some notation: let $A$ be a $C^*$-algebra
endowed with a generalized action $\aA$  of a locally compact groupoid
$\Gamma$. Denote by
$\hat\aA$ the space of norm-bounded continuous maps
vanishing at infinity
$g\mapsto a'_g\in{\mathcal{A}}_{g^{-1}}$. $\hat\aA$ is naturally
a $t^*A$-Hilbert module with the module structure
$$(a'a)_g \equal  a'_g\cdot a_g \quad (a'\in\hat\aA,\; a\in t^*A)$$
and the scalar product
$$\langle a',a''\rangle_g\equal  {a'_g}^*a''_g\in A_{t(g)}.$$

\begin{defn}\label{def:equiv correspondence}
Let $A$ and $B$ be $C^*$-algebras endowed with generalized actions
$\aA$ and $\bB$ of a locally compact groupoid $\Gamma$. Let
$\eE$ be a $C^*$-correspondence from $A$ to $B$.
We say that the correspondence $\eE$ is equivariant if
there is an isomorphism of $s^* A$, $t^* B$ correspondences
$$W\colon s^*\eE\otimes_{s^*B} \hat\bB\to \hat\aA\otimes_{t^*A} t^*\eE$$
such that for every $(g,h)\in\Gamma^{(2)}$,
\begin{eqnarray*}
\lefteqn{({\mathrm{Id}}_{\aA_{h^{-1}}}\otimes W_g)
\smalcirc (W_h\otimes {\mathrm{Id}}_{\bB_{g^{-1}}})}\\
&\in& \lL(\eE_{s(h)}\otimes_{B_{s(h)}} \bB_{h^{-1}}
\otimes_{B_{t(h)}} \bB_{g^{-1}},
\aA_{h^{-1}}\otimes_{A_{t(h)}} \aA_{g^{-1}}
\otimes_{A_{t(g)}} \eE_{t(g)})
\end{eqnarray*}
is equal to
$$W_{gh}\in \lL(\eE_{s(h)}\otimes_{B_{s(h)}}\bB_{h^{-1}g^{-1}},
\aA_{h^{-1}g^{-1}}\otimes_{A_{t(g)}}\eE_{t(gh)})$$
via the identifications $\aA_{(gh)^{-1}}\cong
\aA_{h^{-1}}\otimes_{A_{t(h)}} \aA_{g^{-1}}$
and $\bB_{(gh)^{-1}}\cong
\bB_{h^{-1}}\otimes_{B_{t(h)}} \bB_{g^{-1}}$.
\end{defn}

When the action of $\Gamma$ on $A\simeq C_0(M,\tilde\aA)$
is an action in the usual sense, and $\eE$ is
an equivariant correspondence, then $\eE$ is a $(\Gamma,\bB)$-equivariant
$B$-Hilbert module (see Definition~\ref{def:GEmodule}).
Let $\lL(\tilde{\eE})\equal \coprod_{x\in M}\lL(\eE_x)$ be the bundle
defined in the appendix, preceding  Proposition~\ref{prop:LeE}.
Then the map $A\to \lL(\eE)$ induces a $\Gamma$-equivariant
bundle map $\tilde\aA\to \lL(\tilde\eE)$.

Note that there is a category $\cC_\Gamma$ whose objects consist
of $C^*$-algebras endowed with generalized actions of $\Gamma$,
and whose morphisms are equivariant correspondences.
\par\bigskip

To define the $KK$-groups, we first recall that if $F_2\in
\lL(\eE_2)$
then ${\mathrm{Id}}\otimes F_2$ does not
make sense. Instead, one has to use the notion of connection
\cite[Appendix A, pp. 1174-1178]{cs84}:

\begin{defn}
Let $\eE_1$ be a $D$-Hilbert module and $\eE_2$ a $D$, $B$-correspondence.
Let $\eE\equal \eE_1\otimes_{D}\eE_2$, $F_2\in \lL(\eE_2)$ and
$F\in \lL(\eE)$. We say that $F$ is a $F_2$-connection for $\eE_1$
if for every $\xi\in \eE_1$,
\begin{eqnarray*}
T_\xi F_2-(-1)^{\del\xi\del F_2} FT_\xi&\in& \kK(\eE_2,\eE)\\
F_2 T_\xi^*-(-1)^{\del\xi\del F_2} T_\xi^* F &\in&\kK(\eE,\eE_2).
\end{eqnarray*}
\end{defn}

The above operator $T_\xi\in\lL(\eE_2,\eE)$
is defined by
\begin{equation}\label{eqn:Txi}
T_\xi(\eta)\equal \xi\otimes \eta.
\end{equation}

Note the slight ambiguity, since $\eE_2$ does not appear in the notation
$T_\xi$.\par\medskip
Recall also that $F$ is a $F_2$ connection if and only if for all
$\xi\in \eE_1$, the graded commutator $\left[\theta_\xi,
\left(\begin{array}{cc}
F_2&0\\
0& F
\end{array}\right)\right]$
belongs to $\kK(\eE_2\oplus\eE)$, where
$\theta_\xi\equal \left(
\begin{array}{cc}
0&T_\xi^*\\
T_\xi&0
\end{array}
\right)$.
\par\bigskip

Let us now define the $KK$-groups. If $\eE$ is an equivariant
$A$, $B$-correspondence and $F\in\lL(\eE)$. Denote by
$t^* F\in \lL(t^*\eE)$ and $s^*F\in \lL(s^*\eE)$ the pull-backs
of $F$ by $t$ and $s$ respectively. Let
$$\sigma(F)\equal W(s^*F\otimes{\mathrm{Id}})W^*\in
\lL(\hat\aA\otimes_{t^*A}t^*\eE).$$

\begin{defn}\label{def:KK}
Let $A$ and $B$ be  $C^*$-algebras endowed with generalized actions
of $\Gamma$. An equivariant Kasparov $A$,
$B$-bimodule\footnote{``Kasparov correspondence'' might be a more
appropriate terminology but is not the usual one.}
is a pair $(\eE,F)$, where $\eE$ is a $\zz/2\zz$-graded,
equivariant $A$, $B$-correspondence and $F\in \lL(\eE)$ is a degree 1
operator such that for all $a\in A$,
\begin{itemize}
\item[(i)] $a(F-F^*)\in \kK(\eE)$;
\item[(ii)] $a(F^2-1)\in \kK(\eE)$;
\item[(iii)] $[a,F]\in \kK(\eE)$;
\item[(iv)] $\sigma(F)$ is a $t^*F$-connection for $\hat\aA$.
\end{itemize}
\end{defn}

If $\Gamma$ is a discrete group, (iv) holds if and only if
$Ad_g(F)-F$ is compact for all $g\in\Gamma$. Thus we will
refer to condition (iv) as the condition of invariance modulo
compacts.

As usual, unitarily equivalent Kasparov bimodules are identified.
Let ${\mathbf{E}}_\Gamma(A,B)$ the set of (unitary equivalence
classes of) Kasparov $A$, $B$-bimodules. A \emph{homotopy} in
$E_\Gamma(A,B)$ is given by an element of ${\mathbf{E}}_\Gamma(
A,B[0,1])$. The set of homotopy classes of elements of ${\mathbf{E}}_\Gamma
(A,B)$ is denoted by $KK_\Gamma(A,B)$. Then $KK_\Gamma(A,B)$ is an
abelian group, and $(A,B)\mapsto KK_\Gamma(A,B)$ is
a bifunctor, covariant in $B$, contravariant in $A$ (in
the category $\cC_\Gamma$).

\subsection{The technical theorem}
The main ingredient in the construction of the product
$$KK_\Gamma(A,D)\times KK_\Gamma(D,B)\to KK_\Gamma(A,B)$$
is the so-called technical theorem \cite[pp. 108-109]{hig87}.
We first need a lemma:

\begin{lem}\label{lem:lemtech}
Let $J$ and $J'$ be two $C^*$-algebras.
Let $\pi\colon J\to J'$ be a $*$-homomorphism,
$\varepsilon>0$, $h_0\in J_+$ such that $\|h_0\|<1$, $h\in J$,
$h'\in J'$, $\kK\subset Der(J)$ compact, $\kK'\subset
Der(J')$ compact. Then there exists $u\in J$ such that

\begin{tabular}{ll}
1) $h_0\le u$, $\|u\|<1$;\\
2) $\|uh-h\|\le\varepsilon$;&
2') $\|\pi(u)h'-h'\|\le \varepsilon$;\\
3) $\forall d\in \kK$, $\|[d,u]\|\le\varepsilon$;&
3') $\forall d'\in \kK'$, $\|[d',\pi(u)]\|\le \varepsilon$.
\end{tabular}
\end{lem}

The proof is almost the same as in \cite{hig87}.
Let us now come to the technical theorem.

\begin{them}\label{thm:technical theorem}
Let
\begin{itemize}
\item $A_1$ and $A'_1$ be two $C^*$-algebras such that
$A_1$ is $\sigma$-unital;
\item $J$ and $J'$ equivariant ideals in $A_1$ and $A'_1$
respectively;
\item $\pi\colon A_1\to A'_1$ such that $\pi(J)\subset J'$;
\item $\fF$ (resp. $\fF'$) a separable subspace of
$Der(A_1)$ (resp. of $Der(A'_1)$);
\item $a_2 \in M(J)_+$ such that $a_2 A_1\subset J$;
\item $a'_2\in M(J')_+$ such that $a'_2 A'_1\subset J'$.
\end{itemize}
Then there exists an element $M\in M(A_1)$, of degree 0, such that

\begin{tabular}{ll}
1) $M$ and $1-M$ are strictly positive;\\
2) $(1-M)a_2\in J$; &  2') $\pi(1-M)a'_2\in J'$;\\
3) $MA_1\subset J$; & 3') $\pi(M)A'_1\subset J'$;\\
4) $[\fF,M]\subset J$; & 4') $[\fF',\pi(M)]\subset J'$.
\end{tabular}
\end{them}

Again, the proof is almost the same as in
\cite{hig87}.

\subsection{The Kasparov product}
\begin{them}\label{thm:F1F2F}
Let $A$, $D$ and $B$ be separable $C^*$-algebras endowed with generalized
actions of a groupoid $\Gamma$.
Let $(\eE_1,F_1)\in {\mathbf E}_\Gamma(A,D)$ and
$(\eE_2,F_2)\in {\mathbf E}_\Gamma(D,B)$. Denote by $\eE$ the equivariant
$A$, $B$-correspondence $\eE\equal \eE_1\hat\otimes_D \eE_2$. Then the set
$F_1\hat\#_\Gamma F_2$ of operators $F\in \lL(\eE)$ such that
\begin{itemize}
\item $(\eE,F)\in {\mathbf{E}}_\Gamma(A,B)$;
\item $F$ is a $F_2$-connection for $\eE_1$;
\item $\forall a\in A$, $a[F_1\hat\otimes_D 1,F]a^*\ge 0$
modulo $\kK(\eE)$
\end{itemize}
is non-empty.
\end{them}

\begin{pf}
Choose a $F_2$-connection $T$ for $\eE_1$, and define
\begin{eqnarray*}
J&\equal & \kK(\eE)\\
A_1&\equal & \kK(\eE_1)\hat\otimes_D Id_{\eE_2} +J \subset \lL(\eE)\\
J'&\equal & \{S\in \lL(\hat\aA\otimes_{t^*A}t^*\eE\oplus t^*\eE)|\;
\forall \chi\in C_0(\Gamma),\; \chi S\in \kK\}\\
A'_1&\equal & \{S\in \lL(\hat\aA\otimes_{t^*A}t^*\eE\oplus t^*\eE)|\;\\
&&\qquad
\forall \chi\in C_0(\Gamma),\; \chi S\in \kK(\hat\aA\otimes_{t^*A}
t^*\eE_1\oplus t^*\eE_1)\otimes Id_{\eE_2}\}+J'\\
\fF&\equal &{\mathrm{Vect}}(Ad(F_1\hat\otimes Id_{\eE_2}),
Ad(T), Ad(a)\, (a\in A))\subset Der(A_1).
\end{eqnarray*}
Let $\pi\colon \lL(\eE)\to \lL(\hat\aA\otimes_{t^*A}t^*\eE
\oplus t^*\eE)$ defined by
$$\pi(S)\equal \left(
\begin{array}{cc}
\sigma(S)&0\\
0& t^*S
\end{array}
\right).$$
Then $\pi(A_1)\subset A'_1$ and $\pi(J)\subset J'$.
Let $a'_2$ be a strictly positive element of the $C^*$-algebra
generated by $[\theta_{a'},\pi(T)]$ ($a'\in \hat\aA$), where
$$\theta_{a'}\equal \left(
\begin{array}{cc}
0& T_{a'}\\
T_{a'}^* &0
\end{array}
\right)
\in \lL(\hat\aA\otimes_{t^* A} t^*\eE\oplus t^*\eE).$$
Let $\fF'\equal \{Ad_{\theta_a}|\; a'\in\hat\aA\}$.
Let $A_2$ be the sub-$C^*$-algebra of $\lL(\eE)$ generated
by
$$\{T-T^*; 1-T^2; [T,F_1\hat\otimes_D Id_{\eE_2}];
[T,a], \forall a\in A\}$$
and let $a_2$ be a strictly positive element of $A_2$.
We can apply the technical theorem~\ref{thm:technical theorem}
and thus obtain an operator
$M\in \lL(\eE)$ which satisfies the properties in
Theorem~\ref{thm:technical theorem}.
Let
$$F\equal M^{1/2}(F_1\hat\otimes_D Id_{\eE_2})+
(1-M)^{1/2}T.$$
As in the non-equivariant case, we have
$$(\eE,F)\in {\mathbf{E}}(A,B)$$
and thus it just remains to prove that $F$ satisfies the condition
of ``invariance modulo compacts'' (Definition~\ref{def:KK}(iv)).

Let $M_1\equal M^{1/2}(F_1\otimes Id)$ and $M_2\equal (1-M)^{1/2}T$. 
We show
that both $M_1$ and $M_2$ satisfy the invariance condition.
Since $F\equal M_1+M_2$, it  follows that $F$ also satisfies
the invariance condition.
\par\bigskip

We have $[\theta_{a'},\pi(M_1)]
\equal [\theta_{a'},\pi(M)^{1/2}]\pi(F_1\otimes 1)
+\pi(M)^{1/2}[\theta_{a'},\pi(F_1\otimes 1)]$.
From  property (4') of $M$, we have $[\theta_{a'},\pi(M)^{1/2}]
\in J'$. Since $J'$ is an ideal,
$[\theta_{a'},\pi(M)^{1/2}]\pi(F_1\otimes 1)\in J'$.
From property (iv)
in Definition~\ref{def:KK} for $F_1$, we have
$[\theta_{a'},\pi(F_1\otimes 1)]\in \kK(\hat\aA\otimes_{t^*A} t^*\eE_1
\oplus t^*\eE_1)\otimes Id_{\eE_2}\subset A'_1$.
Since $\pi(M)A'_1\subset J'$,
we have $\pi(M)^{1/2}[\theta_{a'},\pi(F_1\otimes 1)]\in J'$.
Finally, $[\theta_{a'},\pi(M_1)]\in J'$, which means
that $\sigma(M_1)$ is a $t^* M_1$-connection.
\par\bigskip

Let us show that $M_2$ satisfies the invariance condition
(Definition~\ref{def:KK} (iv)).
We have
$$[\theta_{a'},\pi((1-M)^{1/2}T)]
\equal [\theta_{a'},\pi(1-M)^{1/2}]\pi(T)+\pi(1-M)^{1/2}[\theta_{a'},
\pi(T)].$$
By the property (4') of $M$, $[\theta_{a'},\pi(1-M)]\equal -[\theta_{a'},\pi(M)]
\in J'$, and thus $[\theta_{a'},\pi(1-M)^{1/2}]\in J'$. Since $J'$ is an ideal,
we obtain $[\theta_{a'},\pi(1-M)^{1/2}]\pi(T)\in J'$.
\par\medskip

Since $\pi(1-M)a'_2\in J'$, we have $\pi(1-M)^{1/2}a'_2\in J'$ and
thus $\pi(1-M)^{1/2}[\theta_{a'},\pi(T)]\in J'$.
Finally, we get $[\theta_{a'},\pi(M_2)]\in J'$. This is
equivalent to the fact that $\sigma(M_2)$ is a $t^*M_2$-connection.
\par\bigskip
It now remains to show that the condition $a'_2A'_1\subset J'$
is fulfilled. It suffices to show that for all $a'\in\hat\aA$
and all $T'\in \kK(\hat\aA_{t^*A}\otimes t^*\eE_1\oplus t^*\eE_1)
\otimes Id_{t^*\eE_2}$, the operator $[\theta_{a'},\pi(T)]T'$
is compact.

Since $\kK(\hat\aA\otimes_{t^*A} t^*\eE_1\oplus t^*\eE_1)$ is the closed
vector subspace generated by operators of the form
$T_\zeta\smalcirc T_{\zeta'}^*$ ($\zeta,\zeta'\in\hat\aA\otimes t^*\eE_1
\oplus t^*\eE_1$), it suffices to show that
$[\theta_{a'},\pi(T)] T_\zeta \in\kK$ for all $\zeta$, i.e. that
\begin{itemize}
\item[(a)] $(T_{a'}t^* T -(-1)^{\del a'}\sigma(T) T_{a'})
T_{\xi'_1}\in\kK$, $\forall \xi'_1\in t^*\eE_1$;
\item[(b)] $(T_{a'}^*\sigma(T)-(-1)^{\del a'} (t^*T)T_{a'}^*)
T_{a''\otimes \xi'_1}\in \kK$, $\forall a''\otimes \xi'_1
\in \hat\aA\otimes_{t^*A} t^*\eE_1$.
\end{itemize}

Let us show (a).
Since $T$ is a $F_2$-connection, $t^*T$ is a $t^* F_2$-connection.
Hence $(t^*T) T_{\xi'_1}-(-1)^{\del \xi'_1}T_{\xi'_1}t^*F_2\in\kK$,
which in turn implies that
$$T_{a'}(t^*T)T_{\xi'_1}
-(-1)^{\del\xi'_1} T_{a'\otimes \xi'_1} (t^*F_2)\in\kK.$$
Therefore, it suffices to show that
\begin{equation}\label{eqn:commutator in K}
T_{a'\otimes \xi'_1}(t^*F_2)-(-1)^{\del(a'\otimes \xi'_1)}
\sigma(T)T_{a'\otimes \xi'_1}\in \kK.
\end{equation}
Let $W_1\colon s^*\eE_1\otimes_{s^*D}\hat\dD
\stackrel{\sim}{\to}\hat\aA\otimes_{t^*A}t^*\eE_1$ be the isomorphism
induced from the generalized action of $\Gamma$ on $\eE_1$.
Similarly, we introduce the obvious notations $W_2$ and $W$.

To show Eq. (\ref{eqn:commutator in K}), it suffices to prove
that for all $\xi''\otimes d'\in s^*\eE_1\otimes_{s^*D}\hat\dD$,
$$T_{W_1(\xi''_1\otimes d')}(t^*F_2)
-(-1)^{\del(\xi''_1\otimes d')}\sigma(T) T_{W_1(\xi''_1\otimes d')}
\in \kK.$$

Now $T_{W_1(\xi''_1\otimes d')}\equal (W_1\otimes Id)(T_{\xi''_1\otimes d'})$.
Moreover, from the invariance condition on $F_2$
(Definition~\ref{def:KK} (iv)), we get
$T_{d'}\smalcirc t^*F_2-(-1)^{\del d'}\sigma(F_2)\smalcirc T_{d'}\in \kK$.
Thus it suffices  to show that
$$\left(
(W_1\otimes Id)\smalcirc T_{\xi''_1}\smalcirc \sigma(F_2)
-(-1)^{\del\xi''_1}\sigma(T)\smalcirc (W_1\otimes Id)\smalcirc T_{\xi''_1}
\right)\smalcirc T_{d'}\in\kK.$$
Now, recall that
\begin{eqnarray*}
\sigma(T)&\equal & W(s^*T\otimes Id_{\hat\bB})W^{-1}\\
\sigma(F_2)&\equal &W_2 (s^*F_2\otimes Id_{\hat\bB})W_2^{-1}
\end{eqnarray*}
with $W\equal (W_1\otimes Id_{t^*\eE_2})\smalcirc (Id_{s^*\eE_1}\otimes W_2)$.
Hence, we are reduced to
\begin{eqnarray*}
\lefteqn{(W_1\otimes Id)T_{\xi''_1}W_2(s^* F_2\otimes Id_{\hat\bB}) W_2^{-1}}\\
&&-(-1)^{\del\xi''_1} W(s^*T\otimes Id_{\hat\bB}) W^{-1}(W_1\otimes Id)
T_{\xi''_1}\stackrel{?}{\in} \kK
\end{eqnarray*}
and then (multiplying on the left by $(W_1\otimes Id)^{-1}$) to
\begin{eqnarray*}
\lefteqn{T_{\xi''_1}\smalcirc W'_2 (s^* F_2\otimes Id_{\hat\bB})W_2^{-1}}\\
&&-(-1)^{\del \xi''_1} (Id_{s^*\eE_1}\otimes W_2)
(s^*T\otimes Id_{\hat\bB})(Id_{s^*\eE_1}\otimes W_2)^{-1}
T_{\xi''_1}\stackrel{?}{\in}\kK.
\end{eqnarray*}
Now (with the abuse of notation (\ref{eqn:Txi})),
\begin{eqnarray*}
T_{\xi''_1}\smalcirc W_2&\equal &W_2\otimes T_{\xi''_1}\\
(Id_{s^*\eE_1}\otimes W_2)\smalcirc T_{\xi''_1}&\equal &
T_{\xi''_1}\smalcirc (Id_{s^*\eE_1}\otimes W_2).
\end{eqnarray*}
We are finally reduced to showing that
$$T_{\xi''_1}(s^*F_2\otimes Id_{\hat\bB})
-(-1)^{\del\xi''_1}(s^*T\otimes Id_{\hat\bB})T_{\xi''_1}\in \kK,$$
which is true since $T$ is a $F_2$-connection.

This completes the proof of Eq. (\ref{eqn:commutator in K}).
\par\bigskip
Let us now show (b).
Using the fact that $t^*T$ is a $t^*F_2$-connection and that
$T_{a'}^*T_{a''\otimes \xi'_1}\equal T_{\langle a',a''\rangle\xi'_1}$,
we get
$$T_{a'}^*T_{a''\otimes \xi'_1}t^*F_2
-(-1)^{\del a'+\del a''+\del \xi'_1}(t^*T) T_{a'}^* T_{a''\otimes \xi'_1}
\in\kK.$$
Thus we need to show that
$$T_{a'}^*\sigma(T)T_{a''\otimes \xi'_1}
-(-1)^{\del(a'\otimes \xi'_1)} T_{a'}^*T_{a''\otimes \xi'_1}
(t^*F_2)\in\kK.$$
But this is immediate from Eq. (\ref{eqn:commutator in K}).
Thus (b) is proved.
\end{pf}

Theorem~\ref{thm:F1F2F} enables us to construct the Kasparov
product
$$KK_\Gamma(A,D)\times KK_\Gamma(D,B)\to KK_\Gamma(A,B)$$
of $[(\eE_1,F_1)]\in KK_\Gamma(A,D)$ and
$[(\eE_2,F_2)]\in KK_\Gamma(D,B)$ by $[(\eE,F)]\in KK_\Gamma(A,B)$,
where $\eE\equal \eE_1\hat\otimes_D\eE_1$ and $F\in F_1\hat\#_\Gamma F_2$.
As in the case of $C^*$-algebras endowed with an action of $\Gamma$
in the usual sense, the product is well-defined, bilinear,
homotopy-invariant, associative, covariant with respect to $B$
and contravariant with respect to $A$.

More generally, there is an associative product (for $i,j\in\{0,1\}$)
\begin{eqnarray}\nonumber
\lefteqn{KK_{\Gamma}^i(A,B_1\hat\otimes_{C_0(M)}D)\times
KK_{\Gamma}^j(D\hat\otimes_{C_0(M)}A_1,B)}\\
\label{eqn:KK product}
&&\quad
\to KK_{\Gamma,i+j}(A\hat\otimes_{C_0(M)}A_1,B_1\hat\otimes_{C_0(M)}B_1).
\end{eqnarray}

\subsection{Twisted $K$-theory is a $KK$-group}

Assume that $S^1\to R\to\Gamma\toto M$ is an $S^1$-central extension
of Lie groupoids. Recall that the line bundle $L\equal R\times_{S^1}\cc$
can be considered as a Fell bundle over the groupoid $\Gamma$, and thus
the $C^*$-algebra $C_0(M)$ is endowed
with a generalized action of $\Gamma$. Denote by $A_R$  this
$C^*$-algebra. Our goal is to show

\begin{prop}\label{prop:Kalpha equal KK}
If $\Gamma\toto M$ is a proper Lie groupoid
and $M/\Gamma$ is compact, then
for $i\equal 0,1$,
$KK_{\Gamma}^i(C_0(M),A_R)$ is isomorphic to $K^i_\alpha(\Gamma\upcom)$,
where $\alpha\in H^2(\Gamma\upcom,\sS^1)$ denotes the class
of the extension $S^1\to R\to \Gamma\toto M$.
\end{prop}

We will show a more general proposition, which can be considered
as a generalization of the Green-Julg theorem.

\begin{prop}\label{prop:GJ}
Let $\Gamma\toto M$ be a proper locally compact groupoid with a Haar system
such that $M/\Gamma$ is compact. Let $E$ be a u.s.c. Fell bundle
over $\Gamma$ and $A\equal C_0(M;E)$. Then $KK_\Gamma(C_0(M),A)$ and
$K_0(C^*_r(\Gamma;E))$ are isomorphic.
\end{prop}

Note that  Proposition~\ref{prop:GJ} implies
Proposition~\ref{prop:Kalpha equal KK}: take $E\equal L$ if $i\equal 0$
and $E\equal L\otimes C_0(\rr)$ if $i\equal 1$.

\begin{pf}
Let us construct a map $\Phi\colon
KK_\Gamma(C_0(M),A)\to K_0(C^*_r(\Gamma;E))$.
Consider $(\eE,F)\in {\mathbf{E}}_\Gamma(C_0(M),A)$.
Recalling Theorem~\ref{thm:K theory Fredholm}, we have to
construct a generalized Fredholm operator $T\in \fF^0(\Gamma,E)$.
By the stabilization theorem (Proposition~\ref{prop:stabilization}),
we may assume that $\eE\equal L^2(\Gamma;E)\otimes\HH
\oplus L^2(\Gamma;E)\otimes\HH$ with the obvious $\zz_2$-grading.
Then, replacing $F$ by $\frac{1}{2}(F+F^*)$, and then
by $F^\Gamma$ (see notation (\ref{eqn:TGamma})), we may
assume that $F$ is self-adjoint and $\Gamma$-invariant.
Thus, $F$ can be represented as a matrix
$$F\equal \left(
\begin{array}{cc}
0&T^*\\
T&0
\end{array}
\right).$$
Put $\Phi([\eE,F])\equal [T]$. It is routine to check that
$\Phi$ is a well-defined group homomorphism. The only slightly
tricky point is to check that $T$ is invertible modulo
$\kK_\Gamma(L^2(\Gamma;E)\otimes\HH)$. To see this, note that Condition
(ii) in Definition~\ref{def:KK} implies that $TT^*-\mbox{Id}$
and $T^*T-\mbox{Id}$ belong to $\cC(L^2(\Gamma;E)\otimes\HH)$.
By the compactness assumption on $M/\Gamma$ and the fact that
$T$ is $\Gamma$-invariant, we find that
$TT^*-\mbox{Id}$ and $T^*T-\mbox{Id}$ are in
$\kK_\Gamma(L^2(\Gamma;E)\otimes\HH)$.

\par\medskip

We now  construct a map in the other direction
$\Psi\colon K_0(C^*_r(\Gamma;E)) \to KK_\Gamma(C_0(M,A))$.
Let $T\in \fF^0(\Gamma,E)$.
Let $B\equal \lL(L^2(\Gamma;E)\otimes \HH)^\Gamma$ and
$J\equal \kK_\Gamma(L^2(\Gamma;E)\otimes\HH)$. By definition, $T$
is an element in $B$ whose image $\tau$ in $B/J$ is invertible.
Write the polar decomposition $\tau\equal u(\tau^*\tau)^{-1/2}$ and
lift $u$ to an element $T'\in B$. One easily proves by a standard
argument that $T'$ is homotopic to $T$. Therefore,
replacing $T$ by $T'$, we may assume that $T$ is unitary modulo $J$.
That is, $T^*T-\mbox{Id}$ and $TT^*-\mbox{Id}$ belong to $\kK_\Gamma
(\lL^2(\Gamma;E)\otimes \HH)$.
Let
$F\equal \left(
\begin{array}{cc}
0& T^*\\
T& 0
\end{array}
\right),$
which acts on the $\zz_2$-graded Hilbert module
$\eE\equal L^2(\Gamma;E)\otimes\HH\oplus L^2(\Gamma;E)\otimes\HH$.
It is not hard to check that $(\eE,F)\in {\mathbf{E}}_\Gamma(C_0(M),A)$.
Define $\Psi([T])\equal [(\eE,F)]$.  One can verify that
$\Phi$ and $\Psi$ are inverse from each other.
\end{pf}

\subsection{The product $K^i_\alpha(\Gamma\upcom)
\otimes K^j_\beta(\Gamma\upcom)\to K^{i+j}_{\alpha+\beta}
(\Gamma\upcom)$}
Suppose that $S^1\to R_1\to \Gamma\toto M$ and
$S^1\to R_2\to \Gamma\toto M$ are $S^1$-central extensions
of a Lie groupoid $\Gamma$. Denote by $\alpha$ and $\beta$ their
classes in $H^2(\Gamma\upcom,\sS^1)$.
Using the general Kasparov product (\ref{eqn:KK product}), we get
a product
$$KK_{\Gamma}^i(C_0(M),A_R)\otimes KK_{\Gamma}^j(C_0(M),A_{R'})
\to KK_{\Gamma}^{i+j}(C_0(M),A_{R\otimes R'}).$$
If in addition $\Gamma$ is proper and $M/\Gamma$ is compact, then
by  Proposition~\ref{prop:Kalpha equal KK}, we  obtain
a product
\begin{equation}
\label{eqn:product AB}
K^i_\alpha(\Gamma\upcom)
\otimes K^j_\beta(\Gamma\upcom)\to K^{i+j}_{\alpha+\beta}
(\Gamma\upcom).
\end{equation}

From the general properties of the Kasparov product
\cite{black98}, the product defined by Eq.  (\ref{eqn:product AB})
 is associative and graded commutative,
where graded commutativity comes from commutativity of the
diagram

$$\xymatrix{
A_R\otimes_{C_0(M)} A_{R'}\ar[r]\ar[d]^{\mbox{flip}}&
A_{R\otimes R'}\ar@{=}[d]\\
A_{R'}\otimes_{C_0(M)}A_R\ar[r]&
A_{R\otimes R'}.
}$$

\appendix
\section{Fell bundles over groupoids}
In this appendix, we recall the definition  and some
basic properties of a Fell bundle
over a groupoid (Definition~\ref{def:fell over groupoid})
and its reduced $C^*$-algebra.

\subsection{Fields of $C^*$-algebras}
\begin{defn}\label{defi:field banach}
Let $X$ be a Hausdorff
topological space. A continuous (resp. upper semicontinuous)
field of Banach spaces $E$ over $X$
consists of  a family $(E_x)_{x\in X}$ of Banach spaces  together
with a topology on $\tilde{E}\equal \coprod_{x\in X}E_x$
such that
\begin{itemize}
\item[(i)] the topology on $E_x$ induced from  that on
$\tilde{E}$ is the norm-topology;
\item[(ii)] the projection $\pi\colon \tilde{E}\to X$ is continuous
and open;
\item[(iii)] the operations $(e,e')\in \tilde{E}\times_X\tilde{E}
\mapsto e+e'\in\tilde{E}$ and $(\lambda,e)\in\cc\times \tilde{E}\to
\lambda e\in \tilde E$ are continuous;
\item[(iv)] the norm $\tilde E\to \rr_+$ is continuous (resp. u.s.c.);
\item[(v)] if $\|e_i\|\to 0$ and $\pi(e_i)\to x$, then $e_i\to 0_x$;
\item[(vi)] for all $e\in E_x$ there exists a continuous section
$\xi$ such that $\xi(x)\equal e$.
\end{itemize}
\end{defn}

In  \cite{dup-gil83,fd90} only  continuous fields were studied  and  they
are called Banach bundles. We will also use that terminology.
In this paper we are mainly concerned with continuous
fields, but most constructions and results only require the field
to be u.s.c..  In particular, a 
 field of Banach spaces can be constructed in the following way \cite{fd90}:

\begin{prop}\label{prop:section field}
Let $X$ be a topological Hausdorff space. Assume that
  $(E_x)_{x\in X}$
is a family of Banach spaces
 and $\Gam$ is a $C(X)$-module    
of sections of $ \tilde{E}:\equal \coprod_{x\in X}E_x\to X$
such that 
\begin{itemize}
\item[(i)] for every $\xi\in\Gam$, the function
$x\mapsto \|\xi(x)\|$ is continuous (resp. u.s.c.);
\item[(ii)] for all $x\in X$, the set $\{\xi(x)|\;\xi\in\Gam\}$ is
dense in $E_x$.
\end{itemize}
Then there is a unique topology on $\tilde E$ making $\tilde{E}\to X$
into a continuous (resp. u.s.c.) field of Banach spaces such that 
elements of $\Gam$ are exactly  continuous sections.
\end{prop}

In the same way, one defines  fields of Banach
algebras and fields  of $C^*$-algebras. For instance, if $f\colon Y\to X$
is a continuous map between two locally compact spaces, then $C_0(Y)$
may be considered as an u.s.c. field of $C^*$-algebras over $X$
with the  fiber $C_0(f^{-1}(x))$ at $x\in X$. Moreover, the field is continuous
if and only if $f$ is an open map.
\par\medskip

We will use the following conventions.
Denote by $C(X;E)$, $C_0(X;E)$ and $C_c(X;E)$
the space of continuous sections, the space of   continuous sections
vanishing at infinity, the space of compactly supported continuous sections
of the bundle $\tilde{E}\to X$, respectively. We  also use the notations
$C(X,\tilde{E})$, $C_0(X,\tilde{E})$ and $C_c(X,\tilde{E})$.
\par\medskip

Let us explain how pull-backs of fields are constructed. Let
$\coprod_{x\in X} E_x\to X$ be an u.s.c. (resp. continuous) field of
Banach spaces over $X$, and let $f\colon Y\to X$ be a continuous map.
Then the u.s.c. (resp. continuous)
field $f^*E$ is the field with the fiber $E_{f(y)}$ at
$y\in Y$, and whose total space is $Y\times_X\tilde{E}$
with the induced topology from $Y\times \tilde{E}$.
If $E$ is determined by a $C(X)$-module of sections $\Gam
\subset C(X,{\tilde{E}})$ as in Proposition~\ref{prop:section field},
then $f^* E$ is determined by
$f^*\Gam\equal \{\xi\smalcirc f|\; \xi\in \Gam\}$.
\par\medskip

Recall that if $X$ is a locally compact space, then
a $C_0(X)$-algebra is  a $C^*$-algebra $A$ together with  a *-homomorphism
$C_0(X)\to Z(M(A))$ (the center of the multiplier algebra of $A$)
such that $C_0(X)A\equal A$. The proposition below indicates that there is
a bijection  between $C_0(X)$-algebras and u.s.c.
fields of $C^*$-algebras over $X$.

For any  $x\in X$,  by $C_x(X)$, we denote
 the ideal of $C_0(X)$ consisting of
functions that vanish at $x$.

\begin{prop}\label{prop:C(X)-algebra}
Let $X$ be a locally compact space, 
 $A$ a $C_0(X)$-algebra and $A_x\equal A/(C_x(X)A)$. Denote by
$\pi_x\colon A\to A_x$ the projection.  There
is a unique u.s.c. field of $C^*$-algebras $\tilde{A}:\equal 
\coprod_{x\in X}A_x\to X$
such that the map
\begin{eqnarray*}
A&\to& C_0(X,\tilde{A})\\
a&\mapsto&(x\mapsto \pi_x(a))
\end{eqnarray*}
is an isomorphism of $C^*$-algebras.

Conversely, assume that $\tilde{\aA}\equal \coprod_{x\in X}\aA_x\to X$ is a u.s.c.
field of $C^*$-algebras over $X$, and  $A\equal C_0(X,\tilde{\aA})$ is  the
space of continuous sections vanishing at infinity. Then $A$ is obviously
a $C_0(X)$-algebra, and the evaluation map $A\to \aA_x$ induces a
$*$-isomorphism $A_x\to \aA_x$.
\end{prop}

\begin{pf}
This is immediate from \cite[Proposition~2.12 a)]{bla96}.
\end{pf}

Assume that $\Gamma\toto X$ is a topological groupoid. Recall
\cite[Definition~3.3]{legall99} that
a $\gm$-action  on a $C_0  (X)$-algebra $A$ is an isomorphism
of $C_0(\Gamma)$-algebras $\alpha\colon s^*A\to t^*A$ such that
$\alpha_{gh}\equal \alpha_g\alpha_h$ for all $(g,h)\in\Gamma^{(2)}$,
where $\alpha_g\colon (s^*A)_g\cong A_{s(g)} \to (r^*A)_g\cong A_{t(g)}$
is the induced isomorphism.

Let $\tilde{\aA}\equal \coprod_{x\in X}\aA_x\to X$ be an u.s.c. field
of $C^*$-algebras. We say that the groupoid $\Gamma$ acts on
$\tilde{\aA}\to X$ if there is an isomorphism $\alpha\colon
s^*\aA\to t^*\aA$ of fields of $C^*$-algebras over $\Gamma$
such that $\alpha_{gh}\equal \alpha_g\alpha_h$ for all
$(g,h)\in\Gamma^{(2)}$. It is clear that using the dictionary above
(Proposition~\ref{prop:C(X)-algebra}), $\gm$-actions  on
$C_0(X)$-algebras are in bijective correspondence with
$\gm$-actions  on fields of $C^*$-algebras over $X$.

Now, let us explain how $C^*$-modules over a $C_0(X)$-algebra
can be considered as u.s.c. fields of Banach spaces over $X$.

\begin{prop}\label{prop: module field}
Let $\tilde{{\mathcal{A}}}\equal \coprod_{x\in X}
\aA_x\to X$ be an u.s.c. field of $C^*$-algebras
over $X$ and  $A\equal C_0(X,\tilde{\aA})$.
Assume that   $\eE$ is  an $A$-Hilbert
module. Let $\eE_x:\equal \eE\otimes_{A}A_x$ and denote by
$\pi_x\colon \eE\to \eE_x$ the canonical map. Then there is an unique
u.s.c. field of Banach spaces $\tilde{\eE}:\equal \coprod_{x\in X}\eE_x
\to X$ such that the map
\begin{eqnarray*}
\eE&\to&C_0(X,\tilde{\eE})\\
\xi&\mapsto& (x\mapsto \pi_x(\xi))
\end{eqnarray*}
is an isomorphism.
Moreover, if the field $\tilde{\aA}\to X$ is a continuous field,
then $\tilde{\eE} \to X$ is a continuous field as well.
\end{prop}

The proof, which uses Proposition~\ref{prop:section field},
 is straightforward and is  left to the reader. In particular, any
$C_0(X)$-module is the space of continuous sections vanishing
at infinity of a continuous field of Hilbert spaces.

Consider an u.s.c. field of $C^*$-algebras
$\tilde{\aA}\equal \coprod_{x\in X} \aA_x\to X$. Let
$A\equal C_0(X,\tilde{\aA})$. Assume  that $\eE$ is an $A$-Hilbert module.
It is simple to show that there is a unique topology on
${\mathcal{L}}(\tilde{\eE})
:\equal \coprod_{x\in X}{\mathcal{L}}(\eE_x)$ such that for every
net $T_i\in {\mathcal{L}}(\eE_{x_i})$ and $T\in
{\mathcal{L}}(\eE_x)$, $T_i$ converges to $T$ if and only if
for every $\xi\in C(X,\tilde{\eE})$,
\begin{itemize}
\item[(i)] $x_i\to x$;
\item[(ii)] $T_i\xi(x_i)\to T\xi(x)$; and
\item[(iii)] $T_i^*\xi(x_i)\to T^*\xi(x)$.
\end{itemize}
Then the bundle ${\mathcal{L}}(\tilde{\eE})\to X$ satisfies
all the properties of Definition~\ref{defi:field banach},
except that the norm is not necessarily u.s.c.
(in fact, one can show that it is
\emph{lower} semi-continuous if $\tilde{\eE}\to X$ is
a continuous field), and the induced topology on ${\mathcal{L}}(\eE_x)$
is not the norm-topology.

We say that a section
$x\mapsto T_x$ of ${\mathcal{L}}(\tilde{\eE})$ is strongly continuous
if for every $\xi\in C(X,\tilde{\aA})$,
$x\mapsto T_x\xi(x)$ belongs to
$C(X,\tilde{\aA})$, and a section $x\mapsto T_x$ is $*$-strongly continuous
if both  $x\mapsto T_x$ and $x\mapsto T_x^*$ are strongly continuous.
It is not hard to show that a section is $*$-strongly continuous
if and only if it is a continuous section of the bundle defined above.
Denote by $C_b(X,{\mathcal{L}}(\tilde{\eE}))$
the space of continuous and norm-bounded sections.

\begin{prop}\label{prop:LeE}
There is an isomorphism
\begin{eqnarray*}
{\mathcal{L}}(\eE)&\to& C_b(X,{\mathcal{L}}(\tilde{\eE}))\\
T&\mapsto&(x\mapsto T_x),
\end{eqnarray*}
where $T_x\equal T\otimes_A {\mathrm{Id}}\in {\mathcal{L}}(\eE
\otimes_{A}A_x)\equal {\mathcal{L}}(\eE_x)$.
\end{prop}

\begin{pf}
This follows directly  from Proposition~\ref{prop: module field} and the fact
that $\lL(\eE)$ is, by definition,  the space of maps from $\eE$
to $\eE$ admitting an adjoint.
\end{pf}

The analogue of the above proposition for ${\mathcal{K}}(\eE)$ is
less simple. However since we do not need it in full generality in this paper,
we only consider a particular case below.

\begin{prop}\label{prop:KeE}
Let $\tilde{\hH}\equal \coprod_{x\in X}\hH_x\to X$ be a continuous field of
Hilbert spaces, and  $\hH\equal C_0(X,\tilde{\hH})$ be 
 the associated $C_0(X)$-Hilbert
module. Then there exists a unique topology on
$\kK(\tilde{\hH}):\equal \coprod_{x\in X}\kK(\hH_x)$ such that
\begin{itemize}
\item[(i)] the field $\kK(\tilde{\hH})\to X$ is a continuous
field of $C^*$-algebras;
\item[(ii)] for every
$\xi$, $\eta\in C_0(X,\tilde{\hH})$,  we have $(x\mapsto
T_{\xi(x),\eta(x)})\in C_0(X,\kK(\tilde{\hH}))$.
\end{itemize}
Moreover, the map
$$T_{\xi,\eta}\mapsto (x\mapsto T_{\xi(x),\eta(x)})$$
extends uniquely to an isomorphism of $C^*$-algebras
$\kK(\hH)\stackrel{\sim}{\to}C_0(X,\kK(\tilde{\hH}))$.
\end{prop}

\begin{pf}
We sketch the proof in the case that  the field is countably generated.
In this case, by the stabilization theorem (either by
Dixmier and Douady \cite{dix-dou63} or Kasparov \cite{wegge93}),
we may assume that the field is trivial: $\tilde{\hH}\cong
X\times \HH$. It is well-known that $\kK(\hH)$
is isomorphic to $C_0(X,\kK (\HH) )$,  where $\kK (\HH)  $ is endowed with
the norm-topology. Thus it  follows from
Proposition~\ref{prop:section field} that $C_0(X, \kK  (\HH))$ is
the space of continuous sections vanishing at infinity of a
continuous field of $C^*$-algebras over $X$ with fibers isomorphic
to $\kK  (\HH)$.
\end{pf}

\subsection{Fell bundles over groupoids: definition and first properties}

\begin{defn}\label{def:fell over groupoid}
Let $\Gamma\toto M$ be a locally compact groupoid and denote by $m\colon
\Gamma^{(2)}\to \Gamma$ the multiplication map. A continuous
(resp. u.s.c.) Fell bundle over $\Gamma$
is a continuous (resp. u.s.c.)
field of Banach spaces $(E_g)_{g\in\Gamma}$ over $\Gamma$
together with an associative bilinear product
$(\xi,\eta)\in E_g\times E_h\mapsto \xi\eta\in E_{gh}$,  whenever
$(g,h)\in \Gamma^{(2)}$, and an antilinear involution $\xi\in E_g
\mapsto \xi^*\in E_{g^{-1}}$ such that for any  $(g,h)\in \Gamma^{(2)}$,
and $(e_1,e_2)\in E_g\times E_h$,
\begin{itemize}
\item[(i)] $\|e_1e_2\|\le \|e_1\| \|e_2\|$;
\item[(ii)] $(e_1e_2)^*\equal e_2^*e_1^*$;
\item[(iii)] $\|e_1^* e_1\|\equal \|e_1\|^2$;
\item[(iv)] $e_1^*e_1$ is a positive element of the $C^*$-algebra
$E_{s(g)}$;
\item[(v)] the product $(e,e')\in m^*(\tilde{E})\mapsto
ee'\in \tilde{E}$,  and the involution
$e\in\tilde{E}\mapsto e^*$ are continuous;
\item[(vi)] for all $(g,h)\in \gm^{(2)}$,
the image of the product $E_g\times E_h\to E_{gh}$ spans a dense subspace of
$E_{gh}$.
\end{itemize}
\end{defn}

\begin{rmk}
Note that (i)--(iii) imply  that $E_x$, $x\in M$,  is a $C^*$-algebra, so
(iv) makes sense.
\end{rmk}

Continuous Fell bundles were first defined by Yamagami in \cite{yam90}, and
were  called $C^*$-algebras over groupoids. Since continuous
Fell bundles are simply called ``Fell bundles'' in the literature
\cite{fm98,kum96,muh01}, we will follow this convention.
In the literature one also finds the terminology ``full'' Fell
bundle: this refers to Condition (vi).
Note that  $A: \equal C_0(M;E)$, the restriction of $C_0(\Gamma;E)$ to $M$,
 is a $C_0(M)$-algebra, and $A_x\equal E_x$ for all $x\in M$,
by Proposition~\ref{prop:C(X)-algebra}.

\begin{numex}
Let $\Gamma$ be a locally compact groupoid acting on a
$C_0(M)$-algebra $A$, ${\mathcal{A}}$  the
associated u.s.c. field of $C^*$-algebras
(Proposition~\ref{prop:C(X)-algebra}). There is an isomorphism
$\alpha\colon s^*\tilde{\aA}\to t^*\tilde{\aA}$
such that $\alpha_{gh}\equal \alpha_g\smalcirc\alpha_h$
for all $(g,h)\in \Gamma^{(2)}$. Then $E\equal s^*\aA$ is a u.s.c. Fell bundle
over $\Gamma$ with the product $(a,b)\in E_g\times E_h\equal \aA_{s(g)}\times
\aA_{s(h)}\mapsto \alpha_{h^{-1}}(a)b\in E_{gh}\equal \aA_{sh}$ and
the involution $a\in E_g\mapsto \alpha_g(a^*)\in E_{g^{-1}}$.
\end{numex}

Therefore the notion of u.s.c. Fell bundles over $\Gamma$
 generalizes that of  actions of $\Gamma$ on
$C^*$-algebras. In fact, u.s.c. Fell bundles over $\Gamma$ can be viewed as
``actions of $\Gamma$ on $C^*$-algebras by Morita equivalences''
(see \cite{muh01}).
\par\medskip

Now we return to  the discussion on  a general u.s.c. Fell bundle $E$.
Define an $A_{s(g)}$-valued scalar product on $E_g$ by
$\langle e,e'\rangle \equal  e^*e'$. Then $E_g$ becomes an $A_{s(g)}$-Hilbert
module, and the left multiplication by elements of $A_{t(g)}$ defines
a *-homomorphism $A_{t(g)}\to {\mathcal{L}}(E_g)$.
In other words, $E_g$ is an $A_{t(g)}$-$A_{s(g)}$-correspondence.

Note also that the product $E_g\times E_h\to E_{gh}$
induces an isomorphism of $A_{t(g)}-{A_{s(h)}}$ bimodules
$E_g\otimes_{A_{s(g)}}E_h\to E_{gh}$. Indeed, to check that this map
is isometric, we note that $\forall \xi_i\in E_g$, $ \eta_i\in
E_h$,
\begin{eqnarray*}
\lefteqn{\langle\sum_i\xi_i\otimes\eta_i,\sum_i\xi_i\otimes \eta_i\rangle\equal 
\sum_{i,j}\langle\eta_i,
\langle\xi_i,\xi_j\rangle\eta_j\rangle}\\
&\equal &\sum_{i,j}(\xi_i\eta_i)^*(\xi_j\eta_j)\equal 
\langle\sum_i\xi_i\eta_i,\sum_i\xi_i\eta_i\rangle.
\end{eqnarray*}

The surjectivity of $E_g\times E_h\to E_{gh}$ follows from Condition (vi)
of Definition~\ref{def:fell over groupoid}.
\par\medskip

The following proposition justifies the reason that we require the field to be u.s.c.:

\begin{prop}\label{prop:approximation}
If $E$ is  an u.s.c. Fell bundle over the groupoid $\Gamma$, then
sections of the form $(g,h)\mapsto \sum_i\xi_i(g)\eta_i(h)$, where
$\xi_i$, $\eta_i\in C_0(\Gamma;E)$, are dense in $C_0(\Gamma^{(2)},m^*E)$.
\end{prop}

To prove the proposition, we need the following:

\begin{lem}\label{lem:topolem}
Let $K$ and $L$ be two compact spaces, $(\Omega_k)$ an open cover of $K
\times L$. Then there exist finite covers $U_i$ and $V_j$ of $K$ and
$L$ respectively such that $(U_i\times V_j)$ is a refinement of
$(\Omega_k)$.
\end{lem}

\begin{pf}
For every $(x,y)\in K\times L$, there exist
$K_{x,y}$, $L_{x,y}$ compact and $U^1_{x,y}$, $V^1_{x,y}$ open such that
$(x,y)\in {\mbox{Int}}(K_{x,y})\times{\mbox{Int}}(L_{x,y})
\subset {K}_{x,y}\times{L}_{x,y}\subset
{U}^1_{x,y}\times{V}^1_{x,y}\subset\Omega_{k}$ for some $k$. 
Let $U^2_{x,y}\equal K-K_{x,y}$
and $V^2_{x,y}\equal L-L_{x,y}$. By compactness, there exists a finite
family $(x_i,y_i)_{i\in I}$ such that $\cup_{i\in I}
{\mbox{Int}}(K_{x_i,y_i})\times
{\mbox{Int}}(L_{x_i,y_i})$ covers $K\times L$.

For any  $\alpha\equal (\alpha_i)_{i\in I}\in \{1,2\}^I$, let
$U^\alpha\equal \cap U^{\alpha_i}_{x_i,y_i}$ and
$V^\alpha\equal \cap V^{\alpha_i}_{x_i,y_i}$. It is not hard to check that
$(U_\alpha)$ and $(V_\alpha)$ are, respectively, covers of $K$ and $L$
 that satisfy the required properties.
\end{pf}

\begin{pf}
Let $\zeta\in C_0(\Gamma^{(2)},m^* E)$.
We can assume that $\zeta$ is compactly supported.
There exist $K$ and $L\subset \Gamma$ compact such that the support
of $\zeta$ is in the interior of
$KL\equal \{gh|\;(g,h)\in K*L\equal (K\times L)\cap \Gamma^{(2)}\}$.
By the definition of $m^*E$, for every $(g,h)\in K*L$ there exist
$\xi_{g,h}^i$, $\eta_{g,h}^i
\in C_0(\Gamma;E)$ such that
$$\|\sum_i\xi^i_{g,h}(g)\eta_{g,h}^i(h)-\zeta(g,h)\|<\varepsilon.$$
Since the field $m^*E$ is u.s.c., there exists a neighborhood
$\Omega_{g,h}$ of $(g,h)$ such that
$$\|\sum_i\xi^i_{g,h}(g')\eta_{g,h}^i(h')-\zeta(g',h')\|<\varepsilon$$
for all $(g',h')\in \Omega_{g,h}$.

Now, by Lemma~\ref{lem:topolem},
there exist  compactly supported nonnegative
continuous functions  $\varphi_k$ and $\psi_l$  on $\Gamma$
 such that $\sum_{k}\varphi_k\equal 1$
on $K$, $\sum_l\psi_l\equal 1$ on $L$,
$0\le \sum_k\varphi_k\le 1$,
$0\le \sum_l\psi_l\le 1$,
and $(g,h)\mapsto
\varphi_k(g)\psi_l(h)$ is supported
in some $\Omega_{g_k,h_l}$ or in $(K\times L)-K*L$ for all $k,l$.

Thus,
$$\|\zeta(g,h)-\sum_{k,l,i}\varphi_k(g)\psi_l(h)\xi^i_{g_k,h_l}(g)\eta^i_{g_k,h_l}(h)\|<\varepsilon,$$
for all $(g,h)\in K*L$.
Now, choose two compact sets $K'$ and $L'$ whose interior contain
$K$ and $L$ respectively. Applying the above to $K'$ and $L'$ instead
of $K$ and $L$, there exist $\xi_i$, $\eta_i\in E$ such that
\begin{eqnarray}\label{eqn:approx}
\|\zeta(g,h)-\sum_i \xi_i(g)\eta_i(h)\|&<&\varepsilon
\end{eqnarray}
for $(g,h)\in K'*L'$. Replacing
$\xi_i$ by $\varphi\xi$, where $\varphi\in C_c(\Gamma)_+$
has the 
support $\subset{\mbox{Int}}(K')$, $0\le\varphi\le 1$ and $\varphi\equal 1$
on $K$, and replacing
$\eta_i$ by $\psi\eta_i$, where $\psi\in C_c(\Gamma)_+$
has the support $\subset{\mbox{Int}}(L')$, $0\le\psi\le 1$ and $\psi\equal 1$
on $L$, we may assume that Eq.(\ref{eqn:approx})
holds for all $(g,h)\in \Gamma^{(2)}$.
\end{pf}

\subsection{The reduced $C^*$-algebra}\label{subsec:reduced}

In this subsection, we recall the definition of the reduced $C^*$-algebra
associated to an u.s.c. Fell bundle over a groupoid. See~\cite{muh01},
or~\cite[Section 7.7]{ped79} for the definition of the
 crossed-product algebra by a locally compact group, 
or~\cite[Chapter 2]{ren80} for the $C^*$-algebra of a groupoid.

Assume  that $\Gamma$ is a locally compact groupoid with a Haar system,
 and $E$ is an u.s.c. Fell bundle over $\Gamma$.
Let $C_c(\Gamma;E)$ denote the space of compactly supported continuous
sections.
For  $\xi$, $\eta\in C_c(\Gamma;E)$,
define the convolution by
$$(\xi*\eta )(g)\equal \int_{h\in \Gamma^{t(g)}}\xi(h)\eta(h^{-1}g)\,
\lambda^{t(g)}(dh)$$
and the involution by
$\xi^*(g)\equal \xi(g^{-1})^*$.

Let us check that $\xi*\eta$ belongs to $C_c(\Gamma;E)$. By (v) in
Definition~\ref{def:fell over groupoid}, $(g,h)\mapsto
\xi(h)\eta(h^{-1}g)$ is the uniform limit of maps of the form
$\sum_i f_i(h,h^{-1}g)\zeta_i(h(h^{-1}g))$,  where $f_i\in C_c(\Gamma^{(2)})$
and $\zeta_i\in C_c(\Gamma;E)$,
and hence of sums of the form $f(h)f'(g)\zeta(g)$,  where
$f$, $f'\in C_c(\Gamma)$ and $\zeta\in C_c(\Gamma;E)$. Moreover, the
function $f'(g)$ can be assumed to be supported on a fixed
compact subset of $\Gamma$. Now,
$$\int_{h\in\Gamma^{t(g)}}f(h)\,\lambda^{t(g)}(dh) f'(g)\zeta(g)$$
is the product of $\zeta$ by an element of $C_c(\Gamma)$, and hence
belongs to $C_c(\Gamma;E)$.
Therefore $\xi*\eta$ can be uniformly approximated by
elements in $C_c(\Gamma;E)$.

Let $\|\xi\|_1\equal \sup_{x\in M}\int_{\Gamma^x}\|\xi(g)\|\,
\lambda^x(dg)$,
and $\|\xi\|_I\equal \max(\|\xi\|_1,\|\xi^*\|_1)$. Then the completion of
$C_c(\Gamma;E)$ with respect to the norm $\|\cdot\|_I$
is a Banach $*$-algebra, and is denoted by $L^1(\Gamma;E)$.
Its enveloping $C^*$-algebra is
denoted by $C^*(\Gamma;E)$,  and is called the $C^*$-algebra
of the field $E$.

Let $L^2(\Gamma;E)$ be the $A$-Hilbert module obtained by
completing $C_c(\Gamma;E)$ with respect to the $A$-valued scalar product:
$$\langle\xi,\eta\rangle (x) \equal  \int_{g\in \Gamma_x}\langle\xi(g),\eta(g)\rangle
\,\lambda_x(dg)\in A_x.$$
Then for every $\xi\in C_c(\Gamma;E)$, the
map $\pi_l(\xi)\colon \eta\mapsto \xi*\eta$
 belongs to
$\mathcal{L}( L^2(\Gamma;E))$, and $\xi\mapsto \lambda(\xi)$ extends
to a representation of $L^1(\Gamma;E)$, called the left
regular representation.
Its image
$C^*_r(\Gamma;E)\equal \overline{\pi_l(L^1(\Gamma;E))}
\equal \overline{\pi_l(C_c(\Gamma;E))}
\subset \mathcal{L}( L^2(\Gamma;E))$
is called the reduced $C^*$-algebra of the field $E$.
\par\medskip

The $A$-Hilbert module $L^2(\Gamma;E)$ can be considered,
by Proposition~\ref{prop: module field}, as a field of
Banach spaces over $M$ with fiber $L^2(\Gamma;E)\otimes_A
A_x$ at $x\in M$. Denote the total space of this bundle
by $\tilde{L}^2(\Gamma;E)$ and the fiber
by $L^2(\Gamma_x;E)$. To justify our notation,
let $i_x\colon \Gamma_x\to\Gamma$ be the inclusion. Then
$L^2(\Gamma_x;E)$ is the completion of $C_c(\Gamma_x;i_x^*E)$
with respect to the $A_x$-valued scalar product
$\langle\xi,\eta\rangle \equal \int_{g\in\Gamma_x}\langle\xi(g),\eta(g)\rangle
\,\lambda_x(dg)$. Thus $L^2(\Gamma_x;E)$ is an $A_x$-Hilbert module.
\par\medskip

The algebra of compact operators ${\mathcal{K}}(L^2(\Gamma;E))$ is
a field of $C^*$-algebras over $M$ whose total space is denoted
by $\kK(\tilde{L}^2(\Gamma;E))$ (see Proposition~\ref{prop:KeE}).
Its fiber at $x\in M$ is ${\mathcal{K}}(L^2(\Gamma_x;E))$.
If $E$ is a continuous Fell bundle, then ${\mathcal{K}}(L^2(\Gamma;E))$ is
a continuous field of $C^*$-algebras \cite[pp. 76-77]{muh01}.

The $C^*$-algebra ${\mathcal{K}}(L^2(\Gamma;E))$ is
endowed with a continuous action of $\Gamma$:
for every $\gamma \in\Gamma_x^y$, the  map
$$\alpha_\gamma \colon {\mathcal{K}}(L^2(\Gamma_x;E)) 
\stackrel{\sim}{\to} {\mathcal{K}}(L^2(\Gamma_y;E))$$
is obtained as follows:
let $R_{\gamma^{-1}}$ be the right multiplication by $\gamma^{-1}$, and let
$E'\equal (R_{\gamma^{-1}})^*(E_{|\Gamma_y})$, i.e. $E'_g\equal E_{g\gamma^{-1}}\cong
E_g\otimes E_{\gamma^{-1}}$. Then there is an isomorphism
from ${\mathcal{K}}(L^2(\Gamma_x;E))$ to ${\mathcal{K}}(L^2(\Gamma_x;E'))$
given by 
$T\mapsto T\otimes 1$.
  However, $L^2(\Gamma_x;E')$
 and $L^2(\Gamma_y;E)$ are isomorphic under the map
 $\xi\mapsto \eta$, where $\eta(g)\equal \xi(g\gamma)$.
See \cite[pp. 76-77]{muh01} for further details.


\end{document}